\documentclass[a4paper,11pt]{article}
\usepackage[utf8]{inputenc}

\usepackage[a4paper,margin=2.3cm]{geometry}
\usepackage{amsmath, amssymb, amstext, amsfonts, amsthm, bbm, array, enumerate,scalerel}
\usepackage{mathrsfs, dsfont,mathabx}
\usepackage{nicefrac}
\usepackage[dvipsnames]{xcolor}
\usepackage{tikz}
\usepackage{stmaryrd}
\usepackage{tikz-3dplot}
\usepackage{float}
\usepackage{appendix}
\usepackage{longtable}
\usepackage[round]{natbib}
\usetikzlibrary{patterns}
\usetikzlibrary{plotmarks}
\usepackage{shuffle}
\usepackage{ stmaryrd }
 \usepackage{bbm}
 \usepackage{bm}
\usepackage{xcolor}
\usepackage{graphicx}
\usepackage{comment}
\usepackage{framed}

\usepackage{url}

\usepackage[shortlabels]{enumitem}

\newtheoremstyle{boldremark}
    {\dimexpr\topsep/2\relax}
    {\dimexpr\topsep/2\relax}
    {}
    {}
    {\bfseries}
    {.}
    {.5em}
    {}

\RequirePackage[colorlinks,citecolor=blue, urlcolor=blue, linkcolor=blue ]{hyperref}

\theoremstyle{plain}
\newtheorem{theorem}{Theorem}[section]
\newtheorem{lemma}[theorem]{Lemma}
\newtheorem{corollary}[theorem]{Corollary}
\newtheorem{proposition}[theorem]{Proposition}

\theoremstyle{definition}
\newtheorem{definition}[theorem]{Definition}
\newtheorem{example}[theorem]{Example}

\newtheorem*{assumption*}{\assumptionnumber}
\providecommand{\assumptionnumber}{}
\makeatletter

\theoremstyle{boldremark}
\newtheorem{remark}[theorem]{Remark}

\numberwithin{equation}{section}

\def\SS{\mathcal{S}^*}
\def\S{\mathcal{S}}

\def\Q{\mathbf{Q}}
\def\k{\mathbf{k}}
\def\a{\mathbf{a}}
\def\y{\mathbf{y}}
\def\u{\mathbf{u}}
\def\g{\mathbf{g}}
\def\h{\mathbf{h}}
\def\j{\mathbf{j}}
\def\v{\mathbf{v}}
\def\x{\mathbf{x}}
\def\z{\mathbf{z}}
\def\b{\mathbf{b}}
\def\c{\mathbf{c}}

\def\X{\mathbb{X}}
\def\Z{\mathbb{Z}}

\def\N{\mathbb{N}}
\def\Nset{\mathbb{N}}

\newcommand{\cadlag}{c\`adl\`ag~}

\newcommand{\Vcal}{\mathcal V}
\newcommand{\Pol}{\textup{Pol}}

\def\Rset{\mathbb R}
\newcommand{\R}{\mathbb R}
\newcommand{\Y}{\mathbb Y}
\newcommand{\E}{\mathbb E}

\newcommand{\e}{\epsilon}

\title{Entire sig-jump-diffusions, their expected signature and applications to
score-based
generative  modeling}

\author{Christa Cuchiero \thanks{University of Vienna, Department of Statistics and Operations Research, Data Science @ Uni Vienna, Kolingasse 14-16, 1090 Vienna, Austria, christa.cuchiero@univie.ac.at
} 
\and Francesca Primavera  \thanks{ École Polytechnique, CMAP \& Engie Global Markets, France.\\ francesca.primavera@polytechnique.edu
}
\and Sara Svaluto-Ferro \thanks{University of Verona, Department of Economics, Via Cantarane 24, 37129 Verona, Italy, sara.svalutoferro@univr.it \newline
This research was funded in whole or in part by the Austrian Science Fund (FWF)
Y 1235.\\
The second author is grateful for the financial support provided by Engie Global Markets.
}}
\date{}

\usepackage[skip=3pt]{parskip}
\begin{document}
\maketitle

\begin{abstract}

We develop score-based generative models for time series and standard finite-dimensional data, building on the analytical computation of the expected signature of a new class of jump-diffusions, called \textit{entire sig-jump-diffusions}. These are semimartingales whose characteristics are given by linear combinations of potentially infinitely many components of their own signature. We show that their expected signature can be computed by solving linear ODEs in the extended tensor algebra, providing a tractable foundation for learning procedures. Specifically, the learning of a generalized
score function for jump diffusions,
a key ingredient of our generative modeling
framework,  can be reduced to solving linear ODEs and a deterministic optimization problem.

We further prove that polynomial jump-diffusions can be embedded into this framework, and that their expected signature reduces to a finite-dimensional system of linear ODEs.
Building on this result, we develop an alternative approach  for finite-dimensional data, based on 
\textit{inverse signature matching}, a method that
learns the score by matching the expected signature of the reverse process with that of the inverse forward signature.

\end{abstract}

\noindent\textbf{Keywords:}  signature of \cadlag semimartingales, expected signature, holomorphic and polynomial processes, jump-diffusions, score-based generative modeling,  time series\\
\noindent \textbf{MSC (2020) Classification:} 60L10, 60J76, 68T05
\tableofcontents

\section{Introduction}\label{sec:intro}

Score-based generative models have become a leading paradigm for sampling from
complex distributions \cite{sohl2015deep, songermon2019, ho2020denoising,
song2021score, Tang2025Scorebased}. Extending them to
time series data is challenging, as sequential data are path-valued.
Indeed, applying  standard score-based generative modeling via stochastic
differential equations to time series generation requires considering
time-reversed (jump) diffusions on path space, which are inherently
infinite-dimensional and potentially intractable.
A natural way to overcome this complexity is to
handle paths  through their \emph{signature}, a universal feature map given by
iterated integrals \cite{L:98, CK:25, CC:24}, 
which -- in its logarithmic form -- has recently been applied by \cite{salvi} in the context of 
generative time series modeling. In this
paper, we construct score-based generative models for time series as well as
finite-dimensional data, building on an explicit representation of the expected
signature of a broad class of jump-diffusions. Our main contributions are the
following.

\begin{itemize}
  \item \textbf{Entire sig-jump-diffusions.} We introduce the class of
  \textit{entire sig-jump-diffusions}, whose expected signature admits an
  explicit representation via linear ordinary differential equations
  (Definition~\ref{def: sig-jd}, Theorem~\ref{th:expected sig}).
  \item \textbf{Score-based generative modeling with (lifted) entire sig-jump-diffusions  for time series.} We
  exploit this representation to construct generative models for time series, in
  which a generalized score function is learned through a deterministic
  optimization problem
  (Sections~\ref{sec:marginal-reverse-jump}--\ref{sec:invariant}, Theorem \ref{th: score opti}).
  \item \textbf{Finite-dimensional data generation and inverse signature
  matching.} For data in $\mathbb{R}^d$, we adapt this framework to the subclass of polynomial
  processes introduced in \cite{CKT:12, FL:20} (Remark~\ref{rem:appRd}). Additionally we propose an alternative to  score matching,
  called \textit{inverse signature matching} (Section~\ref{sec:alternative}).
\end{itemize}

Our starting point is the expected signature, the path-space analogue of the
sequence of moments of a random variable, which under suitable conditions
characterizes the law of a stochastic process, 
\cite{ChevyrevLyons2016, signaturemoments}. It is therefore a natural tool
whenever laws on path space must be estimated or compared, as in parameter
estimation \cite{PapavasiliouLadroue2011, semnani2025path}, hypothesis testing
\cite{gasteratos2025novelty} and generative modeling of time series
\cite{NSWLB:20, lou2023pcfgan}. It often also reduces stochastic control problems 
\cite{KLP:20, BHRS:21, cuchiero2025signature, GMP:26, abi2026stochastic} and financial derivative pricing
\cite{LyonsNejadPerezArribas2020, CGS:22, AbiJaberHainautMotte2025} to simple deterministic solutions,
since expected costs and prices become linear functionals of the expected
signature.

For these reasons, the computation of expected signatures has received
considerable attention. Explicit or recursive representations are available for
Brownian motion \cite{Fawcett2003, LyonsVictoir2004, LyonsNi2015,
boedihardjo2021expected}, Gaussian processes
\cite{BoedihardjoPapavasiliouQian2015, CassFerrucci2024}, L\'evy processes
\cite{FrizShekhar2017, CPS:UATlevy}, and more general diffusion and semimartingale
models \cite{Ni2012, friz2024expected, CST:23}.

\paragraph{Entire sig-jump-diffusions.}
We contribute to this line of work by introducing a general class of
jump-diffusions whose expected signature is obtained by solving a linear ODE in
the extended tensor algebra. To this end, we define \textit{entire
sig-jump-diffusions} as semimartingales whose characteristics can be written as a
linear combination of \textit{potentially infinitely many components} of their
own signature (Definition~\ref{def: sig-jd}). This framework is very general: by
the universal approximation properties of the signature, continuous path
functionals can be approximated arbitrarily well, on suitable compact sets and
with respect to suitable variation metrics, by linear functionals of the
signature \cite{KLP:20, BHRS:21, CPS:UATlevy},  therefore approximating (in suitable metrics) any semimartingale with sufficiently regular characteristics. The class extends 
signature SDEs introduced in \cite{CST:23}, where only semimartingales with
continuous trajectories are considered. In particular, it includes the 
signature models for asset pricing of \cite{ASS:20, CGS:22, CGMS:23, CPS:UATlevy}.

Based on duality arguments for stochastic processes \cite{JK:14, EK}, as in
\cite{CST:23, CPS:24holo}, we show that expectations of potentially infinitely
many linear combinations of the signature components of such processes can be
computed by solving linear ODEs in the extended tensor algebra
(Theorem~\ref{th:expected sig}), extending to the non-continuous setting the
results of \cite{CST:23}. The main technical challenge is the convergence of the
resulting infinite series, which is considerably harder in the presence of jumps,
as it requires control of the series and all its derivatives (equation \eqref{eq:Lu}).

For certain important cases, we can prove that these convergence conditions are
satisfied. This applies to holomorphic jump-diffusions as of \cite{CPS:24holo}
(Subsection~\ref{sec: holo}), to polynomial processes introduced in \cite{CKT:12, FL:20}
(Subsection~\ref{sec: poly}), and to some jump-diffusions with possibly
path-dependent characteristics (Remark~\ref{rem: beyondpoly}), all of which embed
in the entire sig-jump-diffusion framework. In particular, the truncated
signature of a polynomial jump-diffusion is again a finite-dimensional polynomial
jump-diffusion (Proposition~\ref{pr:expectsigpoly}), so that its expected
signature reduces to an easily computable matrix exponential (Example~\ref{ex:OU-jumps}). Another tractable subclass, whose truncated signature lift
is a finite-dimensional Markov process, is that of entire sig-jump-diffusions
with truncated signature-dependent coefficients (Subsection~\ref{sec:holo}).

\paragraph{\textbf{Score-based generative modeling with (lifted) entire sig-jump-diffusions  for time series.} }
We apply the framework of entire sig-jump-diffusions, specifically the subclass
of Subsection~\ref{sec:holo}, to score-based generative modeling of time series.
We encode the distribution of the observed time series through the distribution
of their truncated signatures and perturb it with a lifted entire
sig-jump-diffusion converging to an invariant measure independent of the initial
data distribution. Starting from this invariant distribution, a suitable
reverse-time process generates new samples approximating the original signature
distribution.

The key ingredient in this construction is a generalized score function, which
extends the classical score function to jump-diffusions. In the L\'evy setting it
corresponds to the score considered in \cite{B:26} for finite jump activity and in
\cite{yoon2023score} for isotropic $\alpha$-stable noise. The latter has been
taken up by \cite{popov2025improved}, who use a parametrized family of
reverse-time processes sharing the prescribed marginal distributions (and not the
whole law on path space), a principle that we extend beyond the $\alpha$-stable
case to general jump-diffusions. The crucial difference to these previous works is
that our definition of the generalized score function is formulated directly in
terms of the generator and the marginal probability measures, without requiring a
particular structure of the jump mechanism. As shown in
Sections~\ref{sec:marginal-reverse-jump} and \ref{sec:liftedsigjump}, this makes
it applicable to state-dependent jump-diffusions and, in particular, to the
lifted truncated signature dynamics of Subsection~\ref{sec:holo}.

Our main contribution is to exploit the analytic structure of entire
sig-jump-diffusions and their explicit expected signature representation to
reduce the learning of this score to linear ODEs and a deterministic
optimization problem (Section~\ref{sec:learning}, see in particular Theorem \ref{th: score opti}). As already mentioned, our approach is related to the
recent work of \cite{salvi}, who introduced score-based generative models
operating on log-signature representations of time series. While their approach
relies on diffusion processes in the vector space of log-signatures, we consider
lifted entire sig-jump-diffusions evolving directly on the space of signatures,
which provides a natural connection to the stochastic dynamics of the underlying
time series and allows the coefficients of the forward process to be calibrated
to the observed dynamics. We refer to Section~\ref{sec:score} for a more detailed explanation. Beyond this, our approach
differs from \cite{salvi} in two further respects. First, we incorporate jumps
into the generative dynamics. Second, thanks to the analytic tractability of
entire sig-jump-diffusions, the score is learned through the explicit expected
signature instead of neural network training.

To show that suitable forward processes exist, in Section~\ref{sec:invariant} we
construct an entire sig-diffusion whose truncated signature lift converges to a
unique invariant measure, independently of the initial data distribution, and
thus provides a natural forward perturbation process.

Since the framework operates on signatures rather than on time series, recovering
a time series from a generated signature requires a separate reconstruction step,
for which several inversion procedures exist
\cite{LyonsXu2017, LyonsXu2018, FermanianChangLyonsBiau2023, salvi}. As truncated
signatures do not determine paths uniquely, this step is subject to approximation
issues, which we do not address here.

\paragraph{Finite-dimensional data generation and inverse signature matching.}
The entire sig-jump-diffusion framework can also be used to generate
finite-dimensional data, i.e., new samples of random variables with values in
$\mathbb{R}^d$. To this end, we use polynomial processes \cite{CKT:12, FL:20},
analyzed in Subsection~\ref{sec: poly}, for the forward dynamics. In this
setting, learning the generalized score function reduces to computing moments and
solving, again, a deterministic optimization problem (Remark~\ref{rem:appRd}).

We also propose an alternative to standard score matching, which we call
\textit{inverse signature matching}, tailored to polynomial diffusions. Since,
under suitable conditions, the expected signature characterizes the law of a
stochastic process  \cite{signaturemoments}, the key idea is to learn the
generalized score by matching the expected signature of the reverse process to
that of the inverse forward signature. Under a suitable parametrization of the
score map, the polynomial structure ensures that both expectations can be
computed analytically, which reduces the problem again to a deterministic optimization
problem (Section~\ref{sec:alternative}).

Finally, we stress that 
the goal of the paper is to
establish the framework and the analytic learning procedure,
while a numerical investigation of the generative models, on time series as well
as on finite-dimensional data, is left to future work.

\paragraph{Related work on generative modeling.}
As the literature on expected signatures was already discussed above,  we here only focus on
generative modeling. In score-based generative models, the reverse-time drift
involves the score of the marginal densities  \cite{HS:86}, which is usually learned by
(denoising) score matching \cite{hyvarinen2005estimation, vincent2011connection}.
Besides the reverse SDE, the probability-flow ODE \cite{song2021score} gives a
deterministic sampler with the same marginal laws. Jump noise has been considered
in \cite{yoon2023score, B:26}, with a family of marginal-preserving reverse
dynamics in \cite{popov2025improved}; we compare with these works in the
description of our contributions above. Other works build on pretrained diffusion
models, for instance for conditional generation and fine-tuning
\cite{han2025stochastic, guo2026conditionaldiffusionguidancehard} or distillation
\cite{sigdegdistillationmakingdiffusion}; these are complementary to the
construction of the forward and reverse dynamics studied here.

For time series, several generative models are based on signatures
\cite{NSWLB:20, lou2023pcfgan}, the closest to our work being \cite{salvi}, which we already
 discussed in detail above. Other approaches that do not rely on signatures include the
Schr\"odinger bridge model of \cite{hamdouche2026nonparametric}, whose
path-dependent drift is estimated nonparametrically, its extension to
jump-diffusions \cite{demarco2026schrodinger}, and the sequential forward-backward
diffusion of \cite{cao2026adaptive}, which generates adapted series by
conditioning on the generated history and comes with statistical guarantees.
\paragraph{Organization of the paper.} In Section~\ref{sec:prelimi}, we introduce the algebraic setting of signature theory and the notion of signature for c\`adl\`ag semimartingales. In Section~\ref{sec:sigjumpdiff}, we introduce the class of entire sig-jump-diffusions, present the main representation results (Theorem~\ref{th:expected sig}), and show that holomorphic jump-diffusions, entire sig-jump-diffusions with truncated signature-dependent coefficients, and polynomial processes can be embedded within this framework. Section~\ref{sec:score} is dedicated to the application to score-based generative modeling.

\section{Preliminaries}\label{sec:prelimi}
\subsection{Notation}

Fix $d\in \mathbb N$ and let $\Rset^d$  be the Euclidean space. The extended tensor algebra over $\Rset^d$ is defined by
	\[ T((\Rset^d)):=\prod_{n=0}^\infty(\Rset^d)^{\otimes n},
	\]
	where $(\Rset^d)^{\otimes n}$ denotes the $n$-fold tensor product of $\Rset^d$ with the convention $(\Rset^d)^{\otimes 0} := \Rset$. We equip
	$T((\Rset^d))$ with the standard addition $+$, tensor multiplication $\otimes$ and scalar multiplication.
	For $N\in\Nset$, the truncated tensor algebra is defined by
	\[
	T^N(\Rset^d):=\bigoplus_{n=0}^N(\Rset^d)^{\otimes n}.
	\]

Let $I = (i_1,\dots, i_n)$ be a multi-index with entries in $\{1,\dots,d\}$. Denoting by $\epsilon_1,\ldots,\epsilon_d$ the canonical basis of $\mathbb R^d$, we use the notation $|I|:=n$ and $\epsilon_I:=\epsilon_{i_1}\otimes \epsilon_{i_2}\otimes\dots\otimes  \epsilon_{i_n}$. Observe that $(\epsilon_I)_I$ is the canonical orthonormal basis of  $(\Rset^d)^{\otimes n}$.
	Furthermore, we denote by $\e_\emptyset$ the basis element of $(\R^d)^{\otimes 0}$ and set $|\emptyset|:=0$.  We also set $I':=(i_1,\dots, i_{n-1})$ for $n>1$, ${I'}=\emptyset$ for $n=1$,   $I^{''}:=(I')'$ for $n>1$, and use the convention $\e_{I'}=0$ for $I=\emptyset$.  For  $\u\in T((\Rset^d)),$ we denote its projection on $T^N(\Rset^{d})$ by $\u^{\leq N}$ and let $\u_I$ denote its component corresponding to $\e_I$.

For two multi-indices $I\in \{1,\dots,d\}^{|I|}$, $J\in \{1,\dots,d\}^{|J|}$, and two indices of length $1$, $a,b\in\{1,\dots,d\}$, the shuffle product $\shuffle$ is defined recursively by
\begin{align*}
	&\e_I\shuffle \e_\emptyset=\e_\emptyset\shuffle \e_I=\e_I,\\
	&\e_{Ia}\shuffle \e_{Jb}=(\e_I\shuffle \e_{ Jb})\otimes \e_a+(\e_{Ia} \shuffle \e_J)\otimes \e_b,
\end{align*}
where $Ia$ denotes the concatenation of multi-indices. 

For $\u\in T(\R^d)$ and $\x\in T((\R^d))$, we write $$\langle\u,\x\rangle:=\sum_{n=0}^\infty\sum_{|I|=n}\u_I\x_I.$$  
Define  $G^N(\R^d):=\{\u^{\leq N}\in T^N(\R^d)\colon \u\in G((\R^d))\}$ for
$$G((\R^d)):=\{\x\in T((\R^d))\colon \langle \e_\emptyset,\x\rangle=1,\ \langle \e_I\shuffle \e_J,\x\rangle=\langle \e_I,\x\rangle\langle  \e_J,\x\rangle,\ \text{for all } I,J\},$$
and for $S\subseteq\Rset^d$, set 
\begin{align}\label{eq: S signature}
    \mathcal{S}:=\Big\{\x\in  G((\Rset^d))  :  \sum_{|I|=1}\x_I\e_I \in S \Big\}.
\end{align}
The set $\S$ is endowed with the sigma-algebra induced by  the product sigma-algebra on $T((\R^d))$. 

Define
\[
\mathcal{S}^*:=\{\u\in T((\Rset^d))+iT((\Rset^d)) \ : \ |\u|_{\x}<\infty \text{ for all } \ \x\in \mathcal{S}\},
\]
where $|\u|_\x:=\sum_{n=0}^\infty|\sum_{|I|=n}\u_I\x_I|$. For $\u\in \mathcal{S}^*$, we again write $$\langle\u,\x\rangle:=\sum_{n=0}^\infty\sum_{|I|=n}\u_I\x_I.$$ Moreover, we set $\u^{(0)}:=\u$ and for $k\in \N$,  $|I|=k$, we define $\u_I^{(k)}\in T((\Rset^d))+iT((\Rset^d))$ via 
\begin{align*}
 \u^{(k)}_I:=\sum_{|J|\geq0}\u_{JI}\e_J,
\end{align*}
where $JI$ denotes the concatenation of the multi indices $J$ and $I$. 
Define $\Vcal:=T((\Rset^d))+iT((\Rset^d))$ and note that 
$\u^{(k)}\in (\Vcal^d)^{\otimes k}$.
 Hence, setting $(\Vcal^d)^{\otimes 0}:=\Vcal$, 
$$(\u^{(0)},\u^{(1)},\u^{(2)},\ldots)\in T((\Vcal^d)).$$

    \subsection{Signature of \cadlag semimartingales}
    \begin{definition}
    For a semimartingale $X=(X_t)_{t\in [0,T]}$, we call the $G((\R^d))$-valued process $\X=(\X_t)_{t\in [0,T]}$
    the \emph{signature of $X$} if $\langle \e_\emptyset,\X_t\rangle=1$ and
    \begin{align}\label{eq: deI}
	\langle\epsilon_{I},\X_t\rangle =&\int_{0}^t \langle\epsilon_{I'},\mathbb{X}_{s^-}\rangle dX_s^{i_{|I|}}1_{\{|I|>0\}}
	+\frac{1}{2}\int_{0}^t\langle\epsilon_{I''},\mathbb{X}_{s^-}\rangle d[X^{i_{|I|-1}},X^{i_{|I|}}]^c_s 1_{\{|I|>1\}}\\
	&\qquad+\sum_{0<s\leq t}\sum_{\e_{I_1}\otimes \e_{I_2}=\e_I}\frac{1}{|I_2|!}1_{\{|I_2|> 1\}}\langle\epsilon_{I_1},\mathbb{X}_{s^-}\rangle \langle \e_{I_2},(\Delta X_s)^{\otimes |I_2|}\rangle,\nonumber 
\end{align}
for each $|I|>0$.
\end{definition}
Since $\mathbb{X}$ takes values in the infinite-dimensional space $G((\mathbb{R}^d))$, it is often convenient - both for practical and technical purposes - to consider a finite-dimensional truncation thereof. Recall that  for  $\u\in T((\Rset^d)),$ we denote its projection on $T^N(\Rset^{d})$, $N\in\mathbb{N}$, by $\u^{\leq N}$.

\begin{definition}
The truncated signature of order $N$ of a semimartingale $X=(X_t)_{t\in [0,T]}$ is the $T^N(\mathbb{R}^d)$-valued process $\mathbb{X}^{\leq N} = (\mathbb{X}_t^{\leq N})_{t\in[0,T]}$ defined by $\langle \epsilon_I, \mathbb{X}_t^{\leq N}\rangle := \langle \epsilon_I, \mathbb{X}_t \rangle$ for every multi-index $I$ with $|I| \leq N$, and $\langle \epsilon_I, \mathbb{X}_t^{\leq N}\rangle := 0$ otherwise.
\end{definition}

\begin{remark}
Since the right-hand side of \eqref{eq: deI} for a given multi-index $I$ only involves coefficients $\langle \epsilon_{I'}, \mathbb{X}_{s^-}\rangle$ indexed by $I'$ with $|I'| < |I|$, the process $\mathbb{X}^{\leq N}$ is itself characterized by the recursion \eqref{eq: deI} restricted to indices $|I| \leq N$, and in particular does not depend on the tensor levels of order greater than $N$. 
\end{remark}
\begin{remark}\label{remMarkus}
This concept of signature corresponds to the solution of the Marcus SDE $\X_t=\X_t\diamond dX_t$, \cite{KP:95, CPS:UATlevy}.
    This in particular implies that if $X-X_0$ is an $S$-valued semimartingale, for some $S\subseteq \Rset^d$, then $\X$ takes values in $\S$ as in \eqref{eq: S signature}, and thus satisfies the shuffle property. 
Moreover, observe that
$\Delta\X_t=\delta(\X_{t-},\Delta X_t)$ for
$$\delta(\x,\xi)_I:=\sum_{\e_{I_1}\otimes \e_{I_2}=\e_I}\frac {1}{|I_2|!}\langle \e_{I_1},\x\rangle\langle \e_{I_2},\xi^{\otimes |I_2|}\rangle1_{\{|I_2|>0\}}.$$

This dual perspective leads to different ways of simulating trajectories of $\X$.

    Fix a partition $\pi=\{0=T_0<\cdots<T_n=T\}$. Consider a trajectory $(X_t)_{t\in [0,T]}$ and fix $(\tau_i)_{i=1}^\infty$ and $(Y^i)_{i=1}^\infty$ such that
        $\sum_{s\leq t} \Delta X_s=\sum_{i=1}^\infty Y^i 1_{\tau_i\leq t}$. Set $\pi^\tau=\{0=t_0<\cdots<t_{n_\tau}=T\}:=\pi\cup\{\tau_i\colon \tau_i\leq T\}$.
    \begin{enumerate}
        \item  From $\X^{\pi^\tau}_{t_0}=\epsilon_\emptyset$, one may use the coordinatewise Euler scheme \[ \begin{aligned} \langle \epsilon_I,\X^{\pi^\tau}_{t_{k+1}}\rangle =&\langle \epsilon_I,\X^{\pi^\tau}_{t_k}\rangle +\langle \epsilon_{I'},\X^{\pi^\tau}_{t_k}\rangle (X^{i_{|I|}}_{t_{k+1}}-X^{i_{|I|}}_{t_{k}})\\
     &+\frac12\langle \epsilon_{I''},\X^{\pi^\tau}_{t_k}\rangle ([X^{i_{|I|-1}},X^{i_{|I|}}]^c_{t_{k+1}}-[X^{i_{|I|-1}},X^{i_{|I|}}]^c_{t_k})
        +\sum_{i=1}^\infty\langle \epsilon_I,\mathbb Y^i\rangle 1_{t_{k+1}=\tau_i}, \end{aligned} \]
        where the contribution from the jumps is given by
        $$ \langle \epsilon _I,\mathbb Y^i\rangle :=\sum_{\e_{I_1}\otimes \e_{I_2}=\e_I}
        \frac{1}{|I_2|!} 1_{\{|I_2|> 1\}}(\langle\epsilon_{I_1},\X^{\pi^\tau}_{\tau_i}\rangle -\langle\epsilon_{I_1},\mathbb Y^i\rangle)\langle\epsilon_{I_2},(Y^i)^{\otimes |I_2|}\rangle.$$ 
\item Alternatively, by Examples 4.20 and 4.21 in \cite{CF:19} one can compute the signature of $X$ by any  signature package based on a piecewise linear approximation of the time series. Given a sequence of deterministic nested partitions, the cited examples establish convergence in probability of the approximation at every time point contained in the union of the partitions.
        \item 
Finally, one can replace each jump of $X$ by a linear segment joining
$X_{\tau_i^-}$ to $X_{\tau_i}$, thereby obtaining a continuous path
whose signature can be computed using a standard signature
package. Since such packages do not typically require a uniform time grid,
this can be achieved by considering the time series
$(X_{ \overline t_k})_{k=0}^{\overline n_\tau}$, where 
$$\overline t_k\in\pi^\tau\cup\{\tau_i-\colon \tau_i\leq T\}$$
and $X_{ \tau_i-}:=X_{\tau_i}-Y^i$. 
If the semimartingale admits deterministic jump times such that the sequence of partitions can be chosen deterministically, then the 
 Wong-Zakai theorem  implies  uniform convergence 
 in probability 
 of the approximated signature to the true one.
    \end{enumerate}
\end{remark}

\section{Entire sig-jump-diffusions}\label{sec:sigjumpdiff}

\subsection{Characteristics of entire sig-jump-diffusions }

We introduce now the set of suitable coefficients of the processes which will be studied in this paper. Fix $S\subseteq\Rset^d$ and recall the set $\S$ introduced in equation \eqref{eq: S signature}.

\begin{definition}\label{defsigcoeff}
We call \emph{sig-coefficients} any  triplet $(b,a,K)$ such that
\begin{align*}
    b:\S\rightarrow\Rset^d, \quad a:\S\rightarrow S^{d}_+,
\end{align*}
and $K$ is a kernel from $\S$ to $\Rset^d$ such that the following conditions are verified. 

For all $\x\in\S$, 
 \begin{itemize}
	\item  $b(\x)=\langle \b,\x\rangle,$ for some $\b\in {\SS}^d$;
 \item $a(\x)+\int_{\Rset^d}\xi^{\otimes 2}   \ K(\x,d\xi)=\langle \a+\k^2,\x\rangle,$ for some $\a+\k^2\in ({\SS}^d)^{\otimes 2}$;
	\item $K(\x,\{0\})=0$ and 

  $$
			\int_{\Rset^d}\xi^{\otimes n} \ K(\x,d\xi)= \langle \k^n,\x\rangle\in (\mathbb{R}^d)^{\otimes n},
		$$
 for  some $\k^n\in ({\SS}^d)^{\otimes n }$,  $n\geq3$.
  \end{itemize}

\end{definition}

  \begin{remark}\phantomsection\label{remark: moments kernel}
\begin{enumerate}
    \item Observe that for every multiindex $I=(i_1,\dots,i_n)$, $(\k^n)^{I}\in \SS$  and thus for all $\x\in \S$, it holds $\langle (\k^n)^{I},\x\rangle\in \mathbb{R}$.
A similar reasoning applies to $\b$ and $\a+\k^2$.
\item  By definition of sig-coefficients, for all $\x\in \S$, $n\geq 3$,
$
\int_{\R^d}\|\xi\|^nK(\x,d\xi)<\infty.$
\end{enumerate}

\end{remark}

\begin{definition}\label{def: sig-jd}
    Fix some sig-coefficients $(b,a,K)$ with $K$ such that $\int_{\R^d}(\|\xi\|\wedge\|\xi\|^2)K(\x,d\xi)<\infty $, for all $\x\in \S$.
    Let $X=(X_t)_{t\in [0,T]}$ be a $(S+x_0)$-valued semimartingale with $X_0=x_0\in \mathbb{R}^d$ on some filtered probability space $(\Omega, \mathcal{F},(\mathcal{F}_t)_t,\mathbb{P})$ and let $\X:=(\X_t)_{t\in [0,T]}$ denote its signature. Assume that 
    the maps $s\mapsto b(\X_{s^-}), a(\X_{s^-}), \int_{\R^d}g(\xi)K(\X_{s^-},d\xi)$ are a.s. left-continuous for 
    $g=1_A$ and $A$ being a Borel set in $\R^d$, and $g(\xi)=\xi^{\otimes n}$ for $n\geq 2$.  
    
    We call $X$ an \emph{entire sig-jump-diffusion} with sig-coefficients $(b,a,K)$ if for all bounded $f\in C^2(\Rset^d)$,  the process $N^f:=(N^f_t)_{t\in[0,T]}$ given by
\begin{align*}
  N_t^f:=  f(X_t)-f(x_0)-\int_0^t \Big(\nabla f(X_{s^-})^\top b(\X_{s^-})+\frac{1}{2}\text{Tr}(a(\X_{s^-})\nabla ^2f(X_{s^-}))\\
     +\int_{\Rset^d}f(X_{s^-}+\xi)-f(X_{s^-})-\nabla f(X_{s^-})^\top\xi  \  K(\X_{s^-},d\xi)\Big)ds,
\end{align*}
defines a local martingale. We also denote by $X^\g$ the entire sig-jump-diffusion with sig-coefficients 
$$ b^{\bm g}(\x):= b({\bm g}\otimes\x),\qquad
a^{\bm g}(\x):= a({\bm g}\otimes\x),\qquad
 K^{\bm g}(\x,d\xi):= K({\bm g}\otimes\x,d\xi).$$
\end{definition}

\begin{remark}\label{rem:characteristics}
\begin{enumerate}
    \item  Notice that by
Proposition II.2.29 and Theorem II.2.42 in \cite{JS:87}, $X$ is a special semimartingale with characteristics with respect to the truncation function $\chi(\xi)=\xi$,
\begin{align*}
    \Bigg(\int_0^\cdot b(\X_{s^-})ds,\int_0^\cdot a(\X_{s^-})ds, K(\X_{t^-},d\xi)dt\Bigg).
\end{align*}
It is important to note that $X$ is not necessarily a jump-diffusion in the sense of Definition~II.2.18 in \cite{JS:87}. This property can nevertheless be shown to hold for the signature process 
$\X$, when viewed as an infinite-dimensional process (Lemma~\ref{lemma: deI}).
    \item The class of entire sig-jump-diffusions constitutes an extension of the class of signature SDEs introduced in \cite{CST:23}, where only processes with continuous trajectories are considered. In particular, it includes all the continuous 
signature models for asset pricing that have been proposed so far \cite{ASS:20, CGS:22, CGMS:23} where only diffusive drivers and finite sums of signature components are considered, as well as the Lévy-type signature models introduced in~\cite{CPS:UATlevy}. 
\end{enumerate}
   
\end{remark}

\subsection{Expected signature of entire sig-jump-diffusions}

In this section, we exploit duality methods for stochastic processes to show that the expected (or rather expected quantities of the) signature of entire sig-jump-diffusions can
be computed analytically by solving certain linear ODEs in the extended tensor algebra.

As a first step, we study the dynamics of the signature of such processes.
\begin{lemma}\label{lemma: deI}
    Let $X=(X_t)_{t\in [0,T]}$ be an entire sig-jump-diffusion  with sig-coefficients $(b,a,K)$ and fix a multiindex $I$. Then, for all $t\in [0,T]$,
    
    $$
        \begin{aligned}
	\langle\epsilon_{I},\X_t\rangle =& \int_0^t\bigg(\langle\epsilon_{I'},\mathbb{X}_{s^-}\rangle \langle \b^{(i_{|I|})},\X_{s^-}\rangle1_{\{|I|>0\}} +
	\frac{1}{2}\langle\epsilon_{I''},\mathbb{X}_{s^-}\rangle \langle (\a+\k^2)^{(i_{|I|-1},i_{|I|})},\X_{s^-}\rangle 1_{\{|I|>1\}} \\
	&\quad +\sum_{\e_{I_1}\otimes \e_{I_2}=\e_I}\frac{1}{|I_2|!}1_{\{|I_2|> 2\}}\langle\epsilon_{I_1},\mathbb{X}_{s^-}\rangle  \int_{\Rset^d}\langle \e_{I_2},\xi^{\otimes |I_2|}\rangle K(\X_{s^-},d\xi)\bigg)ds \\ 
	&+\textup{local martingale}.
\end{aligned}
$$
\end{lemma}
\begin{proof}
First recall that the characteristics of $X$ are given by 
$$\int_0^\cdot b(\X_{s^-})ds,\ \int_0^\cdot a(\X_{s^-})ds, \ K(\X_{s^-},d\xi)ds,$$
(Remark~\ref{rem:characteristics}).
Since by definition of entire sig-jump-diffusion the map $s\mapsto\int_{\R^d}\xi^{\otimes n}K(\X_{s^-},d\xi)$ is a.s. left-continuous  for $n\geq 2$, we have that 
\[
\int_0^T\sum_{\e_{I_1}\otimes \e_{I_2}=\e_I}\frac{1}{|I_2|!}1_{\{|I_2|> 2\}}|\langle\epsilon_{I_1},\mathbb{X}_{t^-}\rangle | \int_{\Rset^d}|\langle \e_{I_2},\xi^{\otimes |I_2|}\rangle| K(\X_{t^-},d\xi)dt<\infty \quad a.s. 
\]
By Theorem II.1.8 in \cite{JS:87}, the compensator of the second line of  equation \eqref{eq: deI} is given by
\begin{align*}
    \int_0^\cdot\int_{\R^d}\sum_{\e_{I_1}\otimes \e_{I_2}=\e_I}\frac{1}{|I_2|!}1_{\{|I_2|> 1\}}\langle\epsilon_{I_1},\mathbb{X}_{s^-}\rangle  \langle \e_{I_2} ,\xi^{\otimes |I_2|}\rangle  K(\X_{s^-},d\xi)ds.
\end{align*}
The claim follows by definition of characteristics.
\end{proof}
Lemma \ref{lemma: deI} can be reformulated in a more compact way with the following notation.
\begin{definition}\label{def: Q}
    Define
\begin{align}
\Q:=\left(0,\b, \frac{1}{2}(\a+\k^2),\frac{1}{3!}\k^3,\dots,\frac{1}{n!}\k^n,\dots\right)\nonumber
\end{align}
such that for each $I$  and all $\x\in\S$, it holds $\langle \Q^I,\x\rangle=Q(\x)_I$, for $ Q:\S\rightarrow T((\R^d))$ defined via
\begin{align}\label{eq: Qx}
Q(\x):=&\left(0,b(\x),\frac{1}{2}(a(\x)+\int_{\Rset^d}\xi^{\otimes 2}K(\x,d\xi)),\dots, \frac{1}{k!}\int_{\Rset^d}\xi^{\otimes k}K(\x,d\xi),\dots\right).
\end{align}
For each $k\in \N_0$, $\u\in \SS$, set 
\begin{align*}
    \u^{(k)}\dot{\shuffle}\ \Q:=\sum_{|I|=k}\u^{(k)}_I\shuffle \Q^I.
\end{align*}
\end{definition}
\begin{remark}\label{remark: Qx}\begin{enumerate}
\item Notice that setting $(\S^d)^{\otimes 0}:=\S$, $\Q \in T(({\SS}^d))$.
\item By Lemma 4.3(ii) in \cite{CST:23}, if for $k\in \N_0$ and all $I$, $\u^{(k)}_I\in \SS$, then $\u^{(k)}\dot{\shuffle}\ \Q\in \SS$.
\end{enumerate}

\end{remark}

\begin{corollary}\label{coro: LuN}
     Let $X=(X_t)_{t\in [0,T]}$ be an entire sig-jump-diffusion  with sig-coefficients $(b,a,K)$. 
Then    the process $M:=(M_t)_{t\in [0,T]}$ defined via
 \begin{equation}\label{eq:M finite}
	M_t:=\langle\u,\X_t\rangle -\langle \u, \e_\emptyset\rangle -\int_0^t \langle  L (\u),\X_s\rangle ds, \quad t\in [0,T],
 \end{equation}
 is a local martingale with $M_0=0$,
  for each  $\u\in T(\R^d)$, where
      \begin{align}\label{eq:Lu}
    L(\u):=\sum_{k=0}^\infty \u^{(k)}\dot{\shuffle}\Q.
 \end{align}

\end{corollary}

\begin{proof}
Observe that for each $I$ by Lemma~\ref{lemma: deI} it holds
$
	L(\e_I)=\sum_{\e_{I_1}\otimes \e_{I_2}=\e_I}\e_{I_1}\shuffle \Q^{I_2}.
$
The result follows by linearity.
\end{proof}

\begin{corollary}\label{Cornew}
    Fix ${\bm g}\in \mathcal S$, some sig-coefficients $(b,a,K)$ and  let $X^\g=(X_t^\g)_{t\in [0,T]}$ be as in Definition~\ref{def: sig-jd}.
Then    the process $M^{\bm g}:=(M_t^{\bm g})_{t\in [0,T]}$ defined via
$$
	M_t^{\bm g}:=\langle\u,{\bm g}\otimes \X_t^{\bm g}\rangle -\langle \u, {\bm g}\rangle -\int_0^t \langle  L (\u),{\bm g}\otimes\X^{\bm g}_s\rangle ds, \quad t\in [0,T],
$$
 is a local martingale with $M_0^\g=0$,
  for each  $\u\in T(\R^d)$, where
$$
    L(\u):=\sum_{k=0}^\infty \u^{(k)}\dot{\shuffle}\Q.
$$
\end{corollary}

\begin{proof}
For each $\u\in T(\R^d)$ define $\u^\g\in T(\R^d)$ as
$\u^\g_I:= \langle\u,\g\otimes \e_I\rangle$,
and note that $\langle \u,\g\otimes \x\rangle=\langle\u^\g,\x\rangle$.
Noting that the operator $\Q^\g$ corresponding to the sig-coefficients
$(b^{\bm g},a^{\bm g},K^{\bm g})$ is given by
$$
\Q^\g:=\left(0,\b^\g, \frac{1}{2}(\a^\g+(\k^2)^\g),\frac{1}{3!}(\k^3)^\g,\dots,\frac{1}{n!}(\k^n)^\g,\dots\right),
$$
by Corollary~\ref{coro: LuN} we obtain that the process $(N_t^\g)_{t\in [0,T]}$ is a local martingale with $N_0^\g=0$, where
$$N_t^\g:=\langle\u^\g, \X_t^{\bm g}\rangle -\langle \u^\g, \e_\emptyset\rangle -\int_0^t \langle  L^\g (\u^\g),\X^{\bm g}_s\rangle ds,$$
for $
    L^\g(\u):=\sum_{k=0}^\infty \u^{(k)}\dot{\shuffle}\Q^\g.
$ To conclude the proof we now just need to show that $N^\g=M^\g$, which, by definition of $\u^\g$, reduces to showing that $L^\g (\u^\g)=(L (\u))^\g$. For $|I|=k$, we have
$$\langle (\u^\g)^{(k)}_I,\x\rangle= \langle \u^\g,\x\otimes \e_I\rangle
=\langle \u,\g\otimes \x\otimes\e_I\rangle=\langle \u^{(k)}_I,\g\otimes \x\rangle
=\langle (\u^{(k)}_I)^\g, \x\rangle$$
and 
$$\langle \u^\g\shuffle \v^\g,\x\rangle
=\langle \u^\g,\x\rangle\langle \v^\g,\x\rangle
=\langle\u,\g\otimes\x\rangle\langle\v,\g\otimes\x\rangle
=\langle\u\shuffle\v,\g\otimes\x\rangle
=\langle(\u\shuffle\v)^\g,\x\rangle,$$
the claim follows.
\end{proof}

\begin{remark}
Notice that in \eqref{eq:Lu} the sum is finite and thus well defined.  In particular for each $I$ it holds  $\langle L(\e_I),\x\rangle=\langle \e_I,\x\otimes Q(\x)\rangle,$ where for each $\x\in \S$, $Q(\x)$ is given in equation \eqref{eq: Qx}.
\end{remark}
Now, we aim to extend the decomposition \eqref{eq:M finite} to some process $\langle \u,\X_\cdot\rangle$, for $\u\in \SS$ with $\u\in T((\R^d))\setminus T(\R^d)+iT((\R^d))\setminus T(\R^d)$. As a first step, we extend the operator $L$ in \eqref{eq:Lu} to a subset of its pointwise closure.
\begin{definition}\label{def:Lbar}
    For each $\u\in \SS$, we write $\u\in \mathcal D(\overline L)$ and
    $$\langle \overline L(\u),\x\rangle:=\lim_{N\to\infty}\langle L(\u^{\leq N}),\x\rangle$$
    whenever the limit is well defined for each $\x\in \S$.
\end{definition}
\begin{theorem}\label{th: polynomial}
 Let $X=(X_t)_{t\in [0,T]}$ be an entire sig-jump-diffusion  with sig-coefficients $(b,a,K)$ and denote by $\X$ its signature. 
Fix $\mathcal{U}\subseteq \mathcal D(\overline L)$ and 
for each $\u\in \mathcal{U}$, define $V^1:=(V^1_t)_{t\in [0,T]}$ and $ V^2:=(V^2_t)_{t\in [0,T]}$ via
 \begin{align*}
V^1_t:=\sup_{N\in\N} |\langle \u^{\leq N},\mathbb{X}_t\rangle |, \qquad V^2_t:=\sup_{N\in\N} \ |\langle  L(\u^{\leq N}),\mathbb{X}_t\rangle|,
 \end{align*}
 for all $t\in [0,T]$. Suppose that $(\sup_{s\leq t}V^1_s)_{t\in [0,T]}$ is locally integrable, and $V^2$ and $\langle \u,\X_\cdot\rangle$  are \cadlag processes, for each $\u\in\mathcal{U}$. Then    the process $M:=(M_t)_{t\in [0,T]}$ defined via
 \begin{equation}\label{eq:M}
	M_t:=\langle\u,\X_t\rangle -\langle \u, \e_\emptyset\rangle -\int_0^t \langle \overline L (\u),\X_s\rangle ds, \quad t\in [0,T],
 \end{equation}
 is a local martingale with $M_0=0$,
  for each $\u\in \mathcal U$.
\end{theorem}

\begin{proof}
 Fix $\u \in \mathcal{U}$.  By Corollary \ref{coro: LuN}, for all $N\in \N$ the process $M^N:=(M^N_t)_{t\in [0,T]}$ given by
   $$
       M^N_t:= \langle \u^{\leq N},\mathbb{X}_t\rangle -\langle \u^{\leq N},\e_\emptyset\rangle -\int_0^t\langle L(\u^{\leq N}),\mathbb{X}_u\rangle \ du,
$$
    is a local martingale. We show that
   $
       \sup_{N\in \N}\sup_{s\leq T}|M^N_{s\wedge\tau_n}|
$ is integrable for $  \tau_n:=\tau^1_n\wedge \tau^2_n$, where 
$(\tau^1_n)_{n\in \N}$ is a localizing sequence for the locally integrable process $V^1$   and 
   $$ \tau^2_n:=\inf\{t>0 \colon V^2_t >n\}.$$
Indeed, this implies that $M^N_{\cdot\land\tau_n}$ is a true martingale for each $N\in \N$ and by the dominated convergence theorem the same holds for the limit process $\lim_{N\to \infty}M^N_{\cdot \land\tau_n}$, which is equal to  $M=(M_t)_{t\in [0,T]}$ given in \eqref{eq:M}, as $\lim_{N\rightarrow\infty}\langle \u^{\leq N},\X_{t}\rangle=\langle \u,\X_{t}\rangle$ and 
        \begin{align*}
    &\lim_{N\rightarrow\infty }\int_0^{t\wedge \tau_n}\langle L(\u^{\leq N}),\mathbb{X}_s \rangle \ ds=\int_0^{t\wedge \tau_n}\langle \overline L(\u),\mathbb{X}_s \rangle \ ds,
        \end{align*}
        a.s.    To see that  $
       \sup_{N\in \N}\sup_{s\leq T}|M^N_{s\wedge\tau_n}|
$ is integrable it suffices to note that 
 \begin{align*}
    \sup_{s\leq T}|M^N_{s\wedge\tau_n}|\leq  C+ \sup_{s\leq T\wedge\tau_n }V^1_{s}+\int_0^{T\wedge\tau_n } V^2_u  du.
     \end{align*}
\end{proof}

Finally, we present the main result of the section. 
Building on an approach used in \cite{CST:23} and \cite{CPS:24holo}, and grounded in duality methods for stochastic processes \cite{JK:14, EK}, we show that certain expected signature quantities of entire sig-jump-diffusions can be represented as solutions to linear ODEs in the extended tensor algebra.
\begin{theorem}\label{th:expected sig}
    Let $X=(X_t)_{t\in [0,T]}$ be an entire sig-jump-diffusion with sig-coefficients $(b,a,K)$, and $\mathcal{U}\subseteq \mathcal D(\overline L)$ so that all the conditions of Theorem \ref{th: polynomial} are verified. Fix $\u\in\mathcal{U}$ and assume that there exists a $\mathcal{U}$-valued solution $(\mathbf{c}(t))_{t\in[0,T]}$ of the linear ODE
    \begin{align}\label{eq: ODE cx}
        \langle \c(t),\x\rangle=\langle \u,\x\rangle+\int_0^t\langle \overline L(\c(s)),\x\rangle ds, \quad t\in [0,T],
    \end{align}
    for all $\x\in \S$. If furthermore
    \begin{align}\label{eq:martingality}
        \E[\sup_{s,t\leq T}|\langle \c(s),\X_t\rangle|]<\infty \qquad \text{and}\qquad  \E[\sup_{s,t\leq T}|\langle \overline L(\c(s)),\X_t\rangle|]<\infty,
    \end{align}
    then it holds that
    \begin{align*}
        \E[\langle \u,\X_T\rangle]=\langle \c(T),\e_\emptyset\rangle.
    \end{align*}
\end{theorem}
\begin{proof}
    The proof follows by Lemma A.1 in \cite{CST:23}.
\end{proof}
\begin{remark}
    Notice that the conditions in \eqref{eq:martingality} ensure that the process $(M_t)_{t\in[0,T]}$ introduced in equation \eqref{eq:M} is a true martingale for each coefficient $\c(s)$, $s\in [0,T]$. 
\end{remark}
Combining Corollary~\ref{Cornew} with the procedure followed in the proof of Theorem~\ref{th: polynomial} and Theorem~\ref{th:expected sig} we obtain the following result.
\begin{corollary}\label{Corsupernew}
     Fix ${\bm g}\in \mathcal S$, some sig-coefficients $(b,a,K)$ and  let $X^\g=(X_t^\g)_{t\in [0,T]}$ be an entire sig-jump-diffusion  with sig-coefficients
$$ b^{\bm g}(\x):= b({\bm g}\otimes\x),\qquad
a^{\bm g}(\x):= a({\bm g}\otimes\x),\qquad
 K^{\bm g}(\x,d\xi):= K({\bm g}\otimes\x,d\xi).$$
 Fix $\mathcal{U}\subseteq \mathcal D(\overline L)$ so that all the conditions of Theorem \ref{th: polynomial} are verified. Fix $\u\in\mathcal{U}$ and assume that there exists a $\mathcal{U}$-valued solution $(\mathbf{c}(t))_{t\in[0,T]}$ of the linear ODE
$$
        \langle \c(t),\x\rangle=\langle \u,\x\rangle+\int_0^t\langle \overline L(\c(s)),\x\rangle ds, \quad t\in [0,T],
 $$
    for all $\x\in \S$. If furthermore
$$
        \E[\sup_{s,t\leq T}|\langle \c(s),\g\otimes\X_t^\g\rangle|]<\infty \qquad \text{and}\qquad  \E[\sup_{s,t\leq T}|\langle \overline L(\c(s)),\g\otimes\X_t^\g\rangle|]<\infty,
$$
    then it holds that
$
      \E[\langle \u,\g\otimes\X_T^\g\rangle]=\langle \c(T),{\bm g}\rangle.
$
\end{corollary}

\subsection{Examples}\label{sec:examples}
\subsubsection{Holomorphic jump-diffusions}\label{sec: holo}

In this section we demonstrate that a broad class of holomorphic jump-diffusions, introduced in \cite{CPS:24holo}, can be embedded within the framework of entire sig-jump-diffusions. The defining property of this class of processes concerns the extended generator, which is required to map a subset of holomorphic functions to themselves. If this subset includes in particular the set of polynomials on $\Rset^d$, by Lemma 3.3 in \cite{CPS:24holo}, one can show that  the characteristics $(\overline b,\overline a, \overline K)$ of such processes with respect to the truncation function $\chi(\xi)=\xi$ satisfy
\begin{equation}\label{eq: charact_tholo}
\begin{aligned}
 &\overline b_i\in\mathcal{O}(\Rset^d),\quad \overline a _{i_1i_2}+\overline m_{i_1i_2}\in \mathcal{O}(\Rset^d),\quad\text{and}\quad \overline m_{i_1\cdots i_n}\in\mathcal{O}(\Rset^d)\text{ for }n\geq 3,
    \end{aligned}
    \end{equation}
where, for each $n\geq 2$,
\begin{align*}
 \overline m_{i_1\cdots i_n}   :=\int_{\Rset^d} \xi_{i_1}\cdots \xi_{i_n}\overline K(\cdot,d\xi),
\end{align*}
and $ \mathcal{O}(\Rset^d):=\Big\{f:\Rset^d\to \mathbb{C}\colon f=h|_{\Rset^d} \text{ for some entire map } h \text{ on }\mathbb{C}^d\Big\}$\footnote{For simplicity, we consider holomorphic jump-diffusions taking values in $\Rset^d$, but our reasoning extends naturally to processes with more general state spaces.}. Note that the argument of each of these maps is the process itself and not its  signature.

To simplify the notation, for a sufficiently differentiable map $f:\R^d\to \R$ we set
\begin{align}\label{eq: notationderivative}
    f^{(\emptyset)}(x):=f(x)\qquad\text{and}\qquad f^{(J)}(x):=\partial_{x_{j_1}}\cdots\partial_{x_{j_n}} f(x),
\end{align}
for $J=(j_1,\ldots, j_n)$.
\begin{lemma}\label{lemma: holocharacteristics}
   Fix $S=\Rset^d$ and let $X=(X_t)_{t\in [0,T]}$ be a holomorphic jump-diffusion such that $X-X_0$ takes values in $S$ and its characteristics $(\overline b,\overline a, \overline K)$ satisfy the conditions in equation \eqref{eq: charact_tholo}.
Then $X$ is an entire sig-jump-diffusion 
 whose sig-coefficients $(b,a,K)$ satisfy the following conditions. 
 
 For all $\x\in \S$, 

 \begin{itemize}
	\item  $b(\x)=\langle \b,\x\rangle,$ for some $\b\in{\SS}^d$ given by
 $$\b^{I}_J=\overline b_{i_1}^{(J)}(X_0)\qquad\text{for  each }  I=(i_1) \text{ and multiindex }J;$$
 \item $a(\x)+\int_{\Rset^d}\xi^{\otimes 2}   \ K(\x,d\xi)=\langle \a+\k^2,\x\rangle,$ for some $\a+\k^2\in ((\SS)^d)^{\otimes 2}$ given by
 $$(\a+\k^2)^{I}_J=(\overline a_{i_1i_2}+\overline m_{i_1i_2})^{(J)}(X_0)\qquad\text{for  each }I=(i_1,i_2) \text{ and multiindex } J;$$
	\item for all $n\geq 3$, $\int_{\Rset^d} \xi^{\otimes n} 
 K(\x,d\xi)= \langle \k^{n},\x\rangle$, for some $\k^{n}\in ((\SS)^d)^{\otimes n}$  given by
 $$(\k^n)^{I}_J=\overline m_{i_1\cdots i_{n}}^{(J)}(X_0)\qquad\text{for each }I=(i_1,\dots,i_n) \text{ and multiindex } J.$$
  \end{itemize}
  \end{lemma}
\begin{proof}
 First of all notice that the maps in $\mathcal{O}(\Rset^d)$ admit a power series expansion centered at any point in $\Rset^d$ \cite[Section D]{CPS:24holo}. 
    Then, the proof follows by performing the Taylor expansion of the maps in  \eqref{eq: charact_tholo} centered at $X_0$, and noticing that by the properties of the shuffle product \cite[Remark 2.8 (iii)]{CST:23}, for each $\x\in\S\subseteq G((\Rset^d))$ and each multiindex $|I|\geq 1$ it holds
$$
    \sum_{J\colon J^{ord}=I^{ord}}\langle \e_J,\x\rangle=\frac 1 {I!} \langle \e_{i_1},\x\rangle \cdots  \langle \e_{i_{|I|}},\x\rangle,
$$
    where $I^{ord}$ is the multi-index obtained by sorting the entries of $I$ and $I!:=\prod_{k=1}^d I(k)!$ for $I(k):=\sum_{j=1}^{|I|}1_{\{i_j=k\}}$.

\end{proof}

\begin{remark}
   
    For holomorphic jump-diffusions in the class of entire sig-jump-diffusions, Theorem \ref{th:expected sig} corresponds to Theorem 3.20 in \cite{CPS:24holo} whenever for each $\u \in \mathcal{U}$ and each multi-index $I$, it holds $\u_{I}=\u_{I^{ord}}$, for $I^{ord}$ introduced in the proof of Lemma \ref{lemma: holocharacteristics}.
    In such a case the local martingality of the process $M$ in equation \eqref{eq:M} follows directly from the definition of holomorphic jump-diffusion (Definition 3.1 in \cite{CPS:24holo}). Moreover, the existence of a solution of the sequence-valued ODE \eqref{eq: ODE cx} can be derived exploiting some tools from complex analysis  \cite[Section 3.5]{CPS:24holo}.
\end{remark}

\subsubsection{Entire sig-jump-diffusions with truncated signature-dependent coefficients}
\label{sec:holo}

Fix two  entire maps
$\overline b:T^N(\mathbb R^d)\to\mathbb R^d$ and
$\overline\sigma:T^N(\mathbb R^d)\to\mathbb R^{d\times d}$ bounded on their real domain.
Let $(E,\mathcal E,F)$ be a finite  measure space and let
$$
\overline\delta:T^N(\mathbb R^d)\times E\to\mathbb R^d
$$
be a bounded measurable map such that
$\overline\delta(\cdot,z)$ is entire for every $z\in E$, with local Lipschitz constants that can be chosen uniformly in $z$.  
The assumption $F(E)<\infty$ can be relaxed, but it simplifies the assumptions needed for the proof of Lemma~\ref{lemUniqueness} later on. Fix ${\bm g}\in G^N(\R^d)$   and set 
$$\overline b^{\bm g}(\x):=\overline b({\bm g}\otimes\x),\qquad
\overline \sigma^{\bm g}(\x):=\overline \sigma({\bm g}\otimes\x),\qquad
\overline \delta^{\bm g}(\x,z):=\overline \delta({\bm g}\otimes\x,z).$$
Let $\mu(dt,dz)$ be a Poisson random measure on
$\mathbb R_+\times E$ with compensator $dt F(dz)$, independent of $W$, and consider the SDE
\begin{equation}\label{entsig}
dX_t^{\bm g}
=
\overline b^{\bm g}(\X^{\bm g}_{t-})dt
+
\overline\sigma^{\bm g}( \X^{\bm g}_{t-})dW_t
+
\int_E
\overline\delta^{\bm g}( \X^{\bm g}_{t-},z)
\bigl(\mu(dt,dz)-F(dz)dt\bigr),
\qquad X_0^{\bm g}=x_0,
\end{equation}
where $W$ is a $d$-dimensional Brownian motion and
$\X^{\bm g}_t$ is the truncated  signature of $X^{\bm g}$. Note that we omit $\leq N$ for notational convenience.

By Remark~\ref{remMarkus} we know that
$
\Delta\X^{\bm g}_t
=
\delta\bigl(\X^{\bm g}_{t-},\Delta X_t^{\bm g}\bigr)
$ for
$$\delta(\x,\xi)_I:=\sum_{\e_{I_1}\otimes \e_{I_2}=\e_I}\frac {1}{|I_2|!}\langle \e_{I_1},\x\rangle\langle \e_{I_2},\xi^{\otimes |I_2|}\rangle1_{\{|I_2|>0\}}.$$

An application of It\^o's formula for semimartingales with jumps shows that, setting
\begin{align*}
b_I^{\bm g}(\x )
&:=
\langle \e_{I'},\x \rangle
\overline b_{i_{|I|}}^{\bm g}(\x )
+
\frac{1}{2}
\langle \e_{I''},\x \rangle
\bigl(\overline\sigma^{\bm g}(\x )\overline\sigma^{\bm g}(\x )^\top\bigr)_
{i_{|I|-1}i_{|I|}}\\
&\qquad
+\int_E
\left(
\delta(\x ,\overline\delta^{\bm g}(\x ,z))_I
-
\langle \e_{I'},\x \rangle
\overline\delta^{\bm g}(\x ,z)_{i_{|I|}}\right)F(dz),
\\
\sigma_{Ik}^{\bm g}(\x )
&:=
\langle \e_{I'},\x \rangle
\overline\sigma_{i_{|I|}k}^{\bm g}(\x ), \qquad \gamma^{\bm g}(\x ,z):=\delta(\x ,\overline\delta^{\bm g}(\x ,z))
\end{align*}
the dynamics of the truncated signature satisfy
\begin{align} \label{eq:Xg}
\begin{split}
d\X^{\bm g}_t
=&
b^{\bm g}(\X^{\bm g}_{t-})dt
+
\sigma^{\bm g}( \X^{\bm g}_{t-})dW_t
+
\int_E
\gamma^{\bm g}(\X^{\bm g}_{t-},z)
\bigl(\mu(dt,dz)-F(dz)dt\bigr),
\end{split}
\end{align}
for $\X^{\bm g}_0=\e_\emptyset$. Setting 
$b:=b^{\e_\emptyset},\sigma:=\sigma^{\e_\emptyset}$ and $\gamma:=\gamma^{\e_\emptyset}$ and noting that
$${\bm g}\otimes b^{\bm g}(\x)=b({\bm g}\otimes \x),\qquad
{\bm g}\otimes \sigma^{\bm g}(\x)=\sigma({\bm g}\otimes \x),\qquad
{\bm g}\otimes \gamma^{\bm g}(\x,z)=\gamma({\bm g}\otimes \x,z),$$
we can conclude that the process $({\bm g}\otimes \X^{\bm g}_t)_{t\in [0,T]}$ is a solution of the SDE
\begin{align}\label{entsig2}
\begin{split}
d\Z_t^{\leq N}
=&
b(\Z_{t-}^{\leq N})dt
+
\sigma( \Z_{t-}^{\leq N})dW_t
+
\int_E
\gamma(\Z_{t-}^{\leq N},z)
\bigl(\mu(dt,dz)-F(dz)dt\bigr),\qquad \Z_0^{\leq N}={\bm g}.
\end{split}
\end{align}

The finite-dimensional nature of this process gives access to classical existence and uniqueness results. Here the pair  $(X^{\g}, \mathbb{Z}^{\leq N})$ is regarded as the unknown
of the coupled system \eqref{entsig}--\eqref{entsig2}.

\begin{lemma}\label{lemUniqueness}
The system of SDEs \eqref{entsig}-\eqref{entsig2} admits a unique non-explosive strong càdlàg solution on $[0,T]$.
\end{lemma}

\begin{proof}
We prove the result for ${\bm g}=\e_\emptyset$; the general case is analogous.
Since $\overline b$ and $\overline\sigma$ are entire, and
$\overline\delta(\cdot,z)$ is entire with local Lipschitz constants that are uniform in $z$, all these maps are locally Lipschitz in their state variable. Since the coordinate maps
$
\x\mapsto\langle \e_I,\x\rangle
$
are linear, it follows that $b$ and $\sigma$ are locally Lipschitz and that, for every $R>0$, there exists a constant $L_R>0$ such that
$$
\int_E
\lVert\gamma(\x ,z)-\gamma(\y ,z)\rVert^2
F(dz)
\leq
L_R\lVert \x -\y \rVert^2,
\qquad
\lVert \x \rVert,\lVert \y \rVert\leq R.
$$

Moreover, the boundedness of $\overline b$, $\overline\sigma$, and $\overline\delta$, together with the finiteness of $F$, implies that
$$
\lVert b(\x )\rVert
+
\lVert\sigma(\x )\rVert
+
\left(
\int_E\lVert\gamma(\x ,z)\rVert^2F(dz)
\right)^{1/2}
\leq
C\bigl(1+\lVert \x \rVert\bigr),
\qquad
\x \in T^N(\mathbb R^d),
$$
for some constant $C>0$. Thus, $b$, $\sigma$, and $\gamma$ satisfy the local Lipschitz and linear-growth conditions for stochastic differential equations with jumps. By
\cite[Theorem~III.2.32]{JS:87}, the system \eqref{entsig}-\eqref{entsig2} admits a unique non-explosive strong càdlàg solution on $[0,T]$.
\end{proof}

\begin{remark} \label{rem:connection}
\begin{enumerate}
\item[]
    \item 
Note that the first level $\mathbb{Z}^{(1)}$ of \eqref{entsig2} is equal to $X^{\g}$ if and only if $\g^{(1)}=x_0$. This is just a consequence of
$\mathbb{Z}^{\leq N}={\bm g}\otimes \X^{\bm g}$, which at the first level reads as $\mathbb{Z}^{(1)}= \g^{(1)}+ X^{\g}-x_0$.

\item Observe furthermore that the
translation invariance of the signature implies that 
\[
X^{\g, y_0}=X^{\g,x_0}-x_0+y_0,
\]
where $X^{\g, x_0}$ denotes the solution of \eqref{entsig} when started at $x_0$. This in turn means that
$\mathbb{X}^{\g}$ and $\mathbb{Z}^{\leq N}$ are independent of the initial value $x_0$ in \eqref{entsig}. 
\end{enumerate}
\end{remark}

\begin{remark}\label{remcontsig}
 As an example, note that for fixed $s\in [0,T]$ we have
\begin{align*}
    dX_{s+t}
& =\overline b(\X_{(s+t)-}^{\leq N})dt
+
\overline\sigma( \X_{(s+t)-}^{\leq N})dW_{s+t}
+
\int_E
\overline\delta( \X_{(s+t)-}^{\leq N},z)
\bigl(\mu_s(dt,dz)-F(dz)dt\bigr)\\
& =\overline b(\X_{s}^{\leq N}\otimes \X_{s,s+t-}^{\leq N})dt
+
\overline\sigma( \X_{s}^{\leq N}\otimes \X_{s,s+t-}^{\leq N})d(W_{s+t}-W_s)\\
&\qquad+
\int_E
\overline\delta( \X_{s}^{\leq N}\otimes \X_{s,s+t-}^{\leq N},z)
\bigl(\mu_s(dt,dz)-F(dz)dt\bigr),
\end{align*}
showing that the solution of \eqref{entsig2} for ${\bm g}=\X_s^{\leq N}$ is given by 
$\X_{s}^{\leq N}\otimes \X_{s,s+t}^{\leq N}=\X_{s+t}^{\leq N}.$
\end{remark}

\begin{remark}\label{Rem:sigisholo}
Observe that each component of $\overline b^\g$, $\overline a^\g:=\overline\sigma^\g(\overline\sigma^\g)^\top$, and $\overline\delta^\g(\cdot,z)$ is an entire map in $\x\in T^N(\R^d)$, for each $z\in E$. Since these maps depend on finitely many signature coordinates, their power series and the shuffle identity give the required sig-coefficient representations for the drift and covariance. Write then the jump kernel as
\[
K^\g(\x,B)
:=\int_E 1_{B\setminus\{0\}}
       \bigl(\overline\delta^\g(\x,z)\bigr)F(dz).
\]
To ensure that $(\overline b^\g,\overline a^\g,K^\g)$ are sig-coefficients and thus \eqref{entsig} is an entire sig-jump-diffusion it remains to obtain sig-coefficient representations for the jump moments
\[
m_n^\g(\x):=\int_E
       \bigl(\overline\delta^\g(\x,z)\bigr)^{\otimes n}F(dz),
\qquad n\geq2.
\]
A sufficient condition is that, for every $n\geq2$ and every compact $C\subset T^N(\mathbb C^d)$,
\[
\int_E\sup_{\x\in C}
       \|\overline\delta(\x,z)\|^n\,F(dz)<\infty.
\]
It holds, in particular, if $\overline\delta$ is locally bounded on the complex domain uniformly in $z\in E$, since $F(E)<\infty$. Under this condition each $m_n^\g$ is entire in \cite[Lemma~A.1 and Remark~3.16(iv)]{CPS:24holo} and the claim follows by the shuffle identity.

Such condition holds, for example, when $\overline\delta(\x,z)=v(z)$ for some bounded map $v$. It also holds when $\overline\delta(\x,z)=z+h(\x)$, $F$ has finite moments, and $h$ is entire and bounded on its real domain.

Under the same condition on $\overline\delta$, the finite dimensional process
$\Z^{\leq N}={\bm g}\otimes\X^{\bm g}$ is also a holomorphic
jump-diffusion in the sense of Definition~3.1 in
\cite{CPS:24holo}, with the space of test functions
chosen as the polynomials on $T^N(\R^d)$. Indeed, its drift and
covariance have entire components, and the same holds for all its
jump moments: the lifted jump size $\gamma(\x,z)$ is a polynomial
in $\overline\delta(\x,z)$ of degree at most $N$, with coefficients
depending linearly on $\x$, so the assumed local integrability
condition also applies to its powers. Larger spaces of holomorphic
test functions can be considered whenever the conditions of Theorem~3.12 in
\cite{CPS:24holo} are satisfied.

In particular, provided the corresponding sequence-valued linear
ODE admits a solution satisfying the compatibility, martingale,
and integrability assumptions of Theorem~3.20 in
\cite{CPS:24holo}, the holomorphic formula established
therein represents expectations of holomorphic functions of
$\Z^{\leq N}$ in terms of that solution.
 This formula then coincides with the one provided in Corollary 
\ref{Corsupernew}, by choosing the coefficients $\u$ in line with the corresponding holomorphic function.
\end{remark}

\subsubsection{Polynomial jump-diffusions and their expected signature}\label{sec: poly}

We now focus on the class of polynomial jump-diffusions introduced in \cite{CKT:12} and further analyzed in \cite{FL:20}. These processes form a subclass of holomorphic jump-diffusions and naturally fit within the entire sig-jump-diffusion framework. For this class, we explicitly compute the sig-coefficients and show that the infinite-dimensional ODE \eqref{eq: ODE cx} reduces to a finite-dimensional system of linear ODEs. This reduction allows for efficient and explicit computation of the expected signature via matrix exponentiation.

First, recall that by Lemma 1 in \cite{FL:20} a polynomial jump-diffusion is a special jump-diffusion such that the coefficients $(\overline a,\overline b,\overline K)$ determining its characteristics with respect to the truncation function $\chi(\xi)=\xi$ satisfy
\begin{equation}\label{eq: chara poly process}
\begin{aligned}
 &\overline b_i\in \Pol_1(\R^d),\quad \overline a _{i_1i_2}+\overline m_{i_1i_2}\in \Pol_2(\R^d),\quad\text{and}\quad \overline m_{i_1\cdots i_n}\in\Pol_n(\Rset^d) \text{ for }n\geq 3,
    \end{aligned}
    \end{equation}
where, for each $n\geq 2$,
\begin{align*}
 \overline m_{i_1\cdots i_n}   :=\int_{\Rset^d} \xi_{i_1}\cdots \xi_{i_n}\overline K(\cdot,d\xi),
\end{align*}
and $\Pol_i(\Rset^d)$ denotes the set of polynomials on $\Rset^d$ up to degree $i\in \Nset$.

In the following lemma, we adopt the same notation for partial derivatives (equation \eqref{eq: notationderivative}) as in Lemma \ref{lemma: holocharacteristics}, and the proof proceeds along similar lines.
\begin{lemma}\label{lemma: polycharacteristics}
   A process $X=(X_t)_{t\in [0,T]}$ is a polynomial jump-diffusion such that $X-X_0$ takes values in $S$ if and only if $X$ is an entire sig-jump-diffusion 
 whose sig-coefficients $(b,a,K)$ satisfy the following conditions. 
 
 For all $\x\in \S$, 

 \begin{itemize}
	\item  $b(\x)=\langle \b,\x\rangle,$ for some $\b\in (T^1(\Rset^d))^d$ given by
 $$\b^{I}_J=\overline b_{i_1}^{(J)}(X_0)\qquad\text{for  each }  I=(i_1) \text{ and } |J|\leq 1;$$
 \item $a(\x)+\int_{\Rset^d}\xi^{\otimes 2}   \ K(\x,d\xi)=\langle \a+\k^2,\x\rangle,$ for some $\a+\k^2\in ((T^2(\Rset^d))^d)^{\otimes 2}$ given by
 $$(\a+\k^2)^{I}_J=(\overline a_{i_1i_2}+\overline m_{i_1i_2})^{(J)}(X_0)\qquad\text{for  each }I=(i_1,i_2) \text{ and } |J|\leq 2;$$
	\item for all $n\geq 3$, $\int_{\Rset^d} \xi^{\otimes n} 
 K(\x,d\xi)= \langle \k^{n},\x\rangle$, for some $\k^{n}\in ((T^{n}(\Rset^d))^d)^{\otimes n}$  given by
 $$(\k^n)^{I}_J=\overline m_{i_1\cdots i_{n}}^{(J)}(X_0)\qquad\text{for each }I=(i_1,\dots,i_n) \text{ and } |J|\leq n.$$
  \end{itemize}
  \end{lemma}

\begin{remark}
A polynomial jump-diffusion is fully characterized by its sig-coefficients. This is not true for the broader class of holomorphic jump-diffusions (Lemma \ref{lemma: holocharacteristics}), where expressing the characteristics in terms of sig-coefficients is not sufficient to guarantee the process is a holomorphic jump-diffusion \cite[Lemma 3.3]{CPS:24holo}.

\end{remark}

Next, we show that for a polynomial jump-diffusion and any initial condition $\u\in T(\Rset^d)$, the infinite-dimensional ODE \eqref{eq: ODE cx} reduces to a finite-dimensional system of linear ODEs allowing the explicit computation of expected signature terms via matrix exponentiation.
\begin{proposition}	\label{pr:expectsigpoly}
	Let $X=(X_t)_{t\in [0,T]}$ be a polynomial jump-diffusion with coefficients $(\overline b,\overline a,\overline K)$ such that $X-X_0$ takes values in $S$. Then, for each $N\in \mathbb N$, the truncated signature $\mathbb{X}^{\leq N}$ of $X$ is a $T^N(\Rset^d)$-valued polynomial jump-diffusion. This in particular implies that all the moments of $\mathbb{X}^{\leq N}$ exist and are finite.
    
Let $G$ denote the matrix representation, with respect to a basis $h_n$, of the restriction of the infinitesimal generator of $\mathbb{X}^{\leq N}$ to the invariant space of polynomials of degree at most $n$. The polynomial property then yields the following conditional moment formula:
    	\begin{equation}\label{eqn4}
	\mathbb{E}[\vec{p}^\top h_n(\mathbb{X}_T^{\leq N})|\mathcal F_t]=\vec{p}^\top \exp((T-t)G^\top)h_n(\mathbb{X}_t^{\leq N}),
	\end{equation}
where $\exp\left((T-t)G\right)$ denotes the matrix exponential. 
\end{proposition}

\begin{proof}  The proof of the polynomial property of $\X^{\leq N}$ is similar to the proof of Proposition~4.5 in \cite{CPS:UATlevy}.
In particular,  we get that $\X^{\leq N}$ is a jump-diffusion whose coefficients 
$$(\beta,\zeta,\mathcal{K}):=((\beta_I)_{|I|\leq N},(\zeta_{I,J})_{|I|,|J|\leq N}, \mathcal{K})$$
satisfy
\begin{align}
\beta_I(\cdot)&=\langle \e_{I'},\cdot\rangle \overline b_{i_{|I|}}(\pi_1^{X_0}(\cdot))
+\frac{1}{2}\langle\e_{I''},\cdot\rangle (\overline a_{(i_{|I|-1},i_{|I|})}+\overline m_{(i_{|I|-1},i_{|I|})})(\pi_1^{X_0}(\cdot))\notag\\
&\qquad + \sum_{\e_{I_1}\otimes \e_{I_2}=\e_I} 1_{\{|I_2|>2\}}\frac {1}{|I_2|!}\langle \e_{I_1} ,\cdot\rangle\overline m_{I_2}(\pi_1^{X_0}(\cdot)),\notag\\
\mathcal{K}(\cdot,A)&=\int_{\R^d}1_{A\setminus\{0\}}(\delta(\cdot,\xi))\overline K(\pi_1^{X_0}(\cdot),d\xi),\notag\\ 
\zeta_{I,J}(\cdot)&=\langle \epsilon _{I'},\cdot\rangle\langle \epsilon _{J'},\cdot\rangle\overline a_{(i_{|I|},j_{|J|})}(\pi_1^{X_0}(\cdot)),\notag
\end{align}
where $\pi_1^{X_0}$ satisfies $\pi_1^{X_0}(\X_t)=X_t$, and $\delta:T^N(\Rset^d)\times \Rset^d\rightarrow T^N(\Rset^d)$ is given by
\begin{align*}
\delta(\x,\xi)_I:=\sum_{\e_{I_1}\otimes \e_{I_2}=\e_I}\frac {1}{|I_2|!}\langle \e_{I_1},\x\rangle\langle \e_{I_2},\xi^{\otimes |I_2|}\rangle1_{\{|I_2|>0\}}.
\end{align*}
By Lemma 1 in \cite{FL:20} it suffices to check that
\begin{align*}
&\beta_I(\cdot)\in  \Pol_1(T^N(\R^d)),\\
    &\zeta_{IJ}(\cdot )+\int_{T^N(\R^d)}\y_I\y_J\mathcal{K}(\cdot,d\y)\in  \Pol_2(T^N(\R^d)),\\
    &\int_{T^N(\R^d)}\y_{I_1}\dots\y_{I_n}\mathcal{K}(\cdot,d\y)\in \Pol_n(T^N(\R^d)), \quad  \text{for }n\geq3,
\end{align*}
where $\Pol_i(\Rset^d)$ denotes the set of polynomials on $\Rset^d$ up to degree $i\in \Nset$.
We now show how to conclude the proof of the second condition; the other conditions can be shown analogously. Observe that
\begin{align*}
    \zeta_{I,J}(\cdot)&+\int_{T^N(\R^d)}\y_I\y_J\mathcal{K}(\cdot,d\y)
    =
    \langle \epsilon _{I'},\cdot\rangle\langle \epsilon _{J'},\cdot\rangle\overline a_{(i_{|I|},j_{|J|})}(\pi_1^{X_0}(\x))\\
    &+
    \sum_{\e_{I_1}\otimes \e_{I_2}=\e_I}\sum_{\e_{J_1}\otimes \e_{J_2}=\e_J}\frac {1}{|I_2|!|J_2|!}\langle \e_{I_1},\x^N\rangle\langle \e_{J_1},\x^N\rangle1_{\{|I_2|>0\}}1_{\{|J_2|>0\}}
\overline m_{I_2J_2}(\pi_1^{X_0}(\cdot))\\
  =
    &\langle \epsilon _{I'},\cdot\rangle\langle \epsilon _{J'},\cdot\rangle(\overline a_{(i_{|I|},j_{|J|})}+\overline m_{i_{|I|}j_{|J|}})(\pi_1^{X_0}(\x))\\
    &+
    \sum_{\e_{I_1}\otimes \e_{I_2}=\e_I}\sum_{\e_{J_1}\otimes \e_{J_2}=\e_J}\frac {1}{|I_2|!|J_2|!}\langle \e_{I_1},\x^N\rangle\langle \e_{J_1},\x^N\rangle1_{\{|I_2|+|J_2|>2\}}1_{\{|I_2|>0\}}1_{\{|J_2|>0\}}
\overline m_{I_2J_2}(\pi_1^{X_0}(\cdot))
\end{align*}
This  implies that for $|I|,|J|\leq N$, setting 
$$(\overline a_{i_{|I|},j_{|J|}}+\overline m_{i_{|I|},j_{|J|}})(\pi_1^{X_0}(\cdot))=\sum_{|H_1|,|H_2|\leq 1}\alpha_{H_1H_2}\langle \e_{H_1}, \cdot\rangle \langle \e_{H_2}, \cdot\rangle, $$ we get  
\begin{align*}
    &\langle \e_{I'},\cdot\rangle\langle  \e_{J'},\cdot\rangle  (\overline a_{i_{|I|},j_{|J|}}+ \overline m_{i_{|I|},j_{|J|}})(\pi_1^{X_0}(\cdot))\\
    &\qquad \qquad \quad  =\sum_{|H_1|,|H_2|\leq 1}\alpha_{H_1H_2}\langle \e_{I'}\shuffle \e_{H_1},\cdot\rangle\langle  \e_{J'}\shuffle \e_{H_2},\cdot\rangle \in \Pol_2(T^N(\R^d)).
\end{align*}
Similarly, exploiting the polynomial property of $\overline m_{i_1\cdots i_{n}}$, for $n\geq 3$ the claim follows. Finally, the representation in equation
\eqref{eqn4} follows then by an application of Theorem~2.5 in \cite{FL:20}.
\end{proof}

As a corollary of Proposition~\ref{pr:expectsigpoly}, the expected signature of a polynomial process, which corresponds to a polynomial of degree 1 on the truncated signature, can be computed analytically. 

\begin{corollary}\label{coro:expectedsigpoly}
Let $X=(X_t)_{t\in [0,T]}$ be a polynomial jump-diffusion with coefficients $(\overline b,\overline a,\overline K)$ such that $X-X_0$ takes values in $S$. 
In the special case of the invariant space given by linear polynomials one can choose the basis $\{\langle \e_I,\cdot\rangle\colon |I|\leq N\}$ and we get $G=L$ for $L$ as in Corollary~\ref{coro: LuN} leading to
		$$
	\mathbb{E}[\mathbb{X}_T^{\leq N}|\mathcal F_t]=
 \exp((T-t)L^\top)\X_t^{\leq N}.
$$
Note in particular that the column of $L$ corresponding to $I$ has the representation
\begin{align}\label{eq:Lpoly}
    L\e_I=\e_{I'}\shuffle \b^{(i_{|I|})} 
+\frac{1}{2}\e_{I''}\shuffle  ((\a+\k^2)^{(i_{|I|-1},i_{|I|})}) 
+ \sum_{\e_{I_1}\otimes \e_{I_2}=\e_I} 1_{\{|I_2|>2\}}\frac {1}{|I_2|!} \e_{I_1} \shuffle (\k^{|I_2|})^{I_2},
\end{align}
 for $\b,\a,\k$ computed in Lemma \ref{lemma: polycharacteristics}.
\end{corollary}

\begin{remark}
Let us remark that 
Proposition~\ref{pr:expectsigpoly} can be extended beyond the polynomial jump-diffusion setting.
To make this precise, let us introduce the following set of graded polynomials.

For $N \in \mathbb{N}$ we define the (ordered) tuples of multi-indices with individual length between $1$ and $N$ and combined length at most $N$ as 
    \begin{align*}
        \mathcal{I}_{N} = \Big\{(I_1,\dots,I_n) \, \Big|\,  n \in \N, 1 \leq |I_1|, \dots, |I_n| \leq N, \sum_{j=1}^n|I_j| \leq N\Big\}.
    \end{align*}
   Graded polynomials are then given as follows: 
    \begin{align*}
        \mathscr{P}_N(\R^d)
        &= 
        \Big\{
        T^N(\mathbb{R}^d) \ni \x  \mapsto c +\sum_{\mathcal{J} \in \mathcal{I}_{N}} c_{\mathcal{J}} \prod_{I \in \mathcal{J}} \x_{I} \; \Big| \; c,c_{\mathcal{J}} \in \mathbb{R}, \,\mathcal{J} \in \mathcal{I}_{N}
        \Big\}.
\end{align*}
For $N = 1$,
$\mathscr{P}_1(\R^d)$ coincides with $\Pol_1(\R^d)$, but for $N \geq 2$ the notions are different.
For simplicity, let $d=2$. Then, for $N = 2$, 
    $\mathscr{P}_2$ is spanned by $\{1,\x_1,\x_2, \x_{11}, \x_{12}, \x_{21}, \x_{22}, \x_1^2, \x_1\x_2,\x_2^2 \}$ which clearly differs from $\Pol_2(\R^d)$ and also  from the set of polynomials with degree at most two in  the variables $\x_1,\x_2, \x_{11}, \x_{12}, \x_{21}$ and $\x_{22}$, as for example the term $\x_{11}^2$ is missing.\\

If the characteristics of a special semimartingale $X$
 with respect to the truncation function $\chi(\xi)=\xi$ now satisfy
\begin{equation}\label{eq: chara beyond poly process}
\begin{aligned}
 &\overline b_i\in \mathscr{P}_1(\R^d)=\Pol_1(\R^d),\quad \overline a _{i_1i_2}+\overline m_{i_1i_2}\in \mathscr{P}_2(\R^d),\quad\text{and}\quad \overline m_{i_1\cdots i_n}\in\mathscr{P}_n(\Rset^d) \text{ for }n\geq 3,
    \end{aligned}
    \end{equation}
then we clearly stay in the entire sig-jump-diffusion setup, but
go beyond the polynomial jump-diffusion case.  Nevertheless, 
a conclusion similar to that of Proposition~\ref{pr:expectsigpoly} remains true.
That is, 
for each $N\in \mathbb N$, the truncated signature $\mathbb{X}^{\leq N}$ of $X$ is a $T^N(\Rset^d)$-valued \emph{graded} polynomial jump-diffusion. This means that  there is a conditional moment formula for every graded polynomial in $\mathscr{P}_N(\R^d)$ similar to \eqref{eqn4},
with $G$ being the matrix representation, with respect to a basis $h_n$, of the restriction of the infinitesimal generator of $\mathbb{X}^{\leq N}$ to the invariant space $\mathscr{P}_N(\R^d)$. The proof proceeds as in Proposition~\ref{pr:expectsigpoly}, but for 
the conclusion of Lemma D.4 in \cite{FL:20}, the true martingale property of $N^f$ (needed to apply Theorem 2.5 in \cite{FL:20}), certain estimates based on a graded norm have to be used. We refer to Lemma 3.19, Lemma 3.20 and Theorem 3.21 in \cite{CP:26}, where the details in the pure diffusion case are carried out.
\end{remark}

\begin{remark}\label{rem: beyondpoly}

For a concrete specification of characteristics of the form \eqref{eq: chara beyond poly process} in the pure diffusion case 
with the second characteristic in $\mathscr{P}_2(\R^d) \setminus \Pol_2(\R^d)$ 
we refer to Example 3.10 in \cite{CP:26}.
This pure diffusion example can be extended by adding jumps of a certain form. Indeed, consider the case of a jump kernel explicitly given via
     \begin{align*}
    K(\cdot,A)=\lambda(\cdot )\int_E 1_{A\setminus\{0\}} ( j(\cdot,y))F(dy),
\end{align*}
where $F$ is a non-negative measure on a measurable space $E$, 
 the jump size map $j$ is given by $  j(\cdot,y):=(\langle \j_1(y),\cdot\rangle,\dots,\langle \j_d(y),\cdot\rangle)$, for some suitable $ \j(y):=(\j_1(y), \dots,\j_d(y))$ belonging to  $ ({\SS})^d$ for all $y\in E$.
 
 One can show that if either of the following conditions, \ref{enu:lambda i} or \ref{enu:lambda ii}, is additionally satisfied,
 \begin{enumerate}
    \item \label{enu:lambda i}$\lambda(\cdot):=\langle \bm{\lambda},\cdot\rangle$, for some $\bm{\lambda}\in \SS$ such that $\langle \bm{\lambda},\x\rangle>0$ for all $\x\in \S$  and $j(\cdot,y)\equiv j(y)$;
    \item \label{enu:lambda ii} $\lambda(\cdot):=\frac{\langle\tilde{\bm{\lambda}},\cdot\rangle}{\langle \j^*,\cdot \rangle^2}1_{\langle \j^*,\cdot\rangle\neq0},$ where $\tilde
    {\bm{\lambda}}\in T^2(\Rset^d)$, ${\j}^*\in T^1(\Rset^d)$, and ${\j}_i(y)=H(y){\j}^*$ for all $i=1,\dots,d$  and $y\in E$,
\end{enumerate}
and the drift and diffusion coefficients $b,a$ are polynomials on $\Rset^d$ of degree one and two, respectively, then under suitable integrability conditions on $K$, the truncated signature process of the jump diffusion corresponding to such sig-coefficients is a polynomial jump-diffusion on the truncated tensor algebra. Note that the parametrizations \ref{enu:lambda i} and \ref{enu:lambda ii} correspond to a path-dependent intensity, meaning that the resulting processes go beyond the classical framework of polynomial processes.

\end{remark}

Finally, we give an example in which we implement the formula in Corollary \ref{coro:expectedsigpoly}.

The corresponding code is available at \url{https://github.com/fprimavera/expected_signature_polyjumpdiffusions.git}.

\begin{example}\label{ex:OU-jumps}
We consider a 3-dimensional process $\widehat{X}_t = (t, X^{(1)}_t, X^{(2)}_t)$, for $t\in[0,T]$, whose last two coordinates follow Ornstein--Uhlenbeck dynamics driven by independent Brownian motions and a common Poisson jump clock:
$$
dX^{(i)}_t = (\alpha_i + \beta_i X^{(i)}_{t-})\,dt + \sigma_i\,dW^{(i)}_t + dJ^{(i)}_t,
\qquad i=1,2,\qquad X_0=(c_1,c_2)\in\mathbb{R}^2,
$$
with $W^{(i)}$ independent 1-dimensional Brownian motions and jumps $J^{(i)}_t = \sum_{k=1}^{N_t} Z^{(i)}_k$, $N_t\sim\mathrm{Poisson}(\lambda t)$. The jump sizes $Z_k=(Z^{(1)}_k,Z^{(2)}_k)$ are drawn from the $2$-dimensional Gaussian distribution  $F=\mathcal N(\mu,\Sigma)$, with mean $\mu=(\mu_1,\mu_2)$, and covariance matrix $\Sigma=\mathrm{diag}(\Sigma_{11},\Sigma_{22})$.

\paragraph{Characteristics.}
The coefficients of the characteristics of $\widehat X$ are
$$
b(x_1,x_2,x_3) = \begin{pmatrix} 1 \\ \alpha_1 + \beta_1 x_2 + \lambda \mu_1 \\ \alpha_2 + \beta_2 x_3 + \lambda \mu_2 \end{pmatrix},
\qquad
a(x_1,x_2,x_3) = \begin{pmatrix} 0 & 0 & 0 \\ 0 & \sigma_1^2 & 0 \\ 0 & 0 & \sigma_2^2 \end{pmatrix},
$$
with a state-independent jump kernel affecting only the OU coordinates,
$$
K(x_1,x_2,x_3,d\xi) = \lambda\,(j_*F)(d\xi), \qquad j(z_1,z_2)=(0,z_1,z_2).
$$

\paragraph{Tensors $Q^{(r)}$.}
For $\boldsymbol x\in\mathcal S$,
$$
Q^{(1)}(\boldsymbol x) = \begin{pmatrix} 1 \\[2pt] \alpha_1 + \beta_1\langle\epsilon_2,\boldsymbol x\rangle + \lambda\mu_1 + \beta_1 c_1 \\[2pt] \alpha_2 + \beta_2\langle\epsilon_3,\boldsymbol x\rangle + \lambda\mu_2 + \beta_2 c_2 \end{pmatrix},
$$
$$
Q^{(2)}(\boldsymbol x) = \frac12\begin{pmatrix} 0 & 0 & 0 \\ 0 & \sigma_1^2+\lambda(\mu_1^2+\Sigma_{11}) & \lambda\mu_1\mu_2 \\ 0 & \lambda\mu_1\mu_2 & \sigma_2^2+\lambda(\mu_2^2+\Sigma_{22}) \end{pmatrix},$$
$$
Q^{(r)}(\boldsymbol x) = \frac{\lambda}{r!}\,\mathbb{E}_F[j(Z)^{\otimes r}], \quad r\geq 3,\quad Z\sim F.
$$
Using the notation of \eqref{eq: notationderivative},
$$
\alpha_1+\lambda\mu_1+\beta_1 c_1 = b_2(\widehat X_0), \qquad \beta_1 = b_2^{(2)}(\widehat X_0),
$$
$$
\alpha_2+\lambda\mu_2+\beta_2 c_2 = b_3(\widehat X_0), \qquad \beta_2 = b_3^{(3)}(\widehat X_0).
$$

\paragraph{Truncation level and jump moments.}
We compute $\mathbb E\big[\widehat{\mathbb X}_T^{\leq 5}\big]$, the expected signature truncated at level $N=5$, by assembling the matrix $L$ from the coefficients of $Q^{(r)}$, $r\le 5$. The non-null jump moments of order 3 appearing in $Q^{(3)}$ are
$$
\mathbb E_F[\xi_1^3]=\mu_1^3+3\mu_1\Sigma_{11}, \qquad \mathbb E_F[\xi_2^3]=\mu_2^3+3\mu_2\Sigma_{22},
$$
$$
\mathbb E_F[\xi_1^2\xi_2]=\mu_2(\mu_1^2+\Sigma_{11}), \qquad \mathbb E_F[\xi_1\xi_2^2]=\mu_1(\mu_2^2+\Sigma_{22}).
$$
The jump moments of order 4 appearing in $Q^{(4)}$ are
$$
\mathbb E_F[\xi_1^4]=\mu_1^4+6\mu_1^2\Sigma_{11}+3\Sigma_{11}^2, \qquad \mathbb E_F[\xi_2^4]=\mu_2^4+6\mu_2^2\Sigma_{22}+3\Sigma_{22}^2,
$$
$$
\mathbb E_F[\xi_1^3\xi_2]=\mu_2(\mu_1^3+3\mu_1\Sigma_{11}), \qquad \mathbb E_F[\xi_1\xi_2^3]=\mu_1(\mu_2^3+3\mu_2\Sigma_{22}),
$$
$$
\mathbb E_F[\xi_1^2\xi_2^2]=(\mu_1^2+\Sigma_{11})(\mu_2^2+\Sigma_{22}).
$$
Finally, the jump moments of order 5 appearing in $Q^{(5)}$ are
$$
\mathbb E_F[\xi_1^5]=\mu_1^5+10\mu_1^3\Sigma_{11}+15\mu_1\Sigma_{11}^2, \qquad \mathbb E_F[\xi_2^5]=\mu_2^5+10\mu_2^3\Sigma_{22}+15\mu_2\Sigma_{22}^2,
$$
$$
\mathbb E_F[\xi_1^4\xi_2]=\mu_2(\mu_1^4+6\mu_1^2\Sigma_{11}+3\Sigma_{11}^2), \qquad \mathbb E_F[\xi_1\xi_2^4]=\mu_1(\mu_2^4+6\mu_2^2\Sigma_{22}+3\Sigma_{22}^2),
$$
$$
\mathbb E_F[\xi_1^3\xi_2^2]=(\mu_2^2+\Sigma_{22})(\mu_1^3+3\mu_1\Sigma_{11}), \qquad \mathbb E_F[\xi_1^2\xi_2^3]=(\mu_1^2+\Sigma_{11})(\mu_2^3+3\mu_2\Sigma_{22}).
$$

\paragraph{Numerical validation and computational cost.}
This computation is available in a notebook at \url{https://github.com/fprimavera/expected_signature_polyjumpdiffusions}, which also benchmarks the resulting expected signature against a Monte Carlo estimate, obtained by simulating $\widehat X$ on a fine time grid and computing the signature of the resulting piecewise-linear path, as discussed in Remark~\ref{remMarkus}. For readability, Figure~\ref{fig:mc-vs-poly} restricts the componentwise comparison to levels $1$ to $3$ (39 non-trivial components): all of them agree with the Monte Carlo estimate within 2 standard errors. 

The comparison of computational cost, on the other hand, is carried out at the full truncation level $N=5$, where the closed-form approach is most advantageous. Assembling $L$ (a $364\times364$ matrix at $N=5$) and evaluating $\exp(TL)$ takes under a second, whereas obtaining a Monte Carlo estimate of comparable accuracy - simulating $10^4$ paths on a grid of $10^4$ time steps and computing their signatures up to level 5 - takes on the order of a minute, roughly two orders of magnitude slower. Moreover, the cost of the closed-form formula grows only mildly with the truncation level $N$ (well under a second already at $N=5$), unlike the simulation-based estimate, whose signature-computation step alone grows from a few seconds at $N=3$ to about 48 seconds at $N=5$.

\begin{figure}[htbp]
\centering
\includegraphics[width=\textwidth]{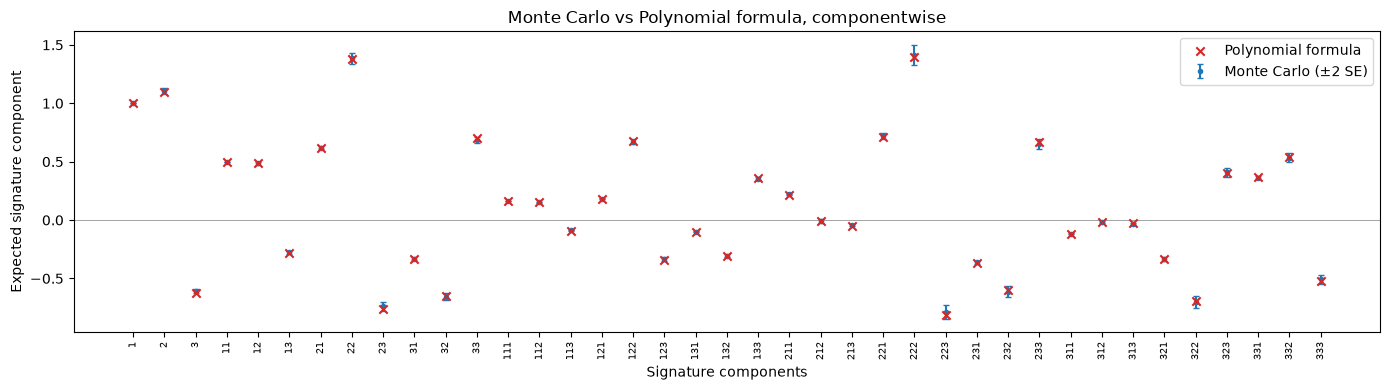}
\caption{Componentwise comparison between the Monte Carlo estimate ($\pm 2$ standard errors) and the polynomial formula for $\mathbb E[\widehat{\mathbb X}_T^{\leq 3}]$ in Example~\ref{ex:OU-jumps}, at $T=1$, restricted to levels $1$--$3$ for readability. The computational-cost comparison discussed in the text refers instead to the full truncation level $N=5$.}
\label{fig:mc-vs-poly}
\end{figure}
\end{example}

\section{Score-based generative modeling with entire sig-jump-diffusions}
\label{sec:score}

The goal of this section is to apply entire sig-jump-diffusions for score-based generative  modeling of time-series data in the spirit of \cite{song2021score}. More precisely, we encode the empirical law of some time series data $(z_t)_{t\in [-S,0]}$ over the time horizon $[-S,0]$ with $S >0$
via the empirical law  of its signature  truncated at level $N$ and denoted by $\z_0$.  In the spirit of \cite{song2021score} we then perturb the data $\z_0$ via a lifted entire sig-jump-diffusion of the form \eqref{entsig2} starting at $\z_0$ such that it converges to some invariant measure\footnote{A probability measure $\pi$ is called invariant for a Markov process with transition semigroup $(P_t)_{t\geq 0}$ if $\int P_tf(x)\pi(dx)=\int f(x)\pi(dx)$ for every $t\geq 0$. Equivalently, if the initial state has law $\pi$, then the law of the process remains $\pi$ at all subsequent times.}. Sampling from this invariant distribution and using a reverse process denoted by $\mathbb{Y}^{\leq N}$ that generates the same marginals as the forward process \eqref{entsig2} then allows us to generate new samples mimicking the law of  $\z_0$. The essential ingredient in 
the reverse process construction is a \emph{generalized score function} (that reduces up to a factor to the gradient of the marginal logarithmic density in the constant diffusion case) introduced in Definition~\ref{def:generalized-score}.
As in standard generative diffusion models this score function has to be learned from data. Here, the framework of entire sig-jump-diffusions in the sense of Section~\ref{sec:holo}  and the expected signature formula become crucial as this allows to reduce the
 matching of the generalized score function  to the solution of linear ODEs and a deterministic optimization problem.

On the level of the original time series data there is also a nice interpretation. Indeed, consider a process given by the data time series in $[-S,0]$ and continued by an entire sig-jump-diffusion afterwards, i.e.,
$$X_t=\begin{cases}
    z_t&t\in [-S,0],\\
    z_0+\int_0^t\bigl(\overline b(\X_{s-}^{\leq N})ds
+
\overline\sigma( \X_{s-}^{\leq N})dW_s
+
\int_E
\overline\delta( \X_{s-}^{\leq N},z)
\bigl(\mu(ds,dz)-F(dz)ds\bigr)\bigr) &t\in [0,T].
\end{cases}$$
Then, as explained in Remark~\ref{remcontsig}, its signature $(\X_{t}^{\leq N})_{t\in [0,T]}$ satisfies \eqref{entsig2} for ${\bm g}=\z_0$.
This means that if the time-series data in $[-S,0]$ are actually samples of the same entire sig-jump-diffusion, $\operatorname{Law}(\mathbb{X}_{-S,0}^{\leq N})=\operatorname{Law}(\mathbb{Y}_T^{\leq N})$ holds true. This consistency relation can be exploited to specify $\overline{b}, \overline{\sigma}, \overline{\delta}$ by additionally calibrating an entire sig-jump-diffusion to the data on $[-S,0]$.

To introduce the whole  generative modeling procedure we start by considering
marginal-preserving reverse dynamics for general jump-diffusions and apply this then subsequently to the lifted entire sig-jump-diffusion of the form \eqref{entsig2}.

\subsection{Marginal-preserving reverse dynamics for general jump-diffusions}
\label{sec:marginal-reverse-jump}

 Let $(Z_t)_{t\in[0,T]}$ be an $\mathbb{R}^n$-valued jump-diffusion satisfying
\begin{equation}
\label{eq:forward-jump-diffusion}
\begin{aligned}
dZ_t
={}&
b(Z_{t-})\,dt
+\sigma(Z_{t-})\,dW_t +
\int_{E}
\gamma(Z_{t-},\xi)
\bigl(
\mu(dt,d\xi)-F(d\xi)\,dt
\bigr), \quad Z_0 \sim p_0
\end{aligned}
\end{equation}
where $W$ is an $\mathbb{R}^m$-valued Brownian motion, $(E,\mathcal E, F)$ some measure space, $\mu$ a Poisson random measure on
$[0,T]\times E$ with compensator $F(d\xi)\,dt$, and
$$
b:\mathbb{R}^n\to\mathbb{R}^n,\quad
\sigma:\mathbb{R}^n\to\mathbb{R}^{n\times m},\quad
\gamma:\mathbb{R}^n\times E \to\mathbb{R}^n$$
 are measurable functions.
We set
$
a(x)=\sigma(x)\sigma(x)^\top
$ and assume that
$$
\int_{E}
\|\gamma(x,\xi)\|^2\,F(d\xi)<\infty,
\qquad \text{for all } x\in\mathbb{R}^n,
$$
and that \eqref{eq:forward-jump-diffusion} admits a unique weak solution, in the sense that its martingale problem is well posed.  These assumptions can, for example, be ensured by standard Lipschitz and linear-growth conditions on $b$, $\sigma$, and $\gamma$. 
We furthermore suppose that its one-time marginal laws admit densities $(p_t)_{t \in [0,T]}$ and the associated forward Kolmogorov equation admits a unique weak solution.\\

For $f\in C_c^2(\mathbb{R}^n)$, the generator of \eqref{eq:forward-jump-diffusion}, with the identity as ``truncation'' function, is
$$
\mathcal G f(x)
=
b(x)^\top\nabla f(x)+\mathcal Kf(x),
$$
where
$$
\begin{aligned}
\mathcal Kf(x)
={}&
\frac12
\operatorname{Tr}
\left(
a(x)\nabla^2 f(x)
\right)
+
\int_{E}
\Big(
f(x+\gamma(x,\xi))-f(x)
-\gamma(x,\xi)^\top\nabla f(x)
\Big)
F(d\xi).
\end{aligned}
$$
Our goal is now to introduce a reverse process whose marginal laws coincide with those of $Z$ at reversed time. To this end we need the notion of a generalized score function.

\begin{definition}
\label{def:generalized-score}
A measurable function
$
S_t:\mathbb{R}^n\to\mathbb{R}^n
$
is called a generalized score function associated with $\mathcal K$ and $p_t$ if
$p_tS_t\in L^1_{\mathrm{loc}}(\mathbb{R}^n;\mathbb{R}^n)$ and
for every
$f\in C_c^2(\mathbb{R}^n)$,
\begin{equation}
\label{eq:weak-score}
\int_{\mathbb{R}^n}
\mathcal Kf(x)\,p_t(x)\,dx
=
-
\int_{\mathbb{R}^n}
\langle S_t(x), \nabla f(x) \rangle\,p_t(x)\,dx.
\end{equation}
\end{definition}

This definition can be viewed as an operator-dependent generalization of the classical score identity underlying score matching as considered in \cite{hyvarinen2005estimation}. 
In the specific setup of L\'evy noise,
 the generalized score defined above corresponds to the score considered in \cite{B:26} for the finite jump activity case and to the one in \cite{yoon2023score} for
the isotropic $\alpha$-stable case, see Remark \ref{rem:special case} for further details.

\begin{remark} \label{rem:special case}
We here show that the generalized score introduced in Definition~\ref{def:generalized-score} reduces to the usual
additional drift component of the reverse process
arising in pure diffusions and  the score function of isotropic $\alpha$-stable L\'evy processes considered in \cite{yoon2023score}.\\

First, consider the purely diffusive case with
\[
\mathcal K f(x)
=
\frac12
\operatorname{Tr}\!\left(
a(x)\nabla^2 f(x)
\right).
\]
Assuming sufficient regularity to perform integration by parts, we have
\begin{align*}
\int_{\mathbb{R}^n}
\mathcal K f(x)\,p_t(x)\,dx
&=
\frac12
\sum_{i,j=1}^n
\int_{\mathbb{R}^n}
a_{ij}(x)
\partial_{ij}f(x)\,
p_t(x)\,dx
\\
&=
-\frac12
\sum_{i,j=1}^n
\int_{\mathbb{R}^n}
\partial_j\!\left(
a_{ij}(x)p_t(x)
\right)
\partial_i f(x)\,dx.
\end{align*}
Comparing with \eqref{eq:weak-score} yields
$$
S_t(x)
=
\frac{1}{2p_t(x)}
\operatorname{div}\!\left(
a(x)p_t(x)
\right)= \frac12\operatorname{div}a(x)
+
\frac12a(x)\nabla\log p_t(x),
$$
where the $i$-th component of $\operatorname{div}(a(x))$ is given by
\begin{align}\label{eq:divergence}
\big(\operatorname{div} a(x)\big)_i
=
\sum_{j=1}^n \partial_{j} a_{ij}(x),
\qquad i=1,\dots,n.
\end{align}
In particular, if $a(x)=\sigma^2 I_n$, then
$
S_t(x)
=
\frac{\sigma^2}{2}\nabla\log p_t(x),
$
so that, up to the natural factor $\sigma^2/2$, the generalized score coincides with the classical score
$\nabla\log p_t$.

Second, let $\mathcal K$ be the generator of an isotropic symmetric $\alpha$-stable L\'evy process with
\[
\mathcal K
=
-(-\Delta)^{\alpha/2},
\qquad
0<\alpha<2.
\]
The operator $\mathcal K$ is self-adjoint, and
\[
-(-\Delta)^{\alpha/2}p
=
\nabla\cdot
\left(
(-\Delta)^{(\alpha-2)/2}\nabla p
\right).
\]
Hence
\begin{align*}
\int_{\mathbb{R}^n}
\mathcal K f(x)\,p_t(x)\,dx
&=
\int_{\mathbb{R}^n}
f(x)\,\mathcal K p_t(x)\,dx
\\
&=
\int_{\mathbb{R}^n}
f(x)\,
\nabla\cdot
\left(
(-\Delta)^{(\alpha-2)/2}
\nabla p_t(x)
\right)
dx
\\
&=
-
\int_{\mathbb{R}^n}
\nabla f(x)^\top
(-\Delta)^{(\alpha-2)/2}
\nabla p_t(x)\,dx.
\end{align*}
Comparison with \eqref{eq:weak-score} therefore gives
$$
S_t(x)
=
\frac{
(-\Delta)^{(\alpha-2)/2}\nabla p_t(x)
}{
p_t(x)
}.
$$
This is precisely the fractional score introduced for isotropic $\alpha$-stable L\'evy diffusion models in \cite{yoon2023score}.

\end{remark}

We now obtain a family of reverse dynamics analogous to the family introduced in \cite{popov2025improved} for isotropic $\alpha$-stable L\'evy processes.

\begin{theorem}
\label{thm:marginal-reverse}
Let $(Z_t)_{t\in[0,T]}$ satisfy \eqref{eq:forward-jump-diffusion}, and suppose that a generalized score $S_t$ satisfying
\eqref{eq:weak-score} exists for every $t\in(0,T]$.
Let
$
\eta:[0,T]\to[0,\infty)
$
be a bounded measurable function and set
$
s=T-t
$.
Let $\overline W$ be an $\mathbb{R}^m$-valued Brownian motion, and let
$\overline\mu_\eta$ be a Poisson random measure with compensator
$
\eta_{T-s}F(d\xi)\,ds
$.
Consider the $\mathbb{R}^n$-valued process $(Y_s)_{s\in[0,T]}$ satisfying
\begin{equation}
\label{eq:reverse-SDE}
\begin{aligned}
dY_s
={}&
\Big(
-b(Y_{s-})
+
(1+\eta_{T-s})
S_{T-s}(Y_{s-})
\Big)\,ds
+
\sqrt{\eta_{T-s}}\,
\sigma(Y_{s-})\,d\overline W_s\\
&+
\int_{E}
\gamma(Y_{s-},\xi)
\Big(
\overline\mu_\eta(ds,d\xi)
-
\eta_{T-s}F(d\xi)\,ds
\Big),
\end{aligned}
\end{equation}
with initial distribution
$
Y_0\sim p_T(x)\,dx$.
Assume that the martingale problem associated with
\eqref{eq:reverse-SDE}
is well posed and that the associated forward Kolmogorov equation admits a unique weak solution. Then
\begin{equation}
\label{eq:marginal-equality}
\operatorname{Law}(Y_s)
=
\operatorname{Law}(Z_{T-s}),
\qquad
0\le s\le T.
\end{equation}
In particular, the reverse process has marginal densities
$
q_s(x)=p_{T-s}(x)
$.
\end{theorem}

\begin{proof}
For $s=T-t$, the generator of the process $Y$ in
\eqref{eq:reverse-SDE} is
\begin{equation}
\label{eq:reverse-generator}
\mathcal G_s^{\mathrm{rev},\eta}f(x)
=
\Big[
-b(x)+(1+\eta_{T-s})S_{T-s}(x)
\Big]^\top\nabla f(x)
+
\eta_{T-s}\mathcal Kf(x).
\end{equation}

Let $q_s$ denote the density of $Y_s$. By uniqueness of both Kolmogorov forward equations, it is enough 
to verify that
$
q_s(x)=p_{T-s}(x)
$
solves the Kolmogorov  forward equation associated with
\eqref{eq:reverse-generator}.
For $f\in C_c^2(\mathbb{R}^n)$, we have by (the weak formulation of) the Kolmogorov  forward equation of the forward process given in \eqref{eq:forward-jump-diffusion}
\[
\frac{d}{dt}
\int_{\mathbb{R}^n}
f(x)p_t(x)\,dx
=
\int_{\mathbb{R}^n}
\left(
b(x)^\top\nabla f(x)+\mathcal Kf(x)
\right)
p_t(x)\,dx.
\]
By the defining identity of $S_t$ in \eqref{eq:weak-score},
this becomes
\begin{equation*}
\frac{d}{dt}
\int_{\mathbb{R}^n} f(x) p_t(x)\,dx
=
\int_{\mathbb{R}^n}
\bigl(
b(x)-S_t(x)
\bigr)^\top
\nabla f(x)\,p_t(x)\,dx.
\end{equation*}
Setting $
q_s=p_{T-s}$, we thus have
\begin{align}
\frac{d}{ds}
\int f(x)q_s(x)\,dx
&=
\frac{d}{ds}
\int f(x)p_{T-s}(x)\,dx
\nonumber\\
&=
\int
\bigl(
-b(x)+S_{T-s}(x)
\bigr)^\top
\nabla f(x)\,p_{T-s}(x)\,dx.
\label{eq:reverse-target}
\end{align}

On the other hand, using \eqref{eq:reverse-generator}, we get

\begin{align*}
\int
\mathcal G_s^{\mathrm{rev},\eta}f(x)
p_{T-s}(x)\,dx=&
\int
\Big(
-b(x)+(1+\eta_{T-s})S_{T-s}(x)
\Big)^\top
\nabla f(x)\,p_{T-s}(x)\,dx\\ &+
\eta_{T-s}
\int
\mathcal Kf(x)\,p_{T-s}(x)\,dx.
\end{align*}
Applying \eqref{eq:weak-score} once more yields
\begin{align}
\int
\mathcal G_s^{\mathrm{rev},\eta}f(x)
p_{T-s}(x)\,dx
&=
\int
\Big(
-b(x)+(1+\eta_{T-s})S_{T-s}(x)-\eta_{T-s}S_{T-s}(x)
\Big)^\top
\nabla f(x)\,p_{T-s}(x)\,dx
\nonumber\\
&=
\int
\bigl(
-b(x)+S_{T-s}(x)
\bigr)^\top
\nabla f(x)\,p_{T-s}(x)\,dx.
\label{eq:eta-cancellation}
\end{align}
Since the right-hand sides of \eqref{eq:reverse-target} and
\eqref{eq:eta-cancellation} coincide, we conclude that $q_s=p_{T-s}$ satisfies the Kolmogorov forward equation of the backward process given in  \eqref{eq:reverse-SDE}, i.e.,
\[
\frac{d}{ds}
\int f(x)q_s(x)\,dx=
\int
\mathcal G_s^{\mathrm{rev},\eta}f(x)
q_{s}(x)\,dx.
\]
The fact that
$
q_0=p_T
$
and the uniqueness of the martingale problem of the backward process imply \eqref{eq:marginal-equality}.
\end{proof}

\begin{remark}
For the special choice  $\eta_t\equiv0$, all stochastic terms in
\eqref{eq:reverse-SDE} disappear and one obtains the following ODE
\begin{equation}
\label{eq:probability-flow}
dY_s
=
\Big(
-b(Y_s)+S_{T-s}(Y_s)
\Big)\,ds.
\end{equation}
This deterministic flow \eqref{eq:probability-flow} has the same one-time marginal distributions $p_{T-s}$ as the forward jump-diffusion.\\

For the other special case, namely $\eta_t\equiv1$, only the drift is modified by the score such that the reverse dynamics have the same local diffusion covariance and the same jump intensity as the forward process:
\begin{equation*}
\begin{aligned}
dY_s
={}&
\bigl(
-b(Y_{s-})+2S_{T-s}(Y_{s-})
\bigr)\,ds
+
\sigma(Y_{s-})\,d\overline W_s
+
\int_{E}
\gamma(Y_{s-},\xi)
\bigl(
\overline\mu(ds,d\xi)-F(d\xi)\,ds
\bigr).
\end{aligned}
\end{equation*}
\end{remark}

\begin{remark} Note that
the process \eqref{eq:reverse-SDE} is in general not the true pathwise time reversal of
\eqref{eq:forward-jump-diffusion}. Rather, it belongs to a family of processes whose one-time marginal distributions agree exactly with those of the forward process under time reversal. This coincides with the notion of reverse dynamics used in Theorem~1 of
\cite{popov2025improved}. In contrast, the pathwise time reversal of a Markov process with jumps 
generally involves a modification of its jump characteristics depending 
on the marginal laws; we refer to \cite{CL:22} for a 
characterization of the semimartingale characteristics of the 
time-reversed process in a pathwise sense.
\end{remark}

\subsection{Marginal-preserving reverse dynamics for entire  sig-jump-diffusions}\label{sec:liftedsigjump}

We now apply the preceding construction to entire sig-jump-diffusions with truncated signature-dependent coefficients as considered in Section \ref{sec:holo}, more precisely to the lifted equation $\mathbb{Z}^{\leq N}$ given in \eqref{entsig2} with some initial value $\z_0 \sim p_0=p_{\mathrm{data}}$. Note that we are in particular assuming that the process defined in \eqref{entsig} is an entire sig-jump-diffusion in the sense of Remark~\ref{Rem:sigisholo}.
With the definition of $b, \sigma, \gamma$ given there and with $a(\x)=\sigma(\x)\sigma(\x)^{\top}$, 
the generator associated with \eqref{entsig2} is given by 
$$
\mathcal G^N f(\x)
=
\langle b(\x),\nabla f(\x)\rangle
+
\mathcal K^Nf(\x),
$$
where
\begin{align*}
\mathcal K^Nf(\x)
=
\frac12
\operatorname{Tr}
\left(
a(\x)\nabla^2f(\x)
\right)
+
\int_E
\Big(
f(\x+\gamma(\x,z))-f(\x)
-
\langle\gamma(\x,z),\nabla f(\x)\rangle
\Big)
F(dz).
\end{align*}
We set $P_t^N=\operatorname{Law}(\mathbb Z_t^{\leq N})
$ and adapt the definition of the
 generalized score function 
 to the current setting, accounting for the fact that we do not  have a density on $T^N(\mathbb{R}^d)$ since the law 
 $P_t^N$ is only supported
on (possibly a subset of) group-like elements $
G^N(\mathbb R^d)$.
We additionally assume that the associated Kolmogorov forward equation  for $P_t^N$ admits a unique weak solution.

\begin{definition}
\label{def:sig-score}
A measurable function
$
S_t^N: T^N(\mathbb R^d)
\longrightarrow
 T^N(\mathbb R^d)
$
is called a generalized score function of $\mathbb{Z}^{\leq N}$ at
time $t$ if it is integrable with respect to $P_t^N$ and
\begin{equation}
\label{eq:weak-sig-score}
\int
\mathcal K^Nf(\x)\,P_t^N(d\x)
=
-
\int
\langle S_t^N(\x),\nabla f(\x)\rangle
\,P_t^N(d\x)
\end{equation}
for every $f\in C_c^2( T^N(\mathbb R^d))$ for which the
above expressions are finite.
\end{definition}

We can now apply the marginal-preserving reverse construction directly
to the signature process.

\begin{corollary}
 Assume that, for every $t\in(0,T]$, there exists a
generalized score $S_t^N$ satisfying
\eqref{eq:weak-sig-score}.
Let
$
\eta:[0,T]\to[0,\infty)
$
be bounded and measurable, and introduce the increasing reverse-time
variable $s=T-t$. Let $\overline W$ be a $d$-dimensional Brownian
motion and let $\overline\mu_\eta$ be a Poisson random measure with
compensator
$
\eta_{T-s}F(dz)\,ds.
$
Consider
\begin{equation}
\label{eq:reverse-signature-SDE}
\begin{aligned}
d\mathbb Y^{\leq N}_s
={}&
\Big(
-b(\mathbb Y^{\leq N}_{s-})
+
(1+\eta_{T-s})
S_{T-s}^N(\mathbb Y^{\leq N}_{s-})
\Big)\,ds
+
\sqrt{\eta_{T-s}}\,
\sigma(\mathbb Y^{\leq N}_{s-})\,d\overline W_s
\\
&+
\int_E
\gamma(\mathbb Y^{\leq N}_{s-},z)
\Big(
\overline\mu_\eta(ds,dz)
-
\eta_{T-s}F(dz)\,ds
\Big),
\end{aligned}
\end{equation}
with
$
\mathbb Y^{\leq N}_0\sim P_T^N.
$
Assume that the martingale problem associated with
\eqref{eq:reverse-signature-SDE} is well posed on $T^N(\mathbb{R}^d)$ and that
 the associated Kolmogorov forward equation  admits a unique weak solution.
Then
$$
\operatorname{Law}(\mathbb Y_s^{\leq N})
=
P_{T-s}^N
=
\operatorname{Law}(\mathbb Z^{\leq N}_{T-s}),
\qquad
0\leq s\leq T.
$$
\end{corollary}

\begin{remark}
Since the law $P_t^N$ of the forward truncated signature is supported on (subsets of) group-like elements
$G^N(\mathbb R^d)$, the marginal identity
$
\operatorname{Law}(\mathbb Y^{\leq N}_s)
=
P_{T-s}^N
$
implies that
$
\mathbb P\bigl(\mathbb Y^{\leq N}_s\in G^N(\mathbb R^d)\bigr)=1
$
for every fixed $s\in[0,T]$. Since $G^N(\mathbb R^d)$ is closed and
$\mathbb Y^{\leq N}$ has c\`adl\`ag paths, it follows by considering a
countable dense set of times and using right-continuity that
\[
\mathbb P\Bigl(
\mathbb Y^{\leq N}_s\in G^N(\mathbb R^d)
\text{ for all }s\in[0,T]
\Bigr)=1.
\]
Thus, the marginal-preserving reverse dynamics remain in $G^N(\mathbb{R}^d)$. 
It is, however, not automatic that
$\mathbb Y^{\leq N}$ is the truncated signature of its own first-level
process,
since the generalized score may
introduce higher-level drift components that do not coincide with the
canonical signature lift of the first-level reverse dynamics. Additional
compatibility conditions on the generalized score would therefore be needed.
\end{remark}

\subsection{Learning the generalized score function of entire  sig-jump-diffusions via an entire  approximation}\label{sec:learning}

We now aim  to use the above reverse process to generate  signatures mimicking the law $p_{\mathrm{data}}$ of the signature $\z_0$ stemming from some time-series data $(z_t)_{t \in [-S,0]}$ and truncated at level $N$. To this end the generalized score function has to be learned. For this purpose we approximate  it by an entire  function in time and the signature components,  and use Corollary~\ref{Corsupernew} for the computation of the expected signature to obtain a tractable training objective. 
 Consider a function
\[
h_{\u}:\mathbb{R} \times T^N(\mathbb{R}^d) \to \mathbb \mathbb{R}
\]
which is jointly entire  in $(t,\x)$. We write its power series representation as
\[
h_{\u}(t,\x)
=
\sum_{k\geq 0}\sum_{\alpha\in\mathbb N_0^{D_N}}
u_{k,\alpha}\,t^k\x^\alpha
\]
with multi-indices $\alpha \in \mathbb{N}^{D_N}_0$ and $D_N := \dim  T^N(\mathbb R^d) $
and define the score approximation by
$$
S^N_t(\x)\approx g_{\u}(t,\x)
:=
\nabla_{\x}h_{\u}(t,\x).
$$

For $\mathbb{Z}^{\leq N}_0=\z_0$, let $S_t(\cdot\mid \z_0)$ denote a conditional generalized score, characterized weakly by
\begin{equation}
\label{eq:conditional-weak-score}
\mathbb E\left[
\mathcal K^N f(\mathbb Z_t^{\leq N})
\mid \mathbb{Z}^{\leq N}_0=\z_0
\right]
=
-
\mathbb E\left[
\left\langle
S_t(\mathbb Z_t^{\leq N}\mid \z_0),
\nabla f(\mathbb Z_t^{\leq N})
\right\rangle
\mid \mathbb{Z}^{\leq N}_0=\z_0
\right]
\end{equation}
for all admissible test functions $f$.  We consider now the following  denoising generalized score-matching objective 
$$
\operatorname{argmin}_{\u}\int_0^T
\lambda(t)
\mathbb E_{\z_0 \sim p_{\mathrm{data}}}
\left[
\mathbb E\left[
\left\|
g_{\u}(t,\mathbb Z_t^{\leq N})
-
S_t(\mathbb Z_t^{\leq N}\mid \z_0)
\right\|^2
\,\middle|\,\Z^{\leq N}_0={\bf z}_0
\right]
\right]dt ,
$$
with $\lambda$ some positive weight function. By Fubini-Tonelli's theorem this is equivalent to
\begin{equation}
\label{eq:hol-denoising-loss}
\operatorname{argmin}_{\u}
\mathbb E_{\z_0 \sim p_{\mathrm{data}}}
\left[
\mathbb E\left[
\int_0^T
\lambda(t)\left\|
g_{\u}(t,\mathbb Z_t^{\leq N})
-
S_t(\mathbb Z_t^{\leq N}\mid \z_0)
\right\|^2 dt
\,\middle|\,\Z^{\leq N}_0={\bf z}_0
\right]
\right],
\end{equation}
and also to
\begin{equation*}
\operatorname{argmin}_{\u}
\mathbb E_{\z_0\sim p_{\mathrm{data}}}
\left[
\mathbb E\left[\int_0^T
\lambda(t)\left(
\left\|
g_{\u}(t,\mathbb Z_t^{\leq N})\right\|^2
-
2  \langle g_{\u}(t,\mathbb Z_t^{\leq N}), S_t(\mathbb Z_t^{\leq N}\mid \z_0) \rangle
\right) dt \, \middle|\,\Z^{\leq N}_0=\z_0
\right]
\right].
\end{equation*}
Applying
\eqref{eq:conditional-weak-score} with
$f(\x)=h_{\u}(t,\x)$ yields
\begin{equation}
\label{eq:tractable-hol-loss}
\operatorname{argmin}_{\u}
\mathbb E_{\z_0\sim p_{\mathrm{data}}}
\left[
\mathbb E
\left[\int_0^T
\lambda(t) \left(
\left\|
\nabla_{\x}h_{\u}
(t,\mathbb Z_t^{\leq N})
\right\|^2
+
2\mathcal K^N
h_{\u}
(t,\mathbb Z_t^{\leq N}) \right) dt \,\middle|\,\Z^{\leq N}_0=\z_0
\right]
\right].
\end{equation}
Here, $\mathcal K^N$ acts only on the signature variable $\x$. In particular,
the unknown generalized score no longer appears in
\eqref{eq:tractable-hol-loss}.\\

The objective \eqref{eq:tractable-hol-loss} can now be computed via Corollary \ref{Corsupernew}. Indeed, define
$$
q_{\u}(t,\x)
:=
\left\|
\nabla_{\x}h_{\u}(t,\x)
\right\|^2
+
2\mathcal K^N h_{\u}(t,\x).
$$
Under suitable conditions on $h_\u$ \cite[Lemma~A.1]{CPS:24holo}, the conditions on $\mathcal{K}^N$ yield that  $q_{\u}$ is  a jointly entire  function. 
To incorporate the explicit time dependence, we shall use the time-extended version of \eqref{entsig2}, which we denote by $\widehat{\mathbb{Z}}^{\leq N}$. Its starting value is $\widehat{\z}_0$ which is the signature of the time-series data. More precisely, $\widehat{\z}_0^{-1} \otimes \widehat{\mathbb{Z}}^{\leq N}$ is the signature of $\widehat{X}_t^{{\widehat\z_0}}:=(t, X_t)^{\widehat{\z}_0}$ given by \eqref{entsig} extended to $\mathbb{R}^{d+1}$   with $\overline{b}_0=1$ and $\overline{\sigma}_{0i}=\overline{\delta}_0=0$. Here, the $0$-index refers to time.
Using this we can express
$$q_{\u}(t,\mathbb Z_t^{\leq N})=\langle \v(\u), \widehat{\mathbb{Z}}_t \rangle,$$ where $\v(\u)$ denotes  the corresponding  coefficient sequence of
$q_{\u}$.
Note that the right-hand side depends on the non-truncated signature $\widehat{\mathbb{Z}}$, but just through the terms that can be expressed as a polynomial in $\widehat{\mathbb{Z}}^{\leq N}$.
In particular, the relevant terms of higher level than $N$ of $\widehat{\mathbb{Z}}_0=:\widehat{\z}_0^{\infty}$ correspond to powers of the lower ones and are thus uniquely determined by them.

Let $L$ be the operator associated to $\widehat{\mathbb{X}}^{e_{\emptyset}}$ as defined in \eqref{eq:Lu} and $\overline L$ its extension to the full signature as of Definition \ref{def:Lbar}.
Provided that $\v(\u)\in \mathcal D(\overline L)$
    and that the corresponding integrability
conditions hold and since
$\widehat{\mathbb{Z}}_t^{\leq N}={\widehat\z_0} \otimes \widehat{\mathbb{X}}_t^{{\widehat\z_0}}$, 
Corollary~\ref{Corsupernew}
thus yields
\[
\E[\langle \v(\u), \widehat{\mathbb{Z}}_t \rangle \, |\,  \widehat{\mathbb{Z}}_0=\widehat{\z}_0^{\infty}]= \langle \c_{\u}(t), \widehat{\z}_0^{\infty}\rangle, 
\]
where
\begin{equation}
\label{eq:hol-coeff-ode}
\frac{d}{ds}\c_{\u}(s)
=
\overline L\,\c_{\u}(s),
\qquad
\c_{\u}(0)=\v(\u),
\end{equation}
under the condition that \eqref{eq:hol-coeff-ode} admits a solution.
The inner expectation of the loss  \eqref{eq:tractable-hol-loss} then becomes
$
\int_0^T
\lambda(t)\,
\langle \c_{\u}(t), 
\widehat{\z}_0^{\infty}\rangle\,dt.
$
If $\lambda(t)$ is an entire  function
the above reasoning can be applied to 
\begin{align}\label{tildev}
\langle \tilde{\v}(\u), \widehat{\mathbb{Z}}_T \rangle := \int_0^T \lambda(t) q_{\u}(t,\mathbb{Z}^{\leq N}_t)
dt
\end{align}
so that not even an integral has to be computed. Again, provided that $\tilde{\v}(\u)\in  \mathcal D(\overline L)$ and the corresponding integrability
conditions hold, 
$
\E[\langle \tilde{\v}(\u), \widehat{\mathbb{Z}}_t \rangle \, |\,  \widehat{\mathbb{Z}}_0=\widehat{\z}^{\infty}_0]= \langle \tilde{\c}_{\u}(t), \widehat{\z}_0^{\infty}\rangle
$
where $\tilde{\c}_{\u}(t)$ is the solution of 
\eqref{eq:hol-coeff-ode}
with $\tilde{\c}_{\u}(0)=\tilde \v(\u)$. Then \eqref
{eq:tractable-hol-loss}
becomes
$$
\operatorname{argmin}_{\u}
 \E_{\z_0\sim p_{\mathrm{data}}}\left[\langle \tilde{\c}_{\u}(T), \widehat{\z}_0^{\infty}\rangle\right].
$$
 Thus, generalized score matching is reduced to the
solution of a linear coefficient equation together with a deterministic
optimization over the coefficients $\u$. 
We now summarize this  in the following theorem.

\begin{theorem}\label{th: score opti}
Let $\lambda$ in \eqref{eq:hol-denoising-loss} be an entire function.
Moreover, let $L$ be the operator associated to $\widehat{\mathbb{X}}^{e_{\emptyset}}$ as defined in \eqref{eq:Lu} and $\overline L$ its extension to the full signature as of Definition \ref{def:Lbar}. Assume that all the conditions of Corollary~\ref{Corsupernew} with $\u$ replaced by $\tilde{\v}(\u)$ defined in \eqref{tildev} are satisfied. 
Then, \eqref{eq:hol-denoising-loss} is equivalent to 
\[
\operatorname{argmin}_{\u}\mathbb{E}_{\z_0\sim p_{\mathrm{data}}}[\langle \tilde{\v}(\u), \widehat{\Z}_T \rangle] =\operatorname{argmin}_{\u}
 \mathbb{E}_{\z_0\sim p_{\mathrm{data}}}\left[\langle \bm \tilde{c}_{\u}(´T), \widehat{\z}^{\infty}_0\rangle\right],
\]
with
\begin{equation}
\label{eq:hol-coeff-odetilde}
\frac{d}{ds}\tilde{\c}_{\u}(s)
=
\overline L \tilde{\c}_{\u}(s),
\qquad
\tilde{\c}_{\u}(0)=\tilde{\v}(\u).
\end{equation}
\end{theorem}

Truncating the operator $\overline L$, $\widehat{\z}_0^{\infty}$ and $\tilde{\v}(\u)$ at some level $K$ and denoting them by $\overline L^K$, $\widehat{\z}_0^K$ and  $\tilde{\v}^K(\u)$, the 
\eqref{eq:hol-coeff-odetilde} becomes  
$$
\operatorname{argmin}_{\u}
(\tilde{\v}^K(\u))^{\top} \exp((\overline L^K )^\top T)\E_{\z_0\sim p_{\mathrm{data}}}\left[ \widehat{\z}_0^K\right],
$$
meaning that only one deterministic optimization and one computation of a matrix exponential have to be performed.

Once the score function is learned by this method, we can use $g_{\u}(t, \x)$ as an approximation of the generalized score function in the reverse process $\widehat{\mathbb{Y}}^{\leq N}$ to generate signature samples whose distribution resembles the empirical law  $p_{\mathrm{data}}$.

\begin{remark} \label{rem:appRd}
\begin{enumerate}
\item 
Note that the entire sig-jump-diffusion framework is tailored to generate samples in $G^N(\mathbb{R}^d)$ but it can of course also be used to generate new samples of random variables with values in  $\mathbb{R}^d$. To this end, the initial value $\g$ in \eqref{entsig2} just has to
be chosen such that 
$\g^{(1)}$ corresponds   to the data samples $z_0 \in \mathbb{R}^d$. By Remark 
\ref{rem:connection} (i), the marginal distributions of the first level of the reverse process $\mathbb{Y}^{\leq N}$ then correspond to those of $X^{\g}$ (with $\g^{(1)}=z_0$) given by \eqref{entsig},
and  yield the new $\mathbb{R}^d$-valued samples.  
\item Observe furthermore that for any polynomial or holomorphic jump-diffusion  in the sense of \cite{CKT:12, FL:20, CPS:24holo} with arbitrary state space, the above learning procedure can be applied by approximating the generalized score function either by a polynomial or a holomorphic function.

Let us exemplify this in the case of a $d$-dimensional polynomial diffusion $X$ with $X_0=\z^{(1)}=z_0$. Note that every appearance of $X_t$ in the coefficients can be replaced by $(\z_0\otimes \X_t)^{(1)}$, showing that $X$ falls in the entire sig-jump-diffusion framework. 
Denoting then its  marginal density by $p_t(x)$, we aim to learn the score $(t,x) \mapsto \nabla \log p_t(x)$ via a polynomial $(t,x) \mapsto g_{\u}(t,x)$ with coefficients $\u$ by minimizing 
the following loss
\[
\u^* \in \operatorname{argmin}_{\u} \int_0^T
\lambda(t) \mathbb{E}_{z_0 \sim p_{\mathrm{data}}} \left[\mathbb{E}[ \| g_{\u}(t,X_t) - \nabla_{x} \log p_t(X_t|z_0) \|^2 | z_0] \right] dt,
\]
where $\lambda$ denotes some weight function and $p_t(x|z_0)$ the  conditional density given $z_0$.
This is equivalent to 
\[
\u^* \in \operatorname{argmin}_{\u} \int_0^T
\lambda(t) \mathbb{E}_{z_0 \sim p_{\mathrm{data}}} [
\mathbb{E}[  \sum_{i=1}^{d}\bigl((g_{\u,i}(t,X_t))^2 - 2g_{\u,i}(t,X_t)\partial_{x_i}\log p_t(X_t|z_0)\bigr)  | z_0]] dt.
\]
Using
 integration by parts and the usual boundary assumptions on the density we get
\begin{align*}
\mathbb{E}\left[   -2 g_{\u,i}(t, X_t)\partial_{X_{t,i}} \log p_t(X_t|z_0)  | z_0\right]
&= \int  -2g_{\u,i}(t, x) \partial_{x_i} \log p_t(x|z_0) p_t(x|z_0)dx\\
&=\int  -2g_{\u,i}(t, x) \partial_{x_i}  p_t(x|z_0) dx\\
&=
\int   2\partial_{x_i} g_{\u,i}(t, x)   p_t(x|z_0) dx
\end{align*}
so that the inner expectation becomes
$
\mathbb{E}[  \sum_{i=1}^{d}\bigl((g_{\u,i}(t,X_t))^2 + 2\partial_{x_i}g_{\u,i}(t,X_t)\bigr) | z_0],
$
which is computable by the standard polynomial moment formula.

\end{enumerate}
\end{remark}

\subsection{Invariant measures for the Markovian lift of entire sig-diffusions} 
\label{sec:invariant}

To illustrate the practical applicability of the above approach, we here give an example of an entire sig-diffusion whose lift admits an invariant measure on $G^N(\mathbb{R}^d)$.

Let $m_G$ be a Haar measure\footnote{For an explicit construction through the corresponding free step-$N$
nilpotent Lie algebra 
\cite[Corollary~9.4.21]{LeDonne2025}.} on $G^N(\mathbb R^d)$, i.e., the unique (up to a positive multiplicative factor) nonzero, locally finite measure on \(G^N(\R^d)\)
such that, for every measurable set \(A\) and every
\({\bm h}\in G^N(\R^d)\),
\[
m_G({\bm h}\otimes A)=m_G(A\otimes{\bm h})=m_G(A).
\]
 The following theorem now shows a  construction of  entire sig-diffusions with truncated signature-dependent coefficients 
in spirit of Langevin type dynamics
 that admit an invariant density. 
As we show in the proof its lifted sig-diffusion $\mathbb{Z}^{\leq N}$  belongs to the class of hypoelliptic gradient diffusions on Lie groups studied in \cite{BaudoinHairerTeichmann2008}. Indeed, the free nilpotent group $G^N(\mathbb R^d)$ is a finite-dimensional Lie group, and the vector fields associated with the first-level signature components are left-invariant and satisfy Hörmander's bracket-generating condition. The corresponding generator is therefore of the form considered in \cite{BaudoinHairerTeichmann2008}, where global existence and invariance of a Gibbs-type measure with a density with respect to the Haar measure are established. Moreover, the strong Feller property and irreducibility of the lifted diffusion yield uniqueness of the invariant measure and convergence in total variation by Doob's ergodic theorem, which we in apply in the version provided by \cite{GerlachNittka2012}.

\begin{theorem}\label{thmInvariant}
Let $U:G^N(\mathbb R^d)\to[0,\infty)$ be the restriction of an entire function on $T^N(\R^d)$
such that $
        \int_{G^N(\mathbb R^d)}
        e^{-U(\g)}\, {m_G}(d\g)<\infty.
   $
With the notation introduced for \eqref{entsig}, define $\overline\sigma=\operatorname{id}$, $\overline \delta=0$, $\overline{b}: G^N(\mathbb{R}^d) \to \mathbb{R}^d$ via
$$
    \overline{b}_i(\g)
    =
    -\frac12
    \left.
    \frac{d}{ds}
    U\bigl(\g\otimes\exp_N(s\epsilon_i)\bigr)
    \right|_{s=0},
$$
and consider for every $\z_0 \in G^N(\mathbb{R}^d)$ the corresponding  entire sig-diffusion with truncated signature-dependent coefficients:
\begin{equation}
\begin{aligned}
    dX_t^{\z_0}
    &=
    \overline{b}\left(
        \z_0 \otimes \mathbb{X}^{\z_0}_{t}
    \right)dt+dW_t, \quad 
    X_0^{\z_0}=\z_0^{(1)},
\end{aligned}
\label{eq:original_signature_sde}
\end{equation}
where $W$ is a $d$-dimensional Brownian motion, and the lifted equation
\begin{equation}
    \mathbb Z_t^{\leq N}
    =
    \z_0\otimes \mathbb{X}^{\z_0}_{t}.
\label{eq:signature_lift_definition}
\end{equation}
Then 
$\mathbb{Z}^{\leq N}$ admits a unique global strong solution  and  
    \begin{equation}
        \mu(d\g)= \frac{e^{-U(\g)}\, {m_G}(d\g)}{\int_{G^N(\mathbb R^d)}
        e^{-U(\g)}\, {m_G}(d\g)}
\label{eq:invariant_measure}
    \end{equation}
is the unique invariant probability measure 
such that 
\[
\lim_{t \to \infty}\| p_t( \cdot | \z_0)- \mu\|_{\mathrm{TV}} = 0, \quad \text{for all } \z_0 \in G^N(\mathbb{R}^d),\]
where $p_t( \cdot | \z_0)$ denotes the transition density of $\mathbb{Z}_t^{\leq N}$ given $\z_0$ and $\mathrm{TV}$ the total variation norm.
In particular, the limiting distribution
of $\mathbb Z_t^{\leq N}$ is independent of its
initial condition $\z_0$.

\end{theorem}

\begin{proof}
We apply Proposition 3.1 of \cite{BaudoinHairerTeichmann2008}  to show  global strong existence and uniqueness for the lifted diffusion $\mathbb{Z}^{\leq N}$.

Denote by \(V_1,\ldots,V_d\) the left-invariant vector fields on \(G^N(\mathbb{R}^d)\) 
given by 
\[
V_i(\g)= \left.\frac{d}{ds}\bigl(\mathbf g\otimes\exp_N(s\epsilon_i)\bigr)
\right|_{s=0}= \, \g \otimes \epsilon_i.
\]
Consequently, when acting on a smooth function \(f:G^N(\mathbb R^d)\to\mathbb R\), we have
\[
V_i f(\mathbf g)
=
\left.\frac{d}{ds}
f(\mathbf g\otimes\exp_N(s\epsilon_i))
\right|_{s=0}.
\]
The lifted diffusion $\mathbb{Z}^{\leq N}$ can now be written in terms of these vector fields, i.e.,
\begin{align} 
d\mathbb Z_t^{\leq N}&=\sum_{i=1}^d
    \overline b_i(\mathbb Z_t^{\leq N})
    (\mathbb Z_t^{\leq N}\otimes\epsilon_i)\, dt+\sum_{i=1}^dV_i(\mathbb Z_t^{\leq N})\circ dW_t^i, \notag\\
    &=-\frac12\sum_{i=1}^d
(V_iU)(\mathbb Z_t^{\leq N})
V_i(\mathbb Z_t^{\leq N})\,dt
+
\sum_{i=1}^d
V_i(\mathbb Z_t^{\leq N})\circ dW_t^i. \notag
\end{align}
and its generator reads
\begin{align}
\label{eq:generator_potential}
&\mathcal L^U f(\g)=\frac12\sum_{i=1}^dV_i^2f(\g)
-\frac12\sum_{i=1}^d(V_iU)V_if(\g)\\
&\quad =
\frac12\sum_{i=1}^d \frac{d^2}{ds}f \left.(\mathbf g\otimes\exp_N(s\epsilon_i))
\right|_{s=0}
-\frac12\sum_{i=1}^d  \frac{d}{ds}U\left.(\mathbf g\otimes\exp_N(s\epsilon_i))
\right|_{s=0}\frac{d}{ds}f \left.(\mathbf g\otimes\exp_N(s\epsilon_i))
\right|_{s=0}.\notag
\end{align}
Setting
$
\mathcal L=\frac12\sum_{i=1}^dV_i^2$ and $
\Gamma(f,g)=\sum_{i=1}^d(V_if)(V_ig)
$, we have
$
\mathcal L^U=\mathcal L-\frac12\Gamma(U,\cdot).
$
From this, we see that the lifted diffusion $\mathbb{Z}^{\leq N}$ is precisely of the form considered in Proposition 3.1 of \cite{BaudoinHairerTeichmann2008}. Note that the conditions on $U$ are also satisfied, since it is  entire and $
\int_{G^N(\mathbb{R}^d)} e^{-U(\g)}\,m_G(d\g )<\infty
$ 
by assumption.
Thus it remains to verify that the Hörmander conditions of the diffusion vector fields $V_1, \ldots, V_d$ are satisfied.
To this end note that the Lie bracket of two such vector fields satisfies
\begin{align*}
[V_i,V_j](\g)
&=
\g\otimes\epsilon_i\otimes\epsilon_j
-
\g\otimes\epsilon_j\otimes\epsilon_i=
\g\otimes[\epsilon_i,\epsilon_j].
\end{align*}
By induction, every iterated Lie bracket of
$V_1,\ldots,V_d$ is obtained by left multiplication
of the corresponding iterated Lie bracket of
$\epsilon_1,\ldots,\epsilon_d$, i.e.,
\begin{equation*}
[V_{i_1},[V_{i_2},\ldots,
[V_{i_{k-1}},V_{i_k}]\ldots]](\g)
=
\g\otimes
[\epsilon_{i_1},[\epsilon_{i_2},\ldots,
[\epsilon_{i_{k-1}},\epsilon_{i_k}]\ldots]].
\end{equation*}
Let now $g^N(\mathbb R^d)$ be the free
step-$N$ nilpotent Lie algebra generated by
$\epsilon_1,\ldots,\epsilon_d$, which in particular implies that these iterated
Lie brackets of lengths at most $N$ span
$g^N(\mathbb R^d)$.
Moreover, left multiplication by $\g$ identifies
the Lie algebra with the tangent space at $\g$,
that is,
$
T_\g G^N(\mathbb R^d)
=
\g\otimes g^N(\mathbb R^d).
$
Consequently, the strong H\"ormander bracket condition of the diffusion vector fields
\[
\operatorname{span}
\left\{
V_i(\g),\,
[V_i,V_j](\g),\,
[V_i,[V_j,V_k]](\g),\ldots
\right\}
=
T_\g G^N(\mathbb R^d)
\]
holds for every $\g\in G^N(\mathbb R^d)$,
where brackets of length at most $N$ are sufficient.
Therefore all assumptions of Proposition 3.1 of \cite{BaudoinHairerTeichmann2008} are satisfied, which together with the Lipschitz property of the vector fields yields a   unique global strong solution for $\mathbb{Z}^{\leq N}$ for every initial condition \(\mathbf z_0\in G^N(\mathbb{R}^d)\).

Concerning the invariant measure,
we here give the full argument even though the first part, namely that  $\mu$ is an invariant measure can be deduced from Section 2 in \cite{BaudoinHairerTeichmann2008}.
We start by showing that the generator $\mathcal  L^U$ is formally symmetric
in $L^2(\mu)$.
Let $f,\varphi\in C_c^\infty(G^N(\mathbb R^d))$.
Recall that the Haar measure
$ {m_G}$ is invariant under right $\otimes$ multiplications.
Consequently, for every
$F\in C_c^\infty(G^N(\mathbb R^d))$ and
$i=1,\ldots,d$,
\begin{equation}
\int_{G^N(\mathbb R^d)}
\left.
\frac{d}{ds}
F\bigl(\g\otimes\exp_N(s\epsilon_i)\bigr)
\right|_{s=0}
\,d {m_G}(\g)=0.
\label{eq:haar_ibp}
\end{equation}

For fixed $i$, define
\[
F_i(\g)
=
f(\g)e^{-U(\g)}
\left.
\frac{d}{ds}
\varphi\bigl(\g\otimes\exp_N(s\epsilon_i)\bigr)
\right|_{s=0}.
\]
Since $f$ has compact support, we have
$F_i\in C_c^\infty(G^N(\mathbb R^d))$. 
Applying \eqref{eq:haar_ibp} to $F_i$ gives
\begin{align*}
0
&=
\int_{G^N(\mathbb R^d)}
e^{-U(\g)}
\Bigg[
\left.
\frac{d}{ds}
f\bigl(\g\otimes\exp_N(s\epsilon_i)\bigr)
\right|_{s=0}\left.
\frac{d}{ds}
\varphi\bigl(\g\otimes\exp_N(s\epsilon_i)\bigr)
\right|_{s=0}
\\
&\qquad+
f(\g)
\left.
\frac{d^2}{ds^2}
\varphi\bigl(\g\otimes\exp_N(s\epsilon_i)\bigr)
\right|_{s=0}\\&\qquad
-
f(\g)
\left.
\frac{d}{ds}
U\bigl(\g\otimes\exp_N(s\epsilon_i)\bigr)
\right|_{s=0}
\left.
\frac{d}{ds}
\varphi\bigl(\g\otimes\exp_N(s\epsilon_i)\bigr)
\right|_{s=0}
\Bigg]\,d {m_G}(\g).
\end{align*}

Summing over $i=1,\ldots,d$, dividing  by
${2\int_{G^N(\mathbb R^d)}
        e^{-U(\g)}\, {m_G}(d\g)}$, and recalling \eqref{eq:invariant_measure}
we obtain, by the expression
\eqref{eq:generator_potential} for the generator,
\begin{align*}
&\int_{G^N(\mathbb R^d)}
f(\g)\mathcal L^U \varphi(\g)\,\mu(d\g)
\\
&\quad=
-\frac12\sum_{i=1}^d
\int_{G^N(\mathbb R^d)}
\left.
\frac{d}{ds}
f\bigl(\g\otimes\exp_N(s\epsilon_i)\bigr)
\right|_{s=0}
\left.
\frac{d}{ds}
\varphi\bigl(\g\otimes\exp_N(s\epsilon_i)\bigr)
\right|_{s=0}
\,\mu(d\g).
\end{align*}
Since the right-hand side is symmetric in
$f$ and $\varphi$, it follows that
\[
\mathcal E(f,\varphi):=-\int_{G^N(\mathbb R^d)}
f(\g)\mathcal L^U \varphi(\g)\,\mu(d\g)
=
-\int_{G^N(\mathbb R^d)}
\varphi(\g)\mathcal L^U f(\g)\,\mu(d\g).
\]
Thus, $\mathcal L^U $ is formally symmetric in $L^2(\mu)$.
Since $\mu$ has a smooth, strictly positive
density with respect to Haar measure, the bilinear form $
    \mathcal E(f,\varphi)$
is closable in $L^2(\mu)$,
and its closure is a regular, strongly
local Dirichlet form.
By local uniqueness of the martingale problem,
the diffusion associated to this Dirichlet form coincides with
$\mathbb Z^{\leq N}$ up to its lifetime.
Since we have already established the existence of a global solution, it is thus
conservative. Consequently, its
transition semigroup $(P_t)_{t\geq0}$
is symmetric in $L^2(\mu)$ and satisfies
$P_t1=1$ (which is also proved in Proposition 3.1 of \cite{BaudoinHairerTeichmann2008}). Hence
\begin{align*}
    \int P_tf\,d\mu
    &=
    \langle P_tf,1\rangle_{L^2(\mu)}=
    \langle f,P_t1\rangle_{L^2(\mu)}=
    \int f\,d\mu.
\end{align*}
Consequently, $\mu$
is an invariant probability measure for $\mathbb{Z}^{\leq N}$.

We now prove that $\mu$ is unique and that
\[
    \lim_{t\to\infty}
    \|p_t(\cdot |\z_0)-\mu\|_{\mathrm{TV}}=0,
\]
for every $\z_0\in G^N(\mathbb R^d)$. To this end we shall apply a version of Doob's ergodic theorem \cite{Doob1948} as proved in \cite{GerlachNittka2012}. To verify the conditions of Proposition~4.3(ii) therein, we first show the strong Feller property of the transition semigroup $(P_t)_{t \geq 0}$ i.e., that it 
 maps bounded measurable functions on $G^N(\mathbb{R}^d)$ to bounded continuous functions on  $G^N(\mathbb{R}^d)$. 
 Since the diffusion vector fields satisfy the strong
H\"ormander bracket condition as proved above, 
Theorem~7.1 of 
\cite{ArnaudonThalmaier2010} (see also p.~1286 therein)
then implies that the transition semigroup
is strong Feller.

Proposition 4.3(ii) in \cite{GerlachNittka2012}
now allows us to conclude $L^1(\mu)$ irreducibility in the sense of 
Definition 2.4 therein\footnote{
For a conservative Markov semigroup $(P_t)_{t\geq 0}$
with invariant probability measure $\mu$,
irreducibility on $L^1(\mu)$ is equivalent to the
following property: whenever a measurable set $A$
satisfies
$
    P_t\mathbf{1}_A=\mathbf{1}_A, \mu\text{-a.e.}$  for all $t\geq 0$, 
then necessarily
$  \mu(A)=0 $ or  $\mu(A)=1$.
}. First, $G^N(\mathbb R^d)$ is a connected,
finite-dimensional Lie group and hence a
second-countable Polish space.
Moreover, the invariant probability measure
$\mu(d\g)$
is strictly positive in the sense that
$\mu(O)>0$
for every nonempty open subset
$O\subset G^N(\mathbb R^d)$, since  $e^{-U(\g)}>0$ for every $\g$, and
the Haar measure has full support.
Since $(P_t)_{t \geq 0}$ is conservative, 
$ P_t1=1$ for all  $t\geq0$.
Thus, the constant function $1$ is a strictly
positive fixed point of $(P_t)_{t\geq0}$.
Consequently, all assumptions of
\cite[Proposition 4.3(ii)]{GerlachNittka2012}
are satisfied and  it follows that the extension
of $(P_t)_{t\geq0}$ to $L^1(\mu)$ is irreducible.

Uniqueness of the invariant probability
measure and convergence to equilibrium now follow from
Theorems~4.6
and 4.10 in
\cite{GerlachNittka2012}.
Note that the semigroup is stochastically
continuous because $\mathbb{Z}^{\leq N}$ has
continuous sample paths.
In particular, for every $f\in C_b(G^N(\mathbb R^d))$,
\[
    \lim_{t\downarrow0}P_tf(\g)=f(\g)
\]
for every $\g\in G^N(\mathbb R^d)$. By dominated convergence,
\[
    \lim_{t\downarrow0}
    \int |P_tf-f|\,d\mu=0.
\]
Using the density of $C_b(G^N(\mathbb R^d))$
in $L^1(\mu)$ and the fact that $(P_t)_{t\geq 0}$
acts as a contraction semigroup on $L^1(\mu)$,
the same convergence holds for every
$f\in L^1(\mu)$.
Thus, the continuity assumption of
Theorem 4.6 in \cite {GerlachNittka2012}
is also satisfied.
Moreover, its transition probabilities are Borel probability measures so that in the notation of \cite{GerlachNittka2012} $P_t'\delta_{\z_0}=p_t(\cdot| \z_0)$ is a Borel measure.
With $e\equiv 1$, 
Theorem 4.6 in \cite{GerlachNittka2012} 
therefore yields
\[
    \lim_{t\to\infty}P_tf(\g)
    =
    \int_{G^N(\mathbb R^d)}f(\h)\,\mu(d\h)
\]
for every bounded measurable function $f$, 
uniformly on compact sets of  $G^N(\mathbb R^d)$. 
Finally, by
Theorem 4.10 in \cite{GerlachNittka2012},
this convergence can be strengthened
to convergence in total variation:
\[
    \lim_{t\to\infty}
    \left\|
        p_t(\cdot | \z_0)-\mu
    \right\|_{\text{TV}}=0
\]
for every $\z_0\in G^N(\mathbb R^d)$.
In particular, $\mu$ is the unique invariant
probability measure, and the limiting distribution
of $\mathbb Z_t^{\leq N}$ is independent of its
initial condition $\z_0$.
\end{proof}

In the following, we provide an example of a quadratic potential $U$ that gives rise to a Gaussian-type Gibbs measure on $G^N(\mathbb{R}^d)$ as the invariant measure.

\begin{example}
Fix \(\z_0,{\bm m}\in G^N(\R^d)\) and set $U(\g)=\frac12\|\g-{\bm m}\|^2$. Observe that the corresponding  entire sig-diffusion \eqref{eq:original_signature_sde} reads
\[
dX_t^{\z_0,i}=-\frac12
\langle
\z_0\otimes\X_t^{\z_0}-{\bm m},\,
\pi_{\leq N}\bigl((\z_0\otimes\X_t^{\z_0})\otimes\e_i\bigr)
\rangle dt+dW_t^i,
\qquad i=1,\ldots,d,
\]
where $\X^{\z_0}$ denotes the signature of $X^{\z_0}$ truncated at level $N$. Then, by Theorem~\ref{thmInvariant}
 the process \(\z_0\otimes\X_t^{\z_0}\) admits the unique invariant law $\mu$ given by \eqref{eq:invariant_measure} which reads 
\[
\mu(d{\g})
=\frac{e^{-\lVert{\g}-{\bm m}\rVert^2/2}\,
  m_G(d{\g})}{\int_{G^N(\R^d)}
e^{-\lVert{\g}-{\bm m}\rVert^2/2}\,
m_G(d{\g})}.
\]
\end{example}

\begin{remark}
 It is important to consider the Markovian lifts of
genuine sig-(jump)-diffusions (with, in particular, coefficients depending on the whole truncated signature) rather than restricting
attention to lifts of classical Markovian ones. 

Indeed, consider a classical Markovian SDE
\[
    dX_t=b(X_t)\,dt+\sigma(X_t)\,dW_t,
    \qquad X_0=\z_0^{(1)},
\]
with $\z_0\in G^N(\mathbb R^d)$,
and note that its dynamics can equivalently be written as
\begin{align}\label{eq:underSDE}
     dX_t^{\z_0}=b((\z_0\otimes \X_t^{\z_0})^{(1)})\,dt+\sigma((\z_0\otimes \X_t^{\z_0})^{(1)})\,dW_t,
    \qquad X_0^{\z_0}=\z_0^{(1)}.
\end{align}
Consider then its signature lift
$
    \mathbb Z_t^{\z_0}
    :=
    \z_0\otimes \mathbb{X}_t^{\z_0},
$ truncated at level $N$. Recall that we omit $\leq N$ for readability.
Note that for every $\g\in  G^N(\mathbb R^d)$ with $\g^{(1)}=0$
we get
$
    \mathbb Z_t^{\g\otimes\z_0}
    =
    \g\otimes\mathbb Z_t^{\z_0}
$.
This follows because
$
    (\g\otimes\z_0)^{(1)}=\z_0^{(1)},
$
so that the underlying Markovian SDE \eqref{eq:underSDE} has
the same initial condition in both cases.

This implies that if the law of $\Z_t^{\z_0}$ converges in law to an invariant probability distribution $\mu$, then the law of $\Z_t^{\g\otimes\z_0}$ converges in law to
$$\mu^\g(A):=\int  1_{\{\g\otimes\xi\in A\}}\mu(d\xi)$$
for all measurable sets $A \subseteq G^N(\mathbb{R}^d)$.
Since the noncompactness of  $G^N(\R^d)$ prevents $\mu^\g$ from corresponding to $\mu$ for each $\g$, we can conclude that the law of $\Z_t^{\g\otimes\z_0}$ depends on $\g$ also for  $t$ converging to infinity. Hence,  forgetting the initial condition, which is necessary for generative diffusion modeling, is in this case not possible.

In contrast, signature-dependent coefficients
generally break this equivariance, allowing
for mean reversion in higher signature levels
and the existence of a unique invariant probability
measure independent of the initial condition, as demonstrated by the preceding theorem.

If one wishes 
to consider signature-like lifts of classical Markovian SDE, an alternative is to work with the process $\widetilde \Z^{\z_0}$ defined by
$$\widetilde \Z_t^{\z_0}=D^\lambda_t(\z_0\otimes \widetilde \X_t^{\z_0}),$$
where $\widetilde \X^{\z_0}$ is the signature of $(\int_0^te^{\lambda u}dX_u^{\z_0})_{t\in [0,T]}$ and $D^\lambda_t$ is the dilation operator defined by $\langle \e_I,D^\lambda_t\u\rangle=e^{-\lambda t|I|}\u_I$. Replacing $\z_0$ by the truncated signature 
$\z_0^\lambda$ of the path
$$\bigg(\int_0^r e^{-\lambda(S-u)}\,dX_{-S+u}\bigg)_{r\in[0,S]},$$ the process $\widetilde \Z^{\z_0^\lambda}$ corresponds to the exponentially fading memory signature of $(X_t)_{t\in[-S,T]}$ studied in \cite{Abi:Fading}, where the relevant ergodic properties are also established for $X_t=(t,W_t)$. Since $D^\lambda_t$ is a deterministic linear operator on $T^N(\R^d)$ that preserves $G^N(\R^d)$, the analytical methods developed in this section can still be applied.
\end{remark}

Next, we turn to the computation of  the generalized score function of $\mathbb{Z}^{\leq N}$ in this specific setting.
 Consider $\mathbb{Z}^{\leq N}$ given by \eqref{eq:signature_lift_definition}.
Our aim here is to identify the generalized score function
of $\mathbb{Z}^{\leq N}$ in the sense of Definition~\ref{def:sig-score}, that is in the case 
when $\mathcal{K}^N$ is 
the second-order operator 
defined with respect
to the ambient truncated tensor algebra $T^N(\mathbb{R}^d)$.

Recall that the dynamics of $\mathbb Z_t^{\leq N}$ are given by
$$
\begin{aligned}
    d\mathbb Z_t^{\leq N}
    &=
    \bigg(
        \mathbb Z_t^{\leq N}
        \otimes\overline b(\mathbb Z_t^{\leq N})
        +
        \frac12\sum_{i=1}^d
        \mathbb Z_t^{\leq N}
        \otimes\epsilon_i^{\otimes 2}
    \bigg)dt +
    \sum_{i=1}^d
    \mathbb Z_t^{\leq N}
    \otimes\epsilon_i\,dW_t^i.
\end{aligned}
$$
Hence, its drift and diffusion
coefficients are
$$
    b(\x)
    =
    \x\otimes\overline b(\x)
    +
    \frac12\sum_{i=1}^d
    \x\otimes\epsilon_i^{\otimes 2}
    \qquad\text{and}\qquad
    \sigma_i(\x)
    =
    \x\otimes\epsilon_i,
$$
respectively. Consequently, the generator decomposes as
$
    Lf(\x)
    =
    \left\langle
        b(\x),\nabla f(\x)
    \right\rangle
    +
    \mathcal K^N f(\x),
$
for
$$
    \mathcal K^N f(\x)
    =
    \frac12\sum_{i=1}^d
    \nabla^2 f(\x)
    \left[
        \x\otimes\epsilon_i,
        \x\otimes\epsilon_i
    \right].
$$
Here $\nabla f$ and $\nabla^2f$ denote the
gradient and Hessian in the ambient tensor algebra
$T^N(\mathbb R^d)$ and in particular
$$\nabla^2 f(\x)[{\bf h},\g]=\sum_{|I|,|J|\leq N}\partial^2_{\x_I\x_J} f(\x){\bf h}_I\g_J.$$

\begin{remark}
It is important to distinguish
$\mathcal K^N$ from the second-order operator
defined by differentiation along the signature group. 
Indeed, for fixed $\g\in G^N(\mathbb R^d)$
and $i=1,\ldots,d$, consider the curve
\[
    h_i(s)
    =
    \g\otimes\exp_N(s\epsilon_i).
\]

We have
\begin{equation}
    h_i'(0)=\g\otimes\epsilon_i,
    \qquad
    h_i''(0)=\g\otimes\epsilon_i^{\otimes 2}.
\end{equation}

By the second-order chain rule,
$$
\left.
    \frac{d^2}{ds^2}
    f\bigl(\g\otimes\exp_N(s\epsilon_i)\bigr)
    \right|_{s=0}=
    \nabla^2 f(\g)
    \left[
        \g\otimes\epsilon_i,
        \g\otimes\epsilon_i
    \right]
    +
    \left\langle
        \nabla f(\g),
        \g\otimes\epsilon_i^{\otimes 2}
    \right\rangle.
$$
Consequently,
\begin{equation}
\begin{aligned}
    \mathcal K^N f(\g)
    &=
    \frac12\sum_{i=1}^d
    \left.
    \frac{d^2}{ds^2}
    f\bigl(\g\otimes\exp_N(s\epsilon_i)\bigr)
    \right|_{s=0}-
    \frac12\sum_{i=1}^d
    \left\langle
        \nabla f(\g),
        \g\otimes\epsilon_i^{\otimes 2}
    \right\rangle.
\end{aligned}
\label{eq:score_operator_identity}
\end{equation}
Thus, the intrinsic second-order operator
and the second-order part of the ambient Itô
generator differ by a first-order correction term.
\end{remark}

Let $
    P_t^N
    =
    \mathcal L(\mathbb Z_t^{\leq N})$
denote the law of the lifted process at time $t$.
Since the diffusion vector fields  $D_1,\ldots,D_d$ satisfy Hörmander's
bracket-generating condition as proved in Theorem \ref{thmInvariant}, the transition law
admits a smooth density with respect to Haar measure
due to the standard hypoellipticity assumptions. 
Moreover, they
define a bracket-generating sub-Riemannian structure on $G^N(\mathbb R^d)$. Furthermore,  $G^N(\mathbb R^d)$   is connected and complete
with respect to the associated Carnot--Carath\'eodory distance.
Since the generator can be written as
\[
Lf=\frac12 e^U\sum_{i=1}^d D_i(e^{-U}D_if),
\]
it is a weighted sub-Laplacian associated
with the smooth, strictly positive invariant measure
$\mu$. By the strict positivity of heat kernels for such
sub-Laplacians on connected, complete sub-Riemannian
manifolds, \cite{
Strichartz1986} and  
\cite[Section~2.1, p.~716]{BarilariBoscainNeel2019}),
the transition density satisfies
\[
p_t(\g | \z_0)>0,
\qquad
\forall\,t>0,\quad
\z_0,\g\in G^N(\mathbb R^d).
\]
Hence, we can write
\begin{equation*}
    P_t^N(d\g)=p_t(\g)\, {m_G}(d\g),
\end{equation*}
with $p_t(\g)=\int_{G^N(\mathbb{R}^d)}p_t( \g|\z_0) p_{\mathrm{data}}(d\z_0)$ smooth  and strictly positive. 
Using this, the following result now establishes the form of the generalized score function of $\mathbb{Z}^{\leq N}$
given by \eqref{eq:signature_lift_definition}.

\begin{proposition}
Let $t>0$ and suppose that the map
$S_t^N:G^N(\mathbb R^d)\to T^N(\mathbb R^d)$,
defined by
\begin{equation}
\begin{aligned}
    S_t^N(\g)
    &=
    \frac12\sum_{i=1}^d
    \Bigg[
    \left.
    \frac{d}{ds}
    \log p_t\bigl(
        \g\otimes\exp_N(s\epsilon_i)
    \bigr)
    \right|_{s=0}
    (\g\otimes\epsilon_i)
    +\g\otimes\epsilon_i^{\otimes 2}
    \Bigg],
\end{aligned}
\label{eq:score_transient}
\end{equation}
is integrable with respect to $P_t^N$.
Then $S_t^N$ is a generalized score function
of $\mathbb Z^{\leq N}$, i.e.,
\begin{equation}
    \int
    \mathcal K^N f(\g)\,P_t^N(d\g)
    =
    -
    \int
    \left\langle
        S_t^N(\g),\nabla f(\g)
    \right\rangle
    P_t^N(d\g)
\label{eq:score_weak_identity}
\end{equation}
for every
$f\in C_c^2(T^N(\mathbb R^d))$.
\end{proposition}

\begin{proof}
By \eqref{eq:score_operator_identity},
\begin{align*}
    \int\mathcal K^Nf\,dP_t^N
    &=
    \frac12\sum_{i=1}^d
    \int
    \left.
    \frac{d^2}{ds^2}
    f\bigl(\g\otimes\exp_N(s\epsilon_i)\bigr)
    \right|_{s=0}
    p_t(\g)\, {m_G}(d\g)
    \\
    &\quad-
    \frac12\sum_{i=1}^d
    \int
    \left\langle
        \nabla f(\g),
        \g\otimes\epsilon_i^{\otimes 2}
    \right\rangle
    p_t(\g)\, {m_G}(d\g).
\end{align*}
Since the Haar measure is invariant under right
$\otimes$ multiplications, integration by parts gives
\begin{align*}
    &\int
    \left.
    \frac{d^2}{ds^2}
    f\bigl(\g\otimes\exp_N(s\epsilon_i)\bigr)
    \right|_{s=0}
    p_t(\g)\, {m_G}(d\g)
    \\
    &\quad=
    -\int
    \left.
    \frac{d}{ds}   f\bigl(\g\otimes\exp_N(s\epsilon_i)\bigr)
    \right|_{s=0} 
    \left.
    \frac{d}{ds}p_t\bigl(\g\otimes\exp_N(s\epsilon_i)\bigr)
    \right|_{s=0}
     {m_G}(d\g).
\end{align*}

Moreover, the first-order chain rule yields
\begin{equation}
    \left.
    \frac{d}{ds}
    f\bigl(\g\otimes\exp_N(s\epsilon_i)\bigr)
    \right|_{s=0}
    =
    \left\langle
        \nabla f(\g),
        \g\otimes\epsilon_i
    \right\rangle.
\end{equation}
Since $p_t>0$, we also have
\begin{align*}
    &\left.
    \frac{d}{ds}
    p_t\bigl(\g\otimes\exp_N(s\epsilon_i)\bigr)
    \right|_{s=0}=
    p_t(\g)
    \left.
    \frac{d}{ds}
    \log p_t\bigl(
        \g\otimes\exp_N(s\epsilon_i)
    \bigr)
    \right|_{s=0}.
\end{align*}
Combining these identities, we obtain
\begin{align*}
    \int\mathcal K^Nf\,dP_t^N
    &=
    -\frac12\sum_{i=1}^d
    \int
    \Bigg\langle
        \nabla f(\g),
    \left.
    \frac{d}{ds}
    \log p_t\bigl(
        \g\otimes\exp_N(s\epsilon_i)
    \bigr)
    \right|_{s=0}
    (\g\otimes\epsilon_i)
    +\g\otimes\epsilon_i^{\otimes 2}
    \Bigg\rangle
    P_t^N(d\g).
\end{align*}

By the definition of the score in \eqref{eq:score_transient},
this is precisely
$
    \int\mathcal K^Nf\,dP_t^N
    =
    -\int
    \langle S_t^N,\nabla f\rangle\,dP_t^N,
$
which proves the claim.
\end{proof}

\begin{remark}
The generalized score in
\eqref{eq:score_transient} is defined
$P_t^N$-almost everywhere on
$G^N(\mathbb R^d)$.
Since $P_t^N$ is supported on the 
group, it can be extended measurably to
$T^N(\mathbb R^d)$ without changing
\eqref{eq:score_weak_identity}.
The integrability condition in the proposition
is an additional assumption and does not follow
immediately from the existence of an invariant measure.
\end{remark}

\subsection{An alternative approach for polynomial diffusions based on inverse signature matching}\label{sec:alternative}

In this section we provide an alternative to score matching based on what we call ``inverse signature matching'', tailored to polynomial diffusion processes
for which we can exploit that the expected signature is obtained by solving a finite-dimensional linear ODE.
Here, we focus on generating new data in line with some law $p_{\textrm{data}}$ on $\mathbb{R}^d$ and not on time-series generation via their signature as in the previous sections.  

The  goal is to use a polynomial diffusion process to perturb the data $z_0 \sim p_{\mathrm{data}}$ on $\mathbb{R}^d$ and to generate new data by the corresponding reverse process. In contrast to the approach provided in Remark~\ref{rem:appRd}~(ii), we here do not directly learn the score by a polynomial or an entire function but, exploiting that, under suitable conditions, the expected signature characterizes the law of a stochastic process \cite{signaturemoments}, we rather match the inverse signature process of the forward polynomial diffusion. 

More precisely, we assume that the forward process is given by 
a polynomial diffusion in the sense of Section~\ref{sec: poly}, i.e.,
\begin{align}\label{eq:SDEforward}
   dX_t=\overline{b}(X_t)dt+\sqrt{\overline{a}(X_t)} dW_t, \quad X_0=z_0\sim p_{\mathrm{data}},
\end{align}
where $\overline b:\mathbb R^d\to\mathbb R^d$ and
$\overline a:\mathbb R^d\to\mathbb S^d_+$ have polynomial components
of degree at most one and two, respectively, and $W$ is a Brownian
motion with respect to a filtration $(\mathcal F_t)$ to which $X$
is adapted. We assume that $\sqrt {\overline a}$ have polynomial components
of degree at most one, and that, for every $t>0$, $X_t$ admits a strictly
positive density $p_t$, with
$p\in C^{1,2}((0,\infty)\times\mathbb R^d)$.
We further assume that $X$ has a unique invariant density $p_\infty$ and that its law converges weakly to $p_\infty(x)dx$ from every deterministic initial state.

Let $Y_t:=X_{T-t}$, $0\leq t\leq T$, be the reversed process.
By Theorem~2.1 in \cite{HS:86}, $Y$ satisfies the following SDE
for $0\leq t<T$:
\begin{align}\label{eq:reversedSDE}
    dY_t =\left(-\overline{b}(Y_t) +\operatorname{div}(\overline{a}(Y_t))
    +\overline{a}(Y_t)\nabla\log p_{T-t}(Y_t)\right)dt+\sqrt{\overline{a}(Y_t)} dB_t, \quad Y_0=X_T\sim p_T,
\end{align}
 where $\operatorname{div}$ is defined in \eqref{eq:divergence}
and $B$ is a Brownian motion with respect to a filtration
$(\hat{\mathcal F}_t)$ to which $Y$ is adapted, on an enlargement
of the probability space if necessary.
Under the stated assumptions, the coefficients of the reverse SDE
are locally Lipschitz in space, uniformly on compact time intervals
contained in $[0,T)$. Consequently, pathwise uniqueness holds
and $Y$ is a strong solution on $[0,T)$.
Moreover, $Y$ extends continuously to $T$, with
$Y_T=X_0\sim p_{\mathrm{data}}$.

Note that the drift of \eqref{eq:reversedSDE} corresponds to choosing $\eta \equiv 1$ in \eqref
{eq:reverse-SDE}
since 
$ \frac{1}{2}(\operatorname{div}(\overline{a})
    +\overline{a}\nabla\log p_{T-t})$ is exactly the generalized score as shown in Remark~\ref{rem:special case}.

As above, the goal here is to learn the map $\overline{a}(\cdot)\nabla\log p_\cdot(\cdot)$ from the available data distributed according to $ p_{\mathrm{data}}$ and then to generate new instances of $z_0$ by sampling $Y_0$ from the invariant measure and computing $Y_T$ by numerically solving the SDE \eqref{eq:reversedSDE}. 
Here the main difference is that, relying on Proposition \ref{prop: inverse sig} below and the results in \cite{signaturemoments}, the unknown term $\overline{a}\nabla\log p_\cdot {(\cdot)}$ is learned  by matching the expected signature of the reverse process to that of the inverse signature of the forward polynomial process, obtained by solving a finite-dimensional linear ODE. The approach is made precise in the following steps.

\vspace{0.2cm}

\begin{enumerate}[label=(\roman*), leftmargin=0pt, itemindent=!, labelsep=0.4em,
                  align=left, itemsep=0.8em, parsep=0pt, topsep=0pt, partopsep=0pt]
    \item \label{it1} \textbf{Signature parametrization of the score.} We assume that each component of
$
        \nabla\log p_\cdot(\cdot ),
   $
    can be approximated with a polynomial (in both variables). Moreover, we  use that
$
         Y_t \approx  \langle \tilde \u(Y_0), \hat{\mathbb{B}}_t\rangle,
$
      for some $\tilde \u(Y_0)\in T^N(\Rset^{d+1})^d$ and some $N\in \Nset$, where $\hat{\mathbb{B}}$ is the Stratonovich signature of the time-extended Brownian motion $\hat B$, whose index coresponding to the time component is denoted by $0$. The proposed approximation is motivated by Lyons' universal limit theorem
\cite{L:98}, which ensures, under suitable
regularity assumptions, that almost surely the solution map
$(\mathbb B_t^{\leq 2})_{t\in[0,T]}\mapsto (Y_t)_{t\in [0,T]}$
is continuous with respect to the appropriate rough-path topologies. Hence,
by the universal approximation theorem for signatures
~\cite{LLN:13,CPS:UATlevy, cuchiero2025signature},  the whole trajectory \((Y_t)_{t\in [0,T]}\) can be uniformly approximated on compact
sets of paths by a linear combination of the time-extended signature of $B$.

   Under these assumptions and as 
   $\overline{a}$ is a polynomial of degree at most $2$,
   we get by the shuffle property
   \begin{align}\label{eq: sigapprox}
        \overline{a}(Y_t)\nabla\log p_{T-t}(Y_t) \approx \langle \u(Y_0), \hat{\mathbb{B}}_t\rangle,
    \end{align}
     for some $\u(Y_0)\in T^N(\Rset^{d+1})^d$. Indeed, after substituting the signature representation of \(Y_t\),
every polynomial in \(Y_t\) becomes a polynomial in the components of
\(\hat{\mathbb B}_t\). By repeated application of the shuffle
product, such a polynomial can be rewritten as a linear functional
of the signature, possibly at a higher truncation level. 
  
    \item \textbf{Parametrized reverse process and its expected signature.} By assumption \ref{it1}, one can now define 
\begin{align*}
    dY_t^\u =\left(-\overline{b}(Y_t^\u) +\operatorname{div}(\overline{a}(Y_t^\u))
    + \langle \u(Y_0^\u), \hat{\mathbb{B}}_t\rangle\right)dt+\sqrt{\overline{a}(Y_t^\u)} dB_t, \quad Y^{\u}_0=X_T.
\end{align*}
 For every fixed \(y_0\) with \(Y_0^\u=y_0\), the process $
    Z:=(Y^\u,\hat{\mathbb{B}}^{\leq N})
$
is a finite-dimensional polynomial process in the sense of Subsection~\ref{sec: poly} whose coefficients depend
parametrically on \(y_0\). As explained there, this implies that its expected signature can be computed analytically, i.e.,
\begin{align*}
        \E[{\mathbb{Z}}_T^{\leq M}|Y_0^\u]=\exp(TL^\u(Y_0^\u)^\top)\e_\emptyset
    \end{align*}
     where $L^{\u}(y)$ denotes the generator matrix of 
    the truncated signature at level $M$  of the process $Z$ starting at $z_0:=(y,\e_\emptyset)$ (Corollary \ref{coro:expectedsigpoly}). In particular, this implies that
    \begin{align*}
        \E[({\mathbb{Y}}_T^\u)^{\leq M}|Y_0^\u]=
        f(Y_0^\u,\u(Y_0^\u)),
    \end{align*}
for an explicit deterministic vector-valued function $f$ given by
    \begin{align}\label{eq:f}
        f(y,\u)=\big(\langle \e_I,\exp(TL^\u(y)^\top)\e_\emptyset\rangle\big)_{I\in \mathcal{C}},
    \end{align}
where $\mathcal{C}$ is the set of all multiindices $I$ such that $\langle \e_I,\mathbb{Z}_T\rangle =\langle \e_I,(\mathbb{Y}_T^\u)^{\leq M}\rangle$. 

Note that the columns of the matrix $L^\u(y)$ are given in \eqref{eq:Lpoly}, with $\b^\u$ and $\a$ the sig-coefficients of the polynomial diffusion $Z$ computed in Lemma \ref{lemma: polycharacteristics} via the coefficients of its characteristics (\eqref{eq: chara poly process})
$(b^{\u},a,0)$, which read
\begin{align*}
 b^{\u}_{i}(z)
   &= -\overline b_i(\pi_Y(z))
      +\sum_{k=1}^d\partial_{x_k}\overline a_{ik}(\pi_Y(z))
      +\langle\u_i,\pi_{\hat{\mathbb B}}(z)\rangle,\\
 b^{\u}_{I}(z)
   &= 1_{\{i_{|I|}=0\}}\langle \e_{I'},\pi_{\hat{\mathbb B}}(z)\rangle
     +\frac12 1_{\{i_{|I|-1}=i_{|I|}\neq 0\}}
       \langle \e_{I''},\pi_{\hat{\mathbb B}}(z)\rangle,\\
 a_{ik}(z)
   &=\overline a_{ik}(\pi_Y(z)),\qquad 
   a_{iI}(z)
   = a_{Ii}(z)= 1_{\{i_{|I|}\neq 0\}}
     \overline\sigma_{i,i_{|I|}}(\pi_Y(z))\langle \e_{I'},\pi_{\hat{\mathbb B}}(z)\rangle,\\
 a_{IJ}(z)
   &= 1_{\{
       i_{|I|}=j_{|J|}\neq 0\}}
      \langle \e_{I'},\pi_{\hat{\mathbb B}}(z)\rangle\langle \e_{J'},\pi_{\hat{\mathbb B}}(z)\rangle,
\end{align*}

where $\pi_Y(z)$ and $\pi_{\hat{\mathbb B}}(z)$ denote the components of $z$ corresponding to $Y$ and $\hat{\mathbb B}^{\leq N}$ in $Z$, respectively, and $\overline\sigma=\sqrt {\overline a}$.
Since by Lemma \ref{lemma: polycharacteristics} the components of $\b^\u$ and $\a$  are equal to some partial derivative of $b^{\u}$ and $a$ evaluated at $z_0=(y,\e_\emptyset)$, respectively, and since $\overline b$ and $\overline \sigma$ are polynomials, we can conclude that the map $y\mapsto L^\u(y)$ has polynomial components and thus
 that the components of $f$ in \eqref{eq:f} are entire in $(y,\u)$.

Alternatively, the assumptions of point \ref{it1} also justify a parametrization of the whole drift directly as a linear combination of the components of $\hat{\mathbb B}_t$, thus yielding
    \begin{align*}
         d\overline Y_t^{\overline\u} =\langle  {\overline\u}(\overline Y_0^{\overline\u}), \hat{\mathbb{B}}_t\rangle dt+\sqrt{\overline{a}(\overline Y_t^{\overline\u})} dB_t, \quad \overline Y_0^{\overline\u}=X_T,
    \end{align*}
 which is an even more tractable representation 
$
        \E[(\overline{\mathbb{Y}}_T^{\overline\u})^{\leq M}|\overline Y_0^{\overline\u}]=\tilde f(\overline Y_0^{\overline\u},\overline\u(\overline Y_0^{\overline\u})),
$
    but  losing the specific structure of the drift of the reverse process.

    \item \label{item: inverse sig match}
    \textbf{Learning via inverse signature matching.} We replace the map $Y_0 \mapsto \u(Y_0)$ in \eqref{eq: sigapprox} with some parametric map $\u_\phi$. Our goal is now to learn $\phi$ 
    by matching  the conditional expectations of $\Y_T^{\u_\phi}$ and $\Y_T$ given $X_T$, i.e.,
    \begin{equation}\label{losstarget}\E[\|\E[({\mathbb{Y}}^{\u_\phi}_T)^{\leq M}|X_T]-\E[{\mathbb{Y}}_T^{\leq M}|X_T]\|^2].
    \end{equation}
    By the previous point, we know that 
$
        \E[({\mathbb{Y}}^{\u_\phi}_T)^{\leq M}|X_T]=f(X_T,\u_\phi(X_T)),
$
    while from Proposition~\ref{prop: inverse sig} below, we get that $\Y_T=\X_T^{-1}$ almost surely  and thus
    $$\E[{\mathbb{Y}}_T^{\leq M}|X_T]=\E[(\mathbb{X}_T^{-1})^{\leq M}|X_T].$$

Here $\mathbb{X}_T^{-1}$ is defined as the unique element of $T((\R^d))$ such that  $\mathbb{X}_T^{-1}\otimes \mathbb{X}_T=\e_\emptyset$, namely
\begin{equation}\label{inverse}
    \left\langle
        \e_{(i_1,\ldots,i_m)},(\X_T)^{-1}
    \right\rangle
    =
    (-1)^m
    \left\langle
        \e_{(i_m,\ldots,i_1)},\X_T
    \right\rangle.
\end{equation}

Observe that
\begin{align*}
    \E[\|f(X_T,\u_\phi(X_T))-\E[(\mathbb{X}_T^{-1})^{\leq M}|X_T]\|^2]
&=\E[\|f(X_T,\u_\phi(X_T))-(\mathbb{X}_T^{-1})^{\leq M}\|^2]\\
&\quad -\E[\|({\mathbb{X}}_T^{-1})^{\leq M}-\E[(\mathbb{X}_T^{-1})^{\leq M}|X_T]\|^2],
\end{align*}
showing that minimizing  \eqref{losstarget} with respect to $\phi$ coincides with minimizing  $$\E[\|f(X_T,\u_\phi(X_T))-(\mathbb{X}_T^{-1})^{\leq M}\|^2].$$  

Letting \(x_0^{(1)},\ldots,x_0^{(K)}\)  denote the available training
instances and moving to the empirical
expectation, we get
\begin{align*}
\min_{\phi}\frac{1}{K}\sum_{m=1}^{K}
&\mathbb E\left[\big\|
f\bigl(X_T,\u_\phi(X_T)\bigr)
-
(\mathbb X_T^{-1})^{\leq M}\big\|^2
\,\middle|\,
X_0=x_0^{(m)}
\right].
\end{align*}
By Jensen's inequality, this condition can be relaxed to
 \begin{align}\label{eq: optprob2}
\min_{\phi}\frac{1}{K}\sum_{m=1}^{K}
\Bigg\|
&\mathbb E\left[
f\bigl(X_T,\u_\phi(X_T)\bigr)
\,\middle|\,
X_0=x_0^{(m)}
\right] 
-
\mathbb E\left[
(\mathbb X_T^{-1})^{\leq M}
\,\middle|\,
X_0=x_0^{(m)}
\right]
\Bigg\|^2,
\end{align}
  which, as explained in \ref{deterministic optimization} below, has more tractable properties.

    When the full expected signature determines the conditional law of the
reverse process \cite{signaturemoments} and if the densities are regular enough, equality
of the reverse path laws identifies the reverse drift and therefore
the score.

    \item \label{deterministic optimization}  \textbf{Analytic expected inverse signature and reduction to a deterministic optimization problem.}
 The second conditional expectation in \eqref{eq: optprob2} can be computed analytically. Indeed, by equation \eqref{inverse}, the (algebraic) inverse signature \((\mathbb X_T^{-1})^{\leq M}\) is a linear transformation of \(\mathbb X_T^{\leq M}\). Since the forward process $X$ in \eqref{eq:SDEforward} is a polynomial diffusion, the conditional expectation of its signature can be computed via the matrix-exponential formula for expected signatures established in Section~\ref{sec: poly} (Corollary \ref{coro:expectedsigpoly}).
 Note that the initial condition $x_0^{(m)}$ enters in the computation of the generator matrix $L$ of $\mathbb X^{\leq M}$ via the sig-coefficients of $X$ computed in Lemma \ref{lemma: polycharacteristics} (equation \eqref{eq:Lpoly}).

The first conditional expectation in \eqref{eq: optprob2} can be approximated by Monte Carlo sampling from \(p_T(\cdot \mid x_0)\). Alternatively, it admits an analytic expression whenever the map \(y \mapsto \mathbf u_\phi(y)\) is entire. Indeed, by the definition of \(f\) in \eqref{eq:f}, the map $(y,\u) \mapsto f(y, \mathbf u)$ has entire components. Hence, the composition \(y \mapsto f(y, \mathbf u_\phi(y))\) has  entire components too. Therefore, provided the assumptions of Theorem 3.20 in \cite{CPS:24holo} are satisfied, the holomorphic formula established therein applies to $X$, which is a polynomial (and hence holomorphic) diffusion.

It follows that, if \(y \mapsto \mathbf u_\phi(y)\) is entire and Theorem 3.20 in \cite{CPS:24holo} applies, both conditional expectations in \eqref{eq: optprob2} admit analytic expressions, and the optimization problem \eqref{eq: optprob2} reduces to a deterministic one.

\item \textbf{Sample generation.}
Once $\phi$ is learned, the map
$Y_0 \mapsto \u_{\phi}(Y_0)$ and in turn \eqref{eq: sigapprox} is fully specified.
New samples can now be generated by drawing
\(Y_0\) from the limiting invariant distribution, evaluating
\(\u_{\phi}(Y_0)\), and numerically solving

\begin{align*}
    dY_t^\u =\left(-\overline{b}(Y_t^\u) +\operatorname{div}(\overline{a}(Y_t^\u))
    + \langle \u_{\phi}(Y_0^\u), \hat{\mathbb{B}}_t\rangle\right)dt+\sqrt{\overline{a}(Y_t^\u)} dB_t.
\end{align*}
The terminal value \(Y_T^\u\) is then taken as a generated approximation
of a sample from \(p_{\mathrm{data}}\).

  \end{enumerate}
  \begin{remark}
 Note that if $X$ is a polynomial process with constant diffusion coefficient, $\overline{a}\equiv \mathrm{const}$, as is the case for standard Ornstein--Uhlenbeck processes, then, recalling \eqref{eq: sigapprox}, the above approach reduces to learning, via inverse signature matching, only the score map, parametrized as a linear function of the signature of the Brownian motion depending parametrically on $Y_0$,
\begin{align*}
  \nabla\log p_{T-t}(Y_t) \approx \langle \mathbf{u}_\phi(Y_0), \hat{\mathbb{B}}_t\rangle.
\end{align*}
New samples are then generated by numerically simulating the reverse-time dynamics
\begin{align*}
  dY_t^{\mathbf{u}} = \Big(-\overline{b}(Y_t^{\mathbf{u}}) + \overline{a}\,\langle \mathbf{u}_\phi(Y_0^{\mathbf{u}}), \hat{\mathbb{B}}_t\rangle\Big)\,dt + \sqrt{\overline{a}}\,dB_t.
\end{align*}

\end{remark}
\begin{remark} \label{rem:X0}
    It is important to notice that the law of $\X_T$ contains information about $X_0$. Indeed, even if at time 0 it holds $\X_0=\e_\emptyset$, we have that the law of $\X_T$ corresponds to the law of the signature of the polynomial process $\mathbb{X}^{(1)}=X-X_0$ whose dynamics satisfy
    $$d\mathbb{X}^{(1)}_t=\overline b(\mathbb{X}^{(1)}_t+ X_0)dt+\sqrt{\overline{a}(\mathbb{X}^{(1)}_t+X_0)} dW_t, \quad X_0=x_0\sim p_{\mathrm{data}},\qquad \mathbb{X}^{(1)}_0=0.$$
    One can see this also by noticing that \eqref{eq:SDEforward}    corresponds to \eqref{entsig} with truncation level 1 and some $\g$ such that $\g^{(1)}=X_0$. 
Indeed, with the notation of \eqref{entsig} and  $\overline{b}$ as in \eqref{eq:SDEforward}, we have
\[
\overline{b}^{\g}( \x^{(1)})=\overline{b}^{\g}(x-X_0)=\overline{b}(x)=\overline{b}(\x^{(1)}+X_0)
\]
and analogously for the diffusion part. Therefore $\overline{b}(\x^{(1)}+X_0)$ and $\overline{a}(\x^{(1)}+X_0)$ enter in the dynamics  of the signature process given in \eqref{eq:Xg}.
\end{remark}

The following proposition is needed to justify point \ref{item: inverse sig match} in the above approach. Recall the definition of $\X_T^{-1}$ given by \eqref{inverse}.
\begin{proposition}\label{prop: inverse sig}
Let \(Y_t:=X_{T-t}\), \(t\in[0,T]\), and denote by
\(\mathbb{X}\) and \(\mathbb{Y}\) the Stratonovich signatures of
\(X\) and \(Y\), respectively. Then
$
    \mathbb{Y}_T=\mathbb{X}_T^{-1}$
   almost surely.
\end{proposition}

\begin{proof}
We proceed by verifying the identity on approximating processes and then moving to the limit.
Let
$
    D_n:=\left\{\frac{kT}{2^n}:k=0,\ldots,2^n\right\},
$
and let \(X^n\) be the piecewise linear interpolation of \(X\) along
\(D_n\). Define
\[
    Y^n_t:=X^n_{T-t},
    \qquad t\in[0,T].
\]
Since \(D_n\) is invariant under the map \(t\mapsto T-t\), \(Y^n\) is
the piecewise linear interpolation of \(Y\) along the same partition
\(D_n\). Let \(\mathbb{X}^n\) and \(\mathbb{Y}^n\) denote the
signatures of \(X^n\) and \(Y^n\), respectively, which are well defined as both of them are semimartingales with respect to their own filtration.

Since \(X^n\) is a continuous path of bounded variation and \(Y^n\)
is the same path run backwards, Theorem~1.25 of \cite{CK:25} gives
\begin{equation}\label{eq:PL-signature-reversal}
    \mathbb{Y}^n_T
    =
    \bigl(\mathbb{X}^n_T\bigr)^{-1}
    \qquad\text{a.s.}
\end{equation}

Fix now \(N\in\Nset\). By the Wong--Zakai approximation for continuous
semimartingales, applied to the equations defining the truncated
signatures, we get that 
$
    \bigl(\mathbb{X}^n_T\bigr)^{\leq N}
$ converges to $
    \mathbb{X}^{\leq N}_T$
and
$
    \bigl(\mathbb{Y}^n_T\bigr)^{\leq N}
$ to $
    \mathbb{Y}^{\leq N}_T$
 in probability \cite[Example~4.20]{CF:19}. Since inversion is continuous on the
finite-dimensional set \(G^N(\Rset^d)\), we also have that
$
   (
        (\mathbb{X}^n_T)^{\leq N}
    )^{-1}
    $ converges to $
    (\mathbb{X}^{\leq N}_T)^{-1}
    $ in probability.
As 
$
    ((\mathbb{X}^n_T)^{-1}
  )^{\leq N}
    =
    ((\mathbb{X}^n_T)^{\leq N}
    )^{-1}
$, by \eqref{eq:PL-signature-reversal} and uniqueness of the
limit in probability,
$
    \mathbb{Y}^{\leq N}_T
    =
    \bigl(\mathbb{X}^{\leq N}_T\bigr)^{-1}
$ almost surely.
Since \(N\) is arbitrary, this identity holds at every tensor level.
Consequently,
$
    \mathbb{Y}_T
    =
    \mathbb{X}_T^{-1}$ almost surely,
where the equality is understood coordinatewise.
\end{proof}

\begin{example}
Note that    in the case of a one-dimensional process $X$, 
we have that each term of the signature of $Y$ at time $T$ corresponds to an element of the form $(Y_T-Y_0)^k/k!$ which, by definition of $Y$, corresponds to
    $(X_0-X_T)^k/k!$. The correspondence between laws then follows directly. 
\end{example}


\begin{thebibliography}{}

\bibitem[Abi~Jaber et~al., 2026]{abi2026stochastic}
Abi~Jaber, E., Attal, E., and Sotnikov, D. (2026).
\newblock {Stochastic control with signatures via Riccati equations on the tensor algebra}.
\newblock {\em Preprint, arXiv:2607.03986}.

\bibitem[Abi~Jaber et~al., 2025]{AbiJaberHainautMotte2025}
Abi~Jaber, E., Hainaut, D., and Motte, E. (2025).
\newblock {Signature Approach for Pricing and Hedging Path-Dependent Options with Frictions}.
\newblock {\em Preprint, arXiv:2511.23295}.

\bibitem[Abi~Jaber and Sotnikov, 2025]{Abi:Fading}
Abi~Jaber, E. and Sotnikov, D. (2025).
\newblock {Exponentially Fading Memory Signature}.
\newblock {\em Preprint, arXiv:2507.03700}.

\bibitem[Arnaudon and Thalmaier, 2010]{ArnaudonThalmaier2010}
Arnaudon, M. and Thalmaier, A. (2010).
\newblock {The Differentiation of Hypoelliptic Diffusion Semigroups}.
\newblock {\em Illinois Journal of Mathematics}, 54(4):1285--1311.

\bibitem[Barancikova et~al., 2025]{salvi}
Barancikova, B., Huang, Z., and Salvi, C. (2025).
\newblock {{SigDiffusions}: Score-Based Diffusion Models for Time Series via Log-Signature Embeddings}.
\newblock In {\em The Thirteenth International Conference on Learning Representations}.

\bibitem[Barilari et~al., 2019]{BarilariBoscainNeel2019}
Barilari, D., Boscain, U., and Neel, R.~W. (2019).
\newblock {Heat Kernel Asymptotics on Sub-Riemannian Manifolds with Symmetries and Applications to the Bi-Heisenberg Group}.
\newblock {\em Annales de la Facult{\'e} des Sciences de Toulouse}, 28(4):707--732.

\bibitem[Baudoin et~al., 2008]{BaudoinHairerTeichmann2008}
Baudoin, F., Hairer, M., and Teichmann, J. (2008).
\newblock Ornstein--uhlenbeck processes on lie groups.
\newblock {\em Journal of Functional Analysis}, 255(4):877--890.

\bibitem[Baule, 2026]{B:26}
Baule, A. (2026).
\newblock {Generative modelling with jump-diffusions}.
\newblock {\em Journal of Statistical Mechanics: Theory and Experiment}, 2026(5):054002.

\bibitem[Bayer et~al., 2023]{BHRS:21}
Bayer, C., Hager, P.~P., Riedel, S., and Schoenmakers, J. (2023).
\newblock {Optimal stopping with signatures}.
\newblock {\em The Annals of Applied Probability}, 33(1):238--273.

\bibitem[Boedihardjo et~al., 2021]{boedihardjo2021expected}
Boedihardjo, H., Diehl, J., Mezzarobba, M., and Ni, H. (2021).
\newblock {The expected signature of Brownian motion stopped on the boundary of a circle has finite radius of convergence}.
\newblock {\em Bulletin of the London Mathematical Society}, 53(1):285--299.

\bibitem[Boedihardjo et~al., 2013]{BoedihardjoPapavasiliouQian2015}
Boedihardjo, H., Papavasiliou, A., and Qian, Z. (2013).
\newblock {Expected Signature of {Gaussian} Processes with Strictly Regular Kernels}.
\newblock {\em Preprint, arXiv:1304.4930}.

\bibitem[Cao et~al., 2026]{cao2026adaptive}
Cao, H., Chen, M., Han, Y., and Xu, R. (2026).
\newblock {Diffusion Models for Adaptive Sequential Data Generation}.
\newblock {\em Preprint, arXiv:2606.06007}.

\bibitem[Cass and Ferrucci, 2024]{CassFerrucci2024}
Cass, T. and Ferrucci, E. (2024).
\newblock {On the {Wiener} Chaos Expansion of the Signature of a {Gaussian} Process}.
\newblock {\em Probability Theory and Related Fields}, 189:909--947.

\bibitem[Cass and Salvi, 2024]{CC:24}
Cass, T. and Salvi, C. (2024).
\newblock {Lecture notes on rough paths and applications to machine learning}.
\newblock {\em Preprint, arXiv:2404.06583}.

\bibitem[Chevyrev and Friz, 2019]{CF:19}
Chevyrev, I. and Friz, P.~K. (2019).
\newblock {Canonical RDEs and general semimartingales as rough paths}.
\newblock {\em The Annals of Probability}, 47(1):420--463.

\bibitem[Chevyrev and Kormilitzin, 2026]{CK:25}
Chevyrev, I. and Kormilitzin, A. (2026).
\newblock {A primer on the signature method in machine learning}.
\newblock In {\em Signature Methods in Finance: An Introduction with Computational Applications}, pages 3--64. Springer.

\bibitem[Chevyrev and Lyons, 2016]{ChevyrevLyons2016}
Chevyrev, I. and Lyons, T. (2016).
\newblock {Characteristic Functions of Measures on Geometric Rough Paths}.
\newblock {\em The Annals of Probability}, 44(6):4049--4082.

\bibitem[Chevyrev and Oberhauser, 2022]{signaturemoments}
Chevyrev, I. and Oberhauser, H. (2022).
\newblock {Signature Moments to Characterize Laws of Stochastic Processes}.
\newblock {\em Journal of Machine Learning Research}, 23(176):1--42.

\bibitem[Conforti and L{\'e}onard, 2022]{CL:22}
Conforti, G. and L{\'e}onard, C. (2022).
\newblock {Time reversal of {Markov} processes with jumps under a finite entropy condition}.
\newblock {\em Stochastic Processes and their Applications}, 144:85--124.

\bibitem[Cuchiero et~al., 2025a]{CGMS:23}
Cuchiero, C., Gazzani, G., M{\"o}ller, J., and Svaluto-Ferro, S. (2025a).
\newblock {Joint calibration to SPX and VIX options with signature-based models}.
\newblock {\em Mathematical Finance}, 35(1):161--213.

\bibitem[Cuchiero et~al., 2023]{CGS:22}
Cuchiero, C., Gazzani, G., and Svaluto-Ferro, S. (2023).
\newblock {Signature-Based Models: Theory and Calibration}.
\newblock {\em SIAM Journal on Financial Mathematics}, 14(3):910--957.

\bibitem[Cuchiero et~al., 2012]{CKT:12}
Cuchiero, C., Keller-Ressel, M., and Teichmann, J. (2012).
\newblock {Polynomial processes and their applications to mathematical finance}.
\newblock {\em Finance and Stochastics}, 16(4):711--740.

\bibitem[Cuchiero and M{\"o}ller, 2025]{cuchiero2025signature}
Cuchiero, C. and M{\"o}ller, J. (2025).
\newblock {Signature methods in stochastic portfolio theory}.
\newblock {\em SIAM Journal on Financial Mathematics}, 16(4):1239--1303.

\bibitem[Cuchiero and Pachschw{\"o}ll, 2026]{CP:26}
Cuchiero, C. and Pachschw{\"o}ll, J. (2026).
\newblock {Polynomial Signature {McKean-Vlasov} {SDEs}}.
\newblock {\em Working Paper}.

\bibitem[Cuchiero et~al., 2025b]{CPS:UATlevy}
Cuchiero, C., Primavera, F., and Svaluto-Ferro, S. (2025b).
\newblock {Universal approximation theorems for continuous functions of c{\`a}dl{\`a}g paths and {L{\'e}vy}-type signature models}.
\newblock {\em Finance and Stochastics}, 29(2):289--342.

\bibitem[Cuchiero et~al., 2026a]{CPS:24holo}
Cuchiero, C., Primavera, F., and Svaluto-Ferro, S. (2026a).
\newblock {Holomorphic jump-diffusions}.
\newblock {\em Stochastic Processes and their Applications}, 191:104781.

\bibitem[Cuchiero et~al., 2026b]{CST:23}
Cuchiero, C., Svaluto-Ferro, S., and Teichmann, J. (2026b).
\newblock {Signature {SDEs} from an affine and polynomial perspective}.
\newblock {\em Electronic Journal of Probability}, 31:1--57.

\bibitem[De~Marco et~al., 2026]{demarco2026schrodinger}
De~Marco, S., Pham, H., and Zanni, D. (2026).
\newblock Schr{\"o}dinger bridges with jumps for time series generation.
\newblock {\em Preprint, arXiv:2602.20011}.

\bibitem[Doob, 1948]{Doob1948}
Doob, J.~L. (1948).
\newblock {Asymptotic Properties of Markoff Transition Probabilities}.
\newblock {\em Transactions of the American Mathematical Society}, 63(3):393--421.

\bibitem[Ethier and Kurtz, 1986]{EK}
Ethier, S.~N. and Kurtz, T.~G. (1986).
\newblock {\em {Markov processes: characterization and convergence}}.
\newblock Wiley Series in Probability and Mathematical Statistics: Probability and Mathematical Statistics. John Wiley \& Sons Inc.

\bibitem[Fawcett, 2003]{Fawcett2003}
Fawcett, T. (2003).
\newblock {\em {Problems in Stochastic Analysis: Connections between Rough Paths and Non-Commutative Harmonic Analysis}}.
\newblock PhD thesis, University of Oxford.

\bibitem[Fermanian et~al., 2024]{FermanianChangLyonsBiau2023}
Fermanian, A., Chang, J., Lyons, T., and Biau, G. (2024).
\newblock The insertion method to invert the signature of a path.
\newblock In {\em Recent Advances in Econometrics and Statistics: Festschrift in Honour of Marc Hallin}, pages 575--595. Springer.

\bibitem[Filipovi{\'c} and Larsson, 2020]{FL:20}
Filipovi{\'c}, D. and Larsson, M. (2020).
\newblock {Polynomial Jump-Diffusion Models}.
\newblock {\em Stochastic Systems}, 10(1):71--97.

\bibitem[Friz et~al., 2026]{friz2024expected}
Friz, P.~K., Hager, P.~P., and Tapia, N. (2026).
\newblock {On Expected Signatures and Signature Cumulants in Semimartingale Models}.
\newblock In {\em Signature Methods in Finance: An Introduction with Computational Applications}, pages 381--424. Springer.

\bibitem[Friz and Shekhar, 2017]{FrizShekhar2017}
Friz, P.~K. and Shekhar, A. (2017).
\newblock {General Rough Integration, L{\'e}vy Rough Paths and a L{\'e}vy--Kintchine-Type Formula}.
\newblock {\em The Annals of Probability}, 45(4):2707--2765.

\bibitem[Gasteratos et~al., 2025]{gasteratos2025novelty}
Gasteratos, I., Jacquier, A., Lemercier, M., Lyons, T., and Salvi, C. (2025).
\newblock {Novelty detection on path space}.
\newblock {\em Preprint, arXiv:2512.03243}.

\bibitem[Ge et~al., 2026]{sigdegdistillationmakingdiffusion}
Ge, W., Jiang, L., Yang, L., Roveri, L., Yi, M., Chen, P.-Y., Cariou-Kotlarek, N., Buet-Golfouse, F., Mittal, G., and Ni, H. (2026).
\newblock {Sig-DEG for Distillation: Making Diffusion Models Faster and Lighter}.
\newblock {\em Preprint, arXiv:2508.16939}.

\bibitem[Gennaro et~al., 2026]{GMP:26}
Gennaro, A., Mastrolia, T., and Primavera, F. (2026).
\newblock {Signature Methods for Optimal Market Making}.
\newblock {\em Preprint, arXiv:2606.19772}.

\bibitem[Gerlach and Nittka, 2012]{GerlachNittka2012}
Gerlach, M. and Nittka, R. (2012).
\newblock {A New Proof of Doob's Theorem}.
\newblock {\em Journal of Mathematical Analysis and Applications}, 388(2):763--774.

\bibitem[Guo et~al., 2026]{guo2026conditionaldiffusionguidancehard}
Guo, Z., Tang, W., and Xu, R. (2026).
\newblock {Conditional Diffusion Guidance under Hard Constraint: A Stochastic Analysis Approach}.
\newblock {\em Preprint, arXiv:2602.05533}.

\bibitem[Hamdouche et~al., 2026]{hamdouche2026nonparametric}
Hamdouche, M., Henry-Labord{\`e}re, P., and Pham, H. (2026).
\newblock Nonparametric generative modeling for time series via {S}chr{\"o}dinger bridge.
\newblock {\em Journal of Machine Learning Research}, 27(58):1--23.

\bibitem[Han et~al., 2025]{han2025stochastic}
Han, Y., Razaviyayn, M., and Xu, R. (2025).
\newblock Stochastic control for fine-tuning diffusion models: Optimality, regularity, and convergence.
\newblock In {\em Proceedings of the 42nd International Conference on Machine Learning}, volume 267 of {\em Proceedings of Machine Learning Research}, pages 21844--21870.

\bibitem[Haussmann and Pardoux, 1986]{HS:86}
Haussmann, U.~G. and Pardoux, E. (1986).
\newblock {Time reversal of diffusions}.
\newblock {\em The Annals of Probability}, 14(4):1188--1205.

\bibitem[Ho et~al., 2020]{ho2020denoising}
Ho, J., Jain, A., and Abbeel, P. (2020).
\newblock Denoising diffusion probabilistic models.
\newblock In {\em Advances in Neural Information Processing Systems}, volume~33, pages 6840--6851.

\bibitem[Hyv{\"a}rinen, 2005]{hyvarinen2005estimation}
Hyv{\"a}rinen, A. (2005).
\newblock {Estimation of Non-Normalized Statistical Models by Score Matching}.
\newblock {\em Journal of Machine Learning Research}, 6:695--709.

\bibitem[Jacod and Shiryaev, 1987]{JS:87}
Jacod, J. and Shiryaev, A. (1987).
\newblock {\em {Limit Theorems for Stochastic Processes}}.
\newblock Springer Berlin, Heidelberg.

\bibitem[Jansen and Kurt, 2014]{JK:14}
Jansen, S. and Kurt, N. (2014).
\newblock {On the notion(s) of duality for Markov processes}.
\newblock {\em Probability Surveys}, 11:59--120.

\bibitem[Kalsi et~al., 2020]{KLP:20}
Kalsi, J., Lyons, T., and Perez~Arribas, I. (2020).
\newblock {Optimal execution with rough path signatures}.
\newblock {\em SIAM Journal on Financial Mathematics}, 11(2):470--493.

\bibitem[Kurtz et~al., 1995]{KP:95}
Kurtz, T.~G., Pardoux, E., and Protter, P. (1995).
\newblock {Stratonovich stochastic differential equations driven by general semimartingales}.
\newblock {\em Annales de l'I.H.P. Probabilit\'es et statistiques}, 31(2):351--377.

\bibitem[Le~Donne, 2025]{LeDonne2025}
Le~Donne, E. (2025).
\newblock {Nilpotent Lie Groups}.
\newblock In {\em Metric Lie Groups}, volume 306 of {\em Graduate Texts in Mathematics}, pages 243--295. Springer.

\bibitem[Levin et~al., 2013]{LLN:13}
Levin, D., Lyons, T., and Ni, H. (2013).
\newblock {Learning from the past, predicting the statistics for the future, learning an evolving system}.
\newblock {\em Preprint, arXiv:1309.0260}.

\bibitem[Liao et~al., 2024]{NSWLB:20}
Liao, S., Ni, H., Sabat{\'e}-Vidales, M., Szpruch, L., Wiese, M., and Xiao, B. (2024).
\newblock {Sig‐Wasserstein GANs for conditional time series generation}.
\newblock {\em Mathematical Finance}, 34(2):622--670.

\bibitem[Lou et~al., 2023]{lou2023pcfgan}
Lou, H., Li, S., and Ni, H. (2023).
\newblock {PCF-GAN: generating sequential data via the characteristic function of measures on the path space}.
\newblock In {\em Advances in Neural Information Processing Systems}, volume~36, pages 39755--39781.

\bibitem[Lyons, 1998]{L:98}
Lyons, T. (1998).
\newblock {Differential equations driven by rough signals}.
\newblock {\em Revista Matem{\'a}tica Iberoamericana}, 14(2):215--310.

\bibitem[Lyons et~al., 2019]{LyonsNejadPerezArribas2020}
Lyons, T., Nejad, S., and Perez~Arribas, I. (2019).
\newblock {Numerical Method for Model-Free Pricing of Exotic Derivatives in Discrete Time Using Rough Path Signatures}.
\newblock {\em Applied Mathematical Finance}, 26(6):583--597.

\bibitem[Lyons and Ni, 2015]{LyonsNi2015}
Lyons, T. and Ni, H. (2015).
\newblock {Expected signature of Brownian motion up to the first exit time from a bounded domain}.
\newblock {\em Annals of Probability}, 43(5):2729--2762.

\bibitem[Lyons and Victoir, 2004]{LyonsVictoir2004}
Lyons, T. and Victoir, N. (2004).
\newblock {Cubature on {Wiener} Space}.
\newblock {\em Proceedings of the Royal Society of London. Series A: Mathematical, Physical and Engineering Sciences}, 460(2041):169--198.

\bibitem[Lyons and Xu, 2017]{LyonsXu2017}
Lyons, T. and Xu, W. (2017).
\newblock Hyperbolic development and inversion of signature.
\newblock {\em Journal of Functional Analysis}.

\bibitem[Lyons and Xu, 2018]{LyonsXu2018}
Lyons, T. and Xu, W. (2018).
\newblock Inverting the signature of a path.
\newblock {\em Journal of the European Mathematical Society}, 20(7):1655--1687.

\bibitem[Ni, 2012]{Ni2012}
Ni, H. (2012).
\newblock {\em {The Expected Signature of a Stochastic Process}}.
\newblock PhD thesis, University of Oxford.

\bibitem[Papavasiliou and Ladroue, 2011]{PapavasiliouLadroue2011}
Papavasiliou, A. and Ladroue, C. (2011).
\newblock {Parameter Estimation for Rough Differential Equations}.
\newblock {\em The Annals of Statistics}, 39(4):2047--2073.

\bibitem[Perez~Arribas et~al., 2021]{ASS:20}
Perez~Arribas, I., Salvi, C., and Szpruch, L. (2021).
\newblock {Sig-SDEs model for quantitative finance}.
\newblock In {\em Proceedings of the First ACM International Conference on AI in Finance}, ICAIF '20, New York, NY, USA.

\bibitem[Popov et~al., 2025]{popov2025improved}
Popov, V., Yermekova, A., Sadekova, T., Khrapov, A., and Kudinov, M.~S. (2025).
\newblock {Improved Sampling Algorithms for {L{\'e}vy}-{It{\^o}} Diffusion Models}.
\newblock In {\em The Thirteenth International Conference on Learning Representations}.

\bibitem[Semnani et~al., 2025]{semnani2025path}
Semnani, P., Guan, V., Robeva, E., and Lee, D. (2025).
\newblock {Path-Dependent {SDEs}: Solutions and Parameter Estimation}.
\newblock {\em Preprint, arXiv:2505.22646}.

\bibitem[Sohl-Dickstein et~al., 2015]{sohl2015deep}
Sohl-Dickstein, J., Weiss, E., Maheswaranathan, N., and Ganguli, S. (2015).
\newblock Deep unsupervised learning using nonequilibrium thermodynamics.
\newblock In {\em Proceedings of the 32nd International Conference on Machine Learning}, volume~37 of {\em Proceedings of Machine Learning Research}, pages 2256--2265.

\bibitem[Song and Ermon, 2019]{songermon2019}
Song, Y. and Ermon, S. (2019).
\newblock Generative modeling by estimating gradients of the data distribution.
\newblock In {\em Advances in Neural Information Processing Systems}, volume~32.

\bibitem[Song et~al., 2021]{song2021score}
Song, Y., Sohl-Dickstein, J., Kingma, D.~P., Kumar, A., Ermon, S., and Poole, B. (2021).
\newblock {Score-Based Generative Modeling through Stochastic Differential Equations}.
\newblock In {\em International Conference on Learning Representations (ICLR)}.

\bibitem[Strichartz, 1986]{Strichartz1986}
Strichartz, R.~S. (1986).
\newblock {Sub-Riemannian Geometry}.
\newblock {\em Journal of Differential Geometry}, 24(2):221--263.

\bibitem[Tang and Zhao, 2025]{Tang2025Scorebased}
Tang, W. and Zhao, H. (2025).
\newblock Score-based diffusion models via stochastic differential equations.
\newblock {\em Statistical Surveys}, 19:28--64.

\bibitem[Vincent, 2011]{vincent2011connection}
Vincent, P. (2011).
\newblock {A Connection Between Score Matching and Denoising Autoencoders}.
\newblock {\em Neural Computation}, 23(7):1661--1674.

\bibitem[Yoon et~al., 2023]{yoon2023score}
Yoon, E., Park, K., Kim, S., and Lim, S. (2023).
\newblock {Score-based Generative Models with {L{\'e}vy} Processes}.
\newblock In {\em Advances in Neural Information Processing Systems}, volume~36.

\end{thebibliography}
\end{document}